\documentclass[oneside]{book}

\usepackage{graphicx}
\usepackage{color}
\usepackage{hyperref}
\usepackage[text={5.5in,8.375in},top=1.25in]{geometry}
\usepackage[font=small]{caption}

\usepackage{titlesec}
\def\thesection{\arabic{section}}
\def\thesubsection{\arabic{section}.\arabic{subsection}}
\titleformat{\section}{\bfseries\sffamily\Large}{\thesection.}{5pt}{}
\titleformat{\subsection}{\bfseries\sffamily\large}{\thesubsection}{5pt}{}

\makeatletter
\renewenvironment{thebibliography}[1]
     {\section*{\bibname}%
      \smallskip
      \list{\@biblabel{\@arabic\c@enumiv}}%
           {\settowidth\labelwidth{\@biblabel{#1}}%
            \leftmargin\labelwidth
            \advance\leftmargin\labelsep
            \@openbib@code
            \usecounter{enumiv}%
            \let\p@enumiv\@empty
            \renewcommand\theenumiv{\@arabic\c@enumiv}}%
      \sloppy
      \clubpenalty4000
      \@clubpenalty \clubpenalty
      \widowpenalty4000%
      \sfcode`\.\@m}
     {\def\@noitemerr
       {\@latex@warning{Empty `thebibliography' environment}}%
      \endlist}
\makeatother

\usepackage[bb=boondox, bbscaled=1.06]{mathalpha}
\DeclareMathAlphabet{\mathCR}{U}{dsrom}{m}{n}    

\usepackage{listings}
\definecolor{codebgcolor}{rgb}{.87,.94,.98}
\lstdefinestyle{pythonstyle}{                          
   language=Python,
   basicstyle=\ttfamily\footnotesize,
   commentstyle=\ttfamily\footnotesize\slshape,
   backgroundcolor=\color{codebgcolor},
   frame=l,
   tabsize=4,
   breaklines=true,
   breakatwhitespace=true,
   captionpos=b
}
\def\dop{{\rm d}}
\def\eop{{\rm e}}
\def\Bb{{\bf b}}
\def\Be{{\bf e}}
\def\Bf{{\bf f}}
\def\Bu{{\bf u}}
\def\Bv{{\bf v}}
\def\Bx{{\bf x}}
\def\By{{\bf y}}
\def\Bz{{\bf z}}
\def\wh#1{\widehat{#1}}
\def\Bzero{{\bf 0}}
\def\BA{{\bf A}}
\def\BAs{\BA{\kern-1pt}}
\def\BU{{\bf U}}
\def\BV{{\bf V}}
\def\BSigma{\mbox{\boldmath$\Sigma$}}
\def\BI{{\bf I}}
\def\R{\mathCR{R}}

\def\Rmn{\R^{m\times n}}
\def\Rnn{\R^{n\times n}}
\def\eps{\varepsilon}
\def\bnoise{\Bb_{\rm noise}}
\def\frec{\Bf_{\rm rec}}
\def\flam{\Bf_\lambda}

\def\DIVIDER{\medskip\centerline{\rule{120pt}{0.5pt}}\medskip}

\makeatletter
\newcounter{prob}
\newcounter{quest}

\newenvironment{exercises}
          {\section*{Exercises}
           \begin{list}{E\arabic{section}.\arabic{prob}}{\usecounter{prob}}
          }
          {\end{list}}
\newenvironment{reflection}
          {\subsection*{Question for Reflection}
           \begin{list}{R\arabic{section}.\arabic{quest}}{\usecounter{quest}}
          }
          {\end{list}}
\newenvironment{reflections}
          {\subsection*{Questions for Reflection}
           \begin{list}{R\arabic{section}.\arabic{quest}}{\usecounter{quest}}
          }
          {\end{list}}

\makeatother

\mathcode`\@="8000
{\catcode`\@=\active \gdef@{\mkern1mu}}

\begin{document}
\thispagestyle{empty}
\vspace*{3em}

\begin{center}
{\LARGE \sf Regularizing Ill-Posed Inverse Problems:\\[6pt]
   Deblurring Barcodes}

\vspace*{2em}
{\normalsize Mark Embree \\[3pt] Department of Mathematics, Virginia Tech \\[3pt] \normalsize  embree@vt.edu}
\end{center}

\noindent
\begin{quote}
Summary: Construct a mathematical model that describes how an image 
gets blurred.
Convert a calculus problem into a linear algebra problem by discretization.
Inverting the blurring process should sharpen up an image; this requires the
solution of a system of linear algebraic equations.
Solving this linear system of equations turns out to be delicate, 
as deblurring is an example of an \emph{ill-posed inverse problem}.
To address this challenge, we recast the system as a regularized 
least squares problem (also known as ridge regression).

\bigskip
\noindent
Prerequisites:  Multivariable Calculus (integration, gradients, optimization)\\
\phantom{Prerequisites: }Linear Algebra (matrix-vector products, solving linear systems)
\end{quote}

\bigskip
An image is a two-dimensional projection of a three-dimensional reality.
Whether that image comes from a camera, a microscope, or a telescope, 
the optics and environment will impose some distortions,
\emph{blurring} the image.  
Can we improve the image to remove this blur?
When editing your photos, you might simply apply a ``sharpen'' tool.
How might such an operation work?

This manuscript describes a mathematical model for blurring.  
If we know how to blur an image, then we should be able to sharpen 
(or \emph{deblur}) the image by ``inverting'' the blurring operation,
in essence, by doing the blurring \emph{in reverse}.
We will discover that blurring operations are notoriously delicate 
to undo.
Indeed, they provide a great example of what computational scientists
call \emph{ill-posed inverse problems}.
(The reason? Two quite distinct images might become very similar 
when they are blurred:
small changes to a blurred image could correspond 
to big changes in the original image.)
To handle this delicate inversion, we will employ a technique called
\emph{regularization}.  
By dialing in the right amount of regularization, we will be able to
deblur images much more reliably.
To test out this technique, we will deblur some Universal Product Code (UPC) 
barcodes (which can be regarded as one-dimensional bitmap images).

The application of inverse problems to barcode deblurring comes 
from the book \emph{Discrete Inverse Problems: Insight and Algorithms}
by Per Christian Hansen~\cite[sect.~7.1]{Han10}, which has been
an inspiration for this manuscript.
Hansen's book is a wonderful resource for students who would like to 
dig deeper into the subject of regularization.
We also point the interested reader to related work by 
Iwen, Santosa, and Ward~\cite{ISW13} and Santosa and Goh~\cite{SG22}
on the specific problem of decoding UPC symbols, the subject of 
Section~\ref{ME:sec:UPC}.
The modeling exercises described here have grown out of a 
homework assignment for the CMDA~3606 (Mathematical Modeling: Methods and Tools)
course at Virginia Tech.  For additional material from this course, 
see~\cite{Emb3606}.

\section{A model for blurring}

When an image is blurred, we might intuitively say that it is ``smeared out'':
the true value of the image at a point is altered in some way that depends upon
the immediately surrounding part of the image.  
For example, when a light pixel is surrounded by dark pixels,
that light pixel should get darker when the image is blurred.  
How might you develop a mathematical description of this intuitive process?

In this manuscript we focus on simple one-dimensional ``images,''
but the techniques we develop here can be readily applied to 
photographs and other two-dimensional images, such as the kinds
that occur in medical imaging.

\begin{figure}[h!]
\begin{center}
  \includegraphics[height=1.25in]{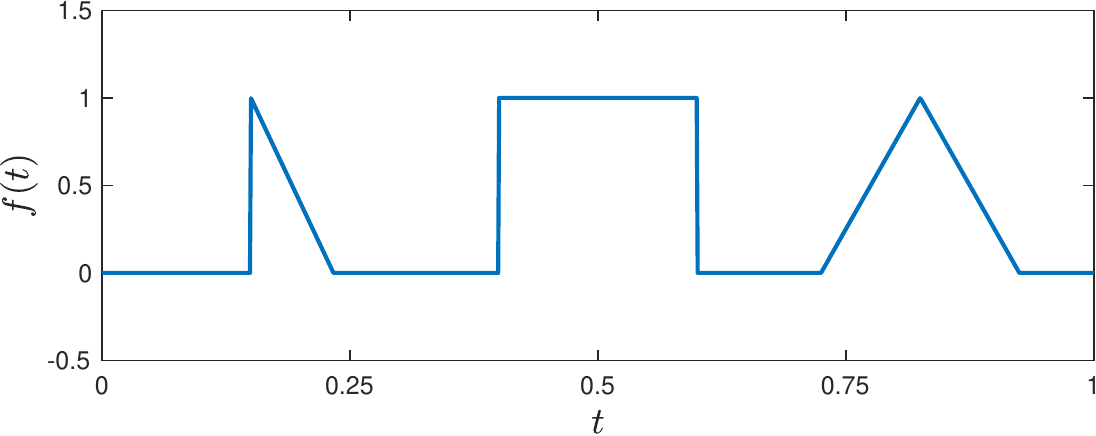}
\end{center}

\vspace*{-1.5em}
\caption{\label{ME:fig:f1}
An example of a one-dimensional ``image'' $f(t)$.  
The vertical lines indicate jump discontinuities in $f$.
}
\end{figure}

When we speak of a one-dimensional ``image,'' we will think of a 
simple function $f(t)$ defined over the interval $0 \le t\le 1$.
Figure~\ref{ME:fig:f1} shows a sample image $f(t)$,
which has a few jump discontinuities to give sharp ``edges.''

\begin{figure}
\begin{center}
  \includegraphics[height=1.25in]{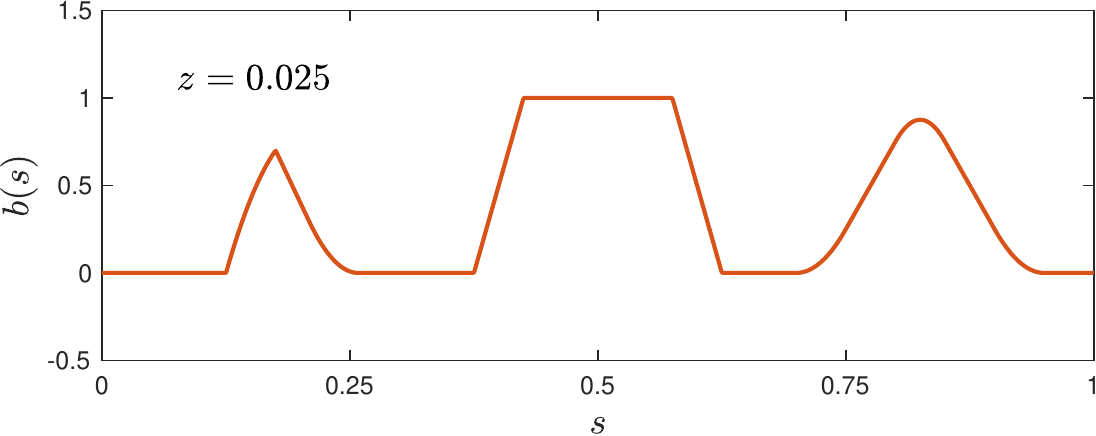}
 
 \includegraphics[height=1.25in]{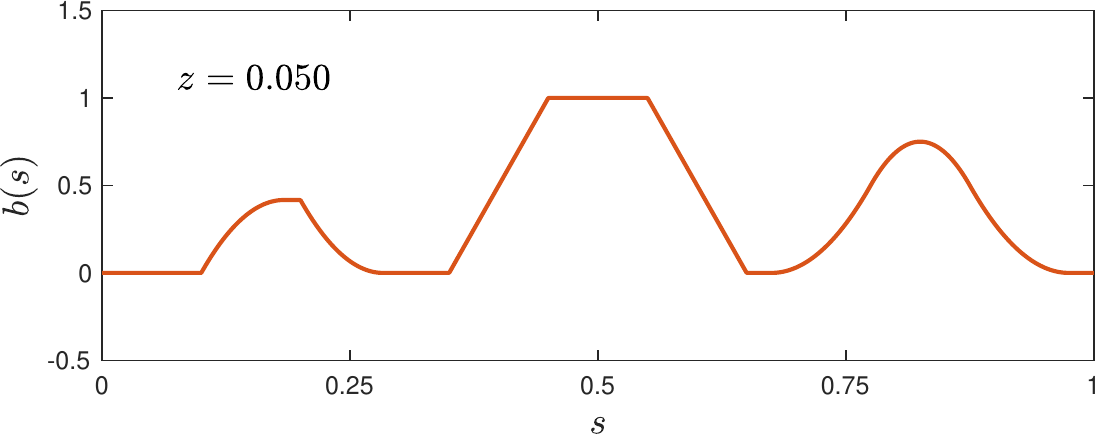}
\end{center}

\vspace*{-1.5em}
\caption{\label{ME:fig:b1}
A blurred version of the image from Figure~\ref{ME:fig:f1},
using the simple blurring function~(\ref{ME:avgker}) 
with parameter $z=0.025$ (top) and $z=0.05$ (bottom):
the larger the value of $z$, the stronger the blurring effect.
}
\end{figure}

To blur this image, we will follow a simple idea:

\smallskip
\centerline{replace the value $f(t)$ with a \emph{local average}.}  
\smallskip
\noindent Since $f$ is a function of $t\in[0,1]$, 
its \emph{local average} will be an integral over a small interval.
Let us focus on some specific point $s\in[0,1]$, and suppose
we want to blur $f(s)$ by averaging $f$ over the region $t\in [s-z, s+z]$ 
for some small value of $z>0$.  
Integrate $f$ over this interval, and divide by the width of the interval
to get our first model for blurring a function:
\begin{equation} \label{ME:bavg}
 b(s) = \frac{\mbox{area under $f(t)$ for $t\in[s-z,s+z]$}}{\mbox{width of interval $[s-z,s+z]$}} = \frac{1}{2z} \int_{s-z}^{s+z} f(t)\,\dop t.
\end{equation}
Think for a moment about the influence of the parameter $z$.  
How does $z$ affect the blurring?  
As $z$ increases, we average $f$ over a wider window; 
more distant points influence $b(s)$, resulting in a more severe blur.  
Figure~\ref{ME:fig:b1} confirms this intuition, showing the affect of
averaging the $f$ from Figure~\ref{ME:fig:f1} over intervals of width $2z$ 
for $z=0.025$ and $z=0.05$.

\begin{figure}[t!]
\begin{center}
  \includegraphics[height=1.2in]{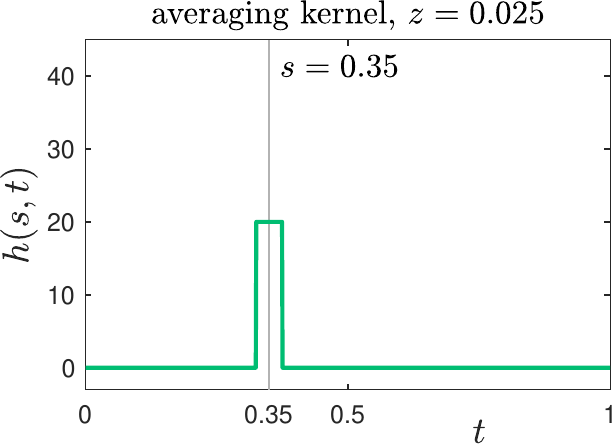}\quad
  \includegraphics[height=1.2in]{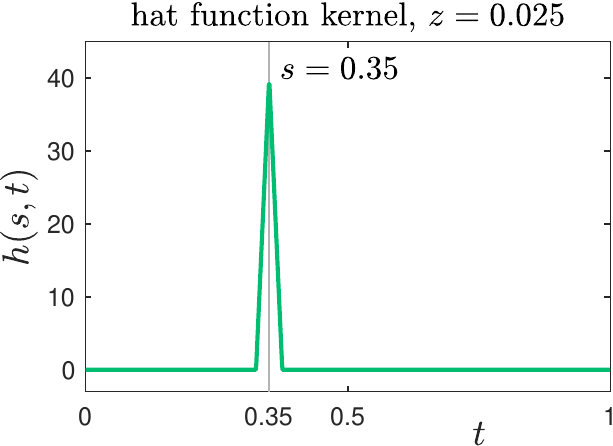}\quad
  \includegraphics[height=1.2in]{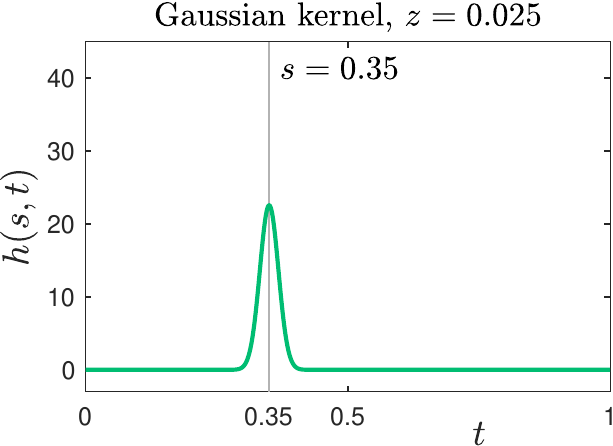}
 
\medskip 
  \includegraphics[height=1.2in]{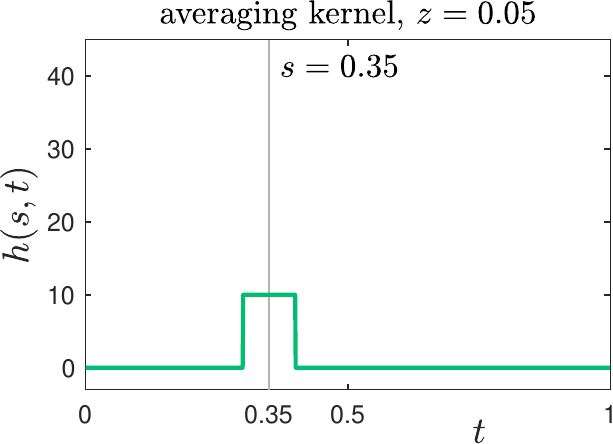}\quad
  \includegraphics[height=1.2in]{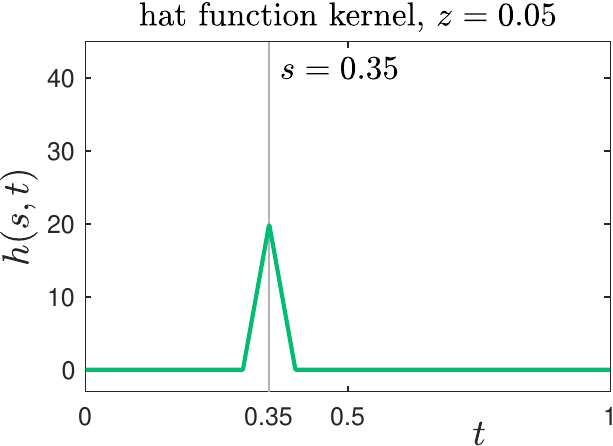}\quad
  \includegraphics[height=1.2in]{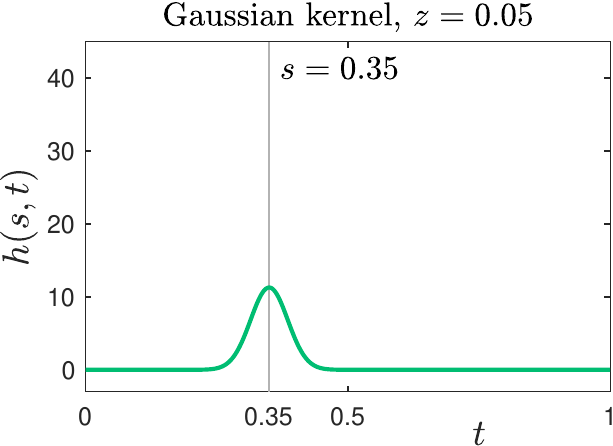}
\end{center}

\vspace*{-1.5em}
\caption{\label{ME:fig:kernels}
Three choices of the kernel function $h(s,t)$, centered at $s=0.35$.
The top plots use the blurring parameter $z=0.025$; 
the bottom plots use $z=0.05$.
As $z$ gets increases, the kernel spreads more
around the point $s$, resulting in more blurring.
}
\end{figure}

We can use $z$ to adjust the severity of the blur, 
but still this model of blurring is quite crude.
To facilitate more general models of blurring, 
realize that the integral in~(\ref{ME:bavg}) can be written as
\begin{equation} \label{ME:conv}
    b(s) = \int_0^1 h(s,t)@f(t)\,\dop t,
\end{equation}
where $h(s,t)$ is a function that controls how the blurring is computed.
Such an $h(s,t)$ is called a \emph{kernel function}.
For the example in~(\ref{ME:bavg}), we have the \emph{averaging kernel}
\begin{equation} \label{ME:avgker}
\begin{array}{ll}
\displaystyle{
 \kern-3.5pt
 h(s,t) = \left\{\begin{array}{rl}
                     \displaystyle{\frac{1}{2z},} & |t-s|\le z; \\[9pt]
                       0,  & |t-s|>z.
            \end{array}\right. }
&
\mbox{\begin{picture}(100,0)(0,5)
\put(24,0){\includegraphics[width=1in]{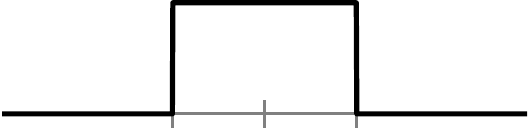}}
\put(39,-4.5){\scriptsize{$s\kern-1.5pt-\kern-1.5pt z$}}
\put(58.2,-4.5){\scriptsize{$s$}}
\put(67.25,-4.5){\scriptsize{$s\kern-1.5pt+\kern-1.5pt z$}}
\end{picture}}
\end{array}
\end{equation}
Consider how different choices of $h(s,t)$ can lead to different kinds of 
blurring.
For example, in our initial model, $b(s)$ is influenced equally by $f(t)$ 
for all $t\in[s-z,s+z]$.  Instead, we might expect $b(s)$ to be more
heavily influenced by $f(t)$ for $t\approx s$, and less influenced by $f(t)$
when $t$ is far from $s$.  We could accomplish this goal by using the 
\emph{hat-function kernel}
\begin{equation} \label{ME:hatker}
\begin{array}{ll}
\displaystyle{
 \kern0pt
 h(s,t) = \frac{1}{z} \max\Big(0, 1 - \frac{|t-s|}{z}\Big)}
&
\mbox{\begin{picture}(100,0)(0,14)
\put(20.3,0){\includegraphics[width=1in]{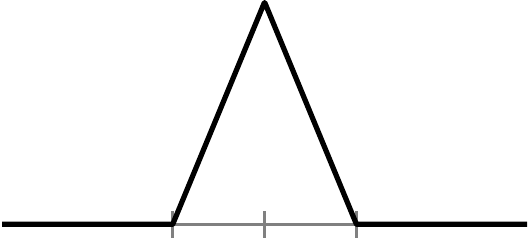}}
\put(35.5,-4.5){\scriptsize{$s\kern-1.5pt-\kern-1.5pt z$}}
\put(54.5,-4.5){\scriptsize{$s$}}
\put(63,-4.5){\scriptsize{$s\kern-1.5pt+\kern-1.5pt z$}}
\end{picture}}
\end{array}
\end{equation}
or the \emph{Gaussian kernel}
\begin{equation} \label{ME:gaussker}
\begin{array}{ll}
\displaystyle{
 \kern-24pt
 h(s,t) = \frac{1}{\sqrt{\pi}@@z} \eop^{-(t-s)^2/z^2}.}
&
\mbox{\begin{picture}(100,0)(0,7)
\put(44.5,0){\includegraphics[width=1in]{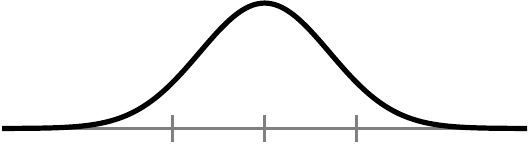}}
\put(59.5,-4.5){\scriptsize{$s\kern-1.5pt-\kern-1.5pt z$}}
\put(78.5,-4.5){\scriptsize{$s$}}
\put(87,-4.5){\scriptsize{$s\kern-1.5pt+\kern-1.5pt z$}}
\end{picture}}
\end{array}
\end{equation}
Figure~\ref{ME:fig:kernels} compares these three choices of $h(s,t)$ 
for $s=0.35$ with the two blurring parameters $z=0.05$ and $z=0.025$.
(The larger $z$, the stronger the blur.)
In a real experiment, we could potentially determine the best choice 
of kernel $h(s,t)$ and blurring factor $z$ by doing experiments with 
our camera, taking pictures of test patterns and measuring the resulting blur.

As the kernel $h(s,t)$ gets more sophisticated, so too does the 
calculus required to convert a signal $f(t)$ into its blurred version.
Even more difficult is the task of taking a blurry image $b(s)$ and 
\emph{inverting} the integration process to determine the unblurred image $f(t)$.
To make such tasks computationally tractable, we will 
approximate the calculus problem with a simpler problem 
involving linear algebra.

\begin{reflections}
\item Can you think of any other good shapes for a blurring function?
\smallskip
\item All the blurring functions described here are \emph{symmetric} about 
the point $s$, which means that $h(s,s-\tau) = h(s,s+\tau)$ for all $\tau >0$.  
Can you design a blurring function $h(s,t)$ so that $b(s)$ in~(\ref{ME:bavg})
only depends on values of $f(t)$ for $t\le s$?
\end{reflections}

\section{Turning calculus into linear algebra}

The trick for turning calculus into linear algebra is a technique called 
\emph{discretization}.  
(This practical approximation aligns very nicely with applications:  
typically we only know our ``image'' $f$ at discrete pixel locations, 
not as a continuous function $f(t)$ for all values of $t\in[0,1]$.)
We shall convert the continuous real intervals 
$s\in [0,1]$ and $t\in[0,1]$ into sets of $n$ discrete points:
\[ s_j = \frac{j - 1/2}{n}, 
    \qquad 
   t_k = \frac{k - 1/2}{n}, \rlap{\qquad $j,k=1,\ldots,n$.}\]
Here $n>1$ is a positive integer describing how many discretization points 
we will use.  For example, when $n=5$ we would have the points

\[ s_1 = 0.10, \quad s_2 = 0.30, \quad s_3 = 0.50, \quad s_4 = 0.70, \quad s_5 = 0.90.\]

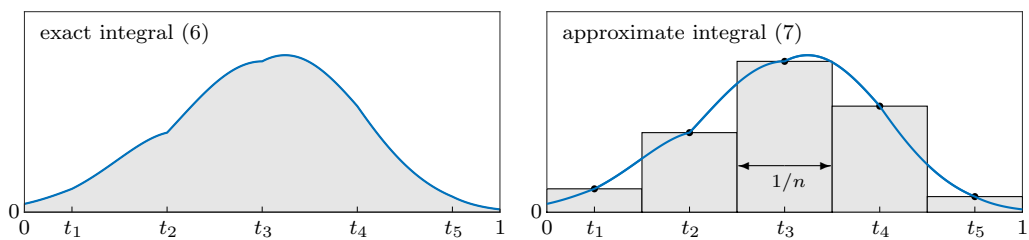
\begin{figure}[h!]
\begin{center}
\hspace*{0pt} 
\begin{minipage}{2.55in}
\begin{center}
\includegraphics[width=2.5in]{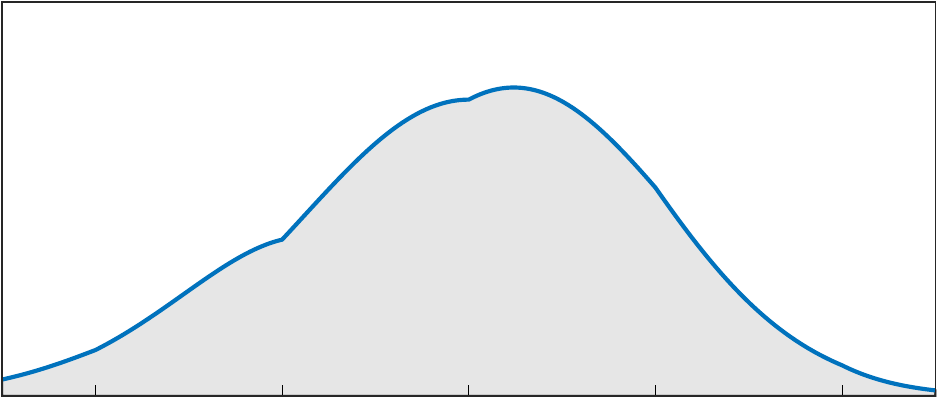} 

\begin{picture}(0,0)
\put(-96,10){\footnotesize $0$}
\put(-92,4){\footnotesize $0$}
\put(-75,4){\footnotesize $t_1$}
\put(-39,4){\footnotesize $t_2$}
\put(-3,4){\footnotesize $t_3$}
\put(33,4){\footnotesize $t_4$}
\put(69,4){\footnotesize $t_5$}
\put(88,4){\footnotesize $1$}
\put(-84,78){\footnotesize exact integral~(\ref{ME:dconv})}
\end{picture}
\end{center}
\end{minipage}
\quad
\begin{minipage}{2.55in}
\begin{center}
\includegraphics[width=2.5in]{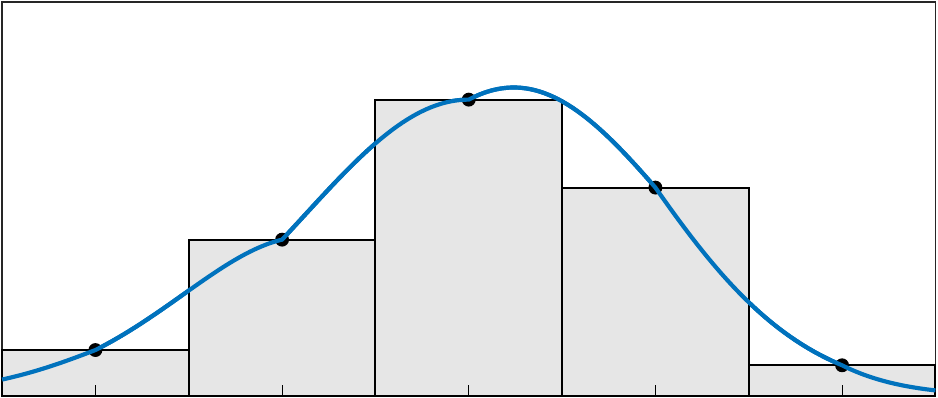}

\begin{picture}(0,0)
\put(-96,10){\footnotesize $0$}
\put(-92,4){\footnotesize $0$}
\put(-75,4){\footnotesize $t_1$}
\put(-39,4){\footnotesize $t_2$}
\put(-3,4){\footnotesize $t_3$}
\put(33,4){\footnotesize $t_4$}
\put(69,4){\footnotesize $t_5$}
\put(88,4){\footnotesize $1$}
\put(-84,78){\footnotesize approximate integral~(\ref{ME:midpt})}
\put(0,30){\vector(-1,0){17.5}}
\put(0,30){\vector(+1,0){17.5}}
\put(-5,22){\scriptsize $1/n$}
\end{picture}
\end{center}
\end{minipage}
\end{center}

\vspace*{-12pt}
\caption{\label{ME:cartoon_quad}
The integral~(\ref{ME:dconv}) is the area under the blue curve $h(s_j,t)f(t)$
in the left plot.  Approximate this quantity with the area of the rectangles
of width $1/n$ that touch $h(s_j,t)f(t)$ at $t=t_1,\ldots t_n$, 
as shown on the left and described in equation~(\ref{ME:midpt}).
This move allows us to replace a calculus problem
with a linear algebra problem.
}
\end{figure}

\vspace*{-4pt}
\begin{center}
\includegraphics[width=4in]{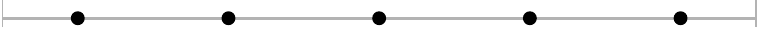}

\begin{picture}(0,0)(0,3)
\put(-145.3,7){\footnotesize $0$}
\put(-119.5,7){\footnotesize $s_1$}
\put(-61,7){\footnotesize $s_2$}
\put(-4,7){\footnotesize $s_3$}
\put(54,7){\footnotesize $s_4$}
\put(111,7){\footnotesize $s_5$}
\put(141.5,7){\footnotesize $1$}
\end{picture}
\end{center}
We will only compute the blurred values $b(s)$ of $f(t)$ at the grid points $s_j$,
so equation~(\ref{ME:conv}) becomes
\begin{equation} \label{ME:dconv}
\kern-29pt    b(s_j) = \int_0^1 h(s_j,t)@f(t)\,\dop t, \rlap{\kern48.7pt $j=1,\ldots,n$.}
\end{equation}
We want to \emph{approximate} $b(s_j)$ by only accessing the original function 
$f(t)$ at the grid points $t_k$.  
You likely encountered this idea when you first learned about integration:
you can approximate the area under the curve by the sum of the areas of some
rectangles that touch the curve at a point.  
This construction is known as \emph{Riemann sum}.  
As the rectangles get narrower, the mismatch in area decreases.
We will use rectangles of width $1/n$, and center them at the 
values of $t_k$, a process known as the \emph{midpoint rule}.
Figure~\ref{ME:cartoon_quad} shows a sketch of the process.
(We could use a fancier method to approximate the integral,
like the trapezoid rule or Simpson's rule~\cite[chap.~7]{SM03}, 
but that would be tangential to our purpose here.)

Applying the \emph{midpoint rule} to~(\ref{ME:dconv}) gives
\begin{equation} \label{ME:midpt}
    b(s_j) \approx \sum_{k=1}^n \frac{1}{n} \big(h(s_j, t_k) f_k\big) =: b_j,
   \rlap{\qquad $j=1,\ldots,n$,}
\end{equation}
where $f_k = f(t_k)$.  Notice that we have introduced the notation $b_j$ to 
refer to the \emph{approximation} to $b(s_j)$ obtained by the midpoint rule.

Given this notation, you can regard $f_k$ as the intensity of the $k$th pixel in our ``image;''
similarly, $b_j$ represents the blurred value of the $j$th pixel.
Our next challenge is to understand how the various $\{f_k\}$ values relate 
to the $\{b_j\}$ values.

Taking $j=1,\ldots, n$ gives us a set of $n$ equations that has a 
very regular structure:
\begin{eqnarray*}
{1\over n} \big( h(s_1,t_1) f_1 + h(s_1,t_2) f_2 + \cdots + h(s_1,t_n) f_n \big) &\;=\;& b_1 \\[3pt]
{1\over n} \big( h(s_2,t_1) f_1 + h(s_2,t_2) f_2 + \cdots + h(s_2,t_n) f_n \big) &=& b_2 \\[3pt]
& \vdots & \\[3pt]
{1\over n} \big( h(s_n,t_1) f_1 + h(s_n,t_2) f_2 + \cdots + h(s_n,t_n) f_n \big) &=& b_n.
\end{eqnarray*}
We envision $n$ being a large number (we will use $n=500$ later), 
and so keeping track of these $n$ equations would be quite tedious.
Thankfully the equations are all \emph{linear} in the values $f_1, \ldots, f_n$, 
and so we can collect these equations in a matrix-vector form:
\begin{equation} \label{ME:blurmat}
 \underbrace{{1\over n} \left[\begin{array}{cccc}
h(s_1,t_1) & h(s_1,t_2) & \cdots & h(s_1, t_n) \\[.25em]
h(s_2,t_1) & h(s_2,t_2) & \cdots & h(s_2, t_n) \\[.25em]
\vdots & \vdots & \ddots & \vdots \\[.25em]
h(s_n,t_1) & h(s_n,t_2) & \cdots & h(s_n, t_n)
\end{array}\right]}_{\mbox{$\BA$}}
\underbrace{\left[\begin{array}{c}
f_1 \\[.25em] f_2 \\[.25em] \vdots \\[.25em] f_n
\end{array}\right]}_{\mbox{$\Bf$}}
=
\underbrace{\left[\begin{array}{c}
b_1 \\[.25em] b_2 \\[.25em] \vdots \\[.25em] b_n
\end{array}\right]}_{\mbox{$\Bb$}}.
\end{equation}
Thus, the $n$ equations in~(\ref{ME:midpt}) can be expressed as the
matrix-vector product $\BA\Bf = \Bb$.  
(Take note of the $1/n$ factor on the far left of~(\ref{ME:blurmat}):
be sure to include it in the matrix $\BA$, 
so that the $(j,k)$ entry is $a_{j,k} = h(s_j,t_k)/n$.)

\medskip
The code in Figure~\ref{ME:fig:buildA} 
shows how to construct the matrix $\BA$ in Python.
Note that the function \verb|build_blur_A| builds the matrix 
\emph{one row at a time}, using only a single \verb|for| loop.
(This routine uses Gaussian blurring, but any other kernel could readily be swapped in.)

\begin{figure}[t!]
\begin{lstlisting}[language=Python]
import numpy as np

def h_average(s,t,z=.025):                      # averaging kernel 
    return np.double(np.abs(s-t)<=z)/(2*z)

def h_hat(s,t,z=.025):                          # hat-function kernel 
    return np.maximum(0,1-np.abs(s-t)/z)/z

def h_gaussian(s,t,z=.025):                     # Gaussian kernel     
    c = 1/(np.sqrt(np.pi)*z)
    return c*np.exp(-np.power(s-t,2)/(z**2))

def build_blur_A(n=100,z=.025):
    A = np.zeros((n,n));                        # n-by-n matrix of zeros
    s = np.array([(j+.5)/n for j in range(n)])  # s_j points
    t = np.array([(k+.5)/n for k in range(n)])  # t_k points
    for j in range(0,n):
        A[j,:] = h_gaussian(s[j],t,z)/n         # A(j,k) = h(s_j,t_k)/n
    return A
\end{lstlisting}
\vspace*{-7pt}
\caption{\label{ME:fig:buildA}
Python code to construct the three blurring functions~(\ref{ME:avgker}), (\ref{ME:hatker}), 
and (\ref{ME:gaussker}), and to build the matrix $\BA$ defined in~(\ref{ME:blurmat}).}
\end{figure}

\medskip
The matrix-vector form~(\ref{ME:blurmat}) has many advantages: for one thing, 
it provides a tidy way to organize the $n$ equations in $n$ unknown variables.
The structure $\BA\Bf=\Bb$ encourages us think about the equation at a higher 
level of abstraction, enabling us to apply tools -- both theoretical ideas 
and mathematical software -- developed for general problems that share 
this same structure.  

\medskip
We can formulate a basic algorithm for approximating the blurring operation:
\begin{enumerate} \setlength{\itemsep}{0pt} \setlength{\parskip}{0pt}
\item Select a grid size, $n$.
\item Construct the grid points, $t_1, \ldots, t_n$.
\item Build the $n\times n$ matrix $\BA$ from the blurring kernel $h(s,t)$.
\item Sample the function $f(t)$ at the grid points $t_k$, $k=1,\ldots, n$.\\
      Collect the results in the vector $\Bf \in \R^n$.
\item Compute the matrix-vector product $\BA\Bf = \Bb$ to find the blurred vector, $\Bb$.
\end{enumerate}

\begin{figure}[h!]
\medskip
\begin{center}
  \includegraphics[width=2.58in]{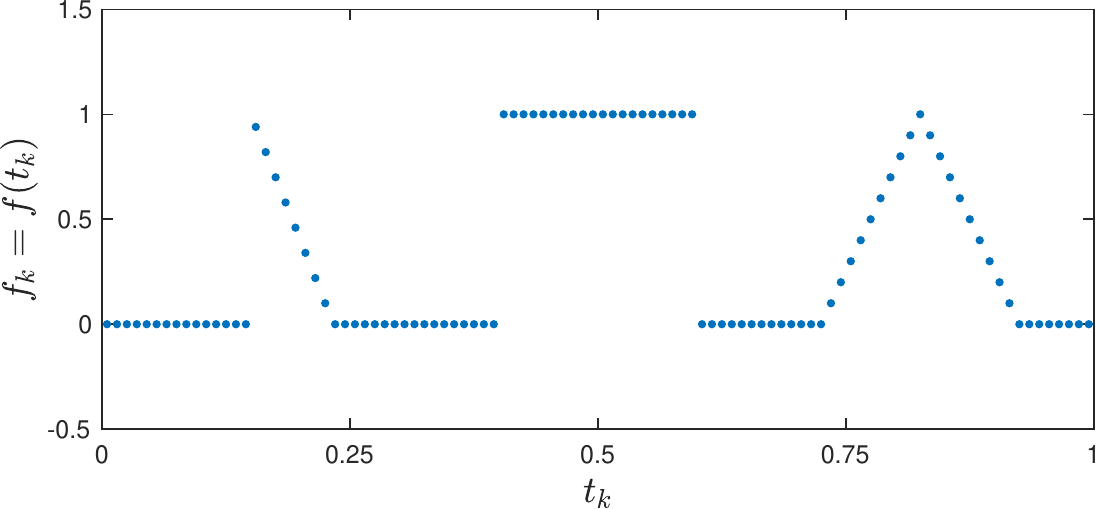}
 
 \includegraphics[width=2.58in]{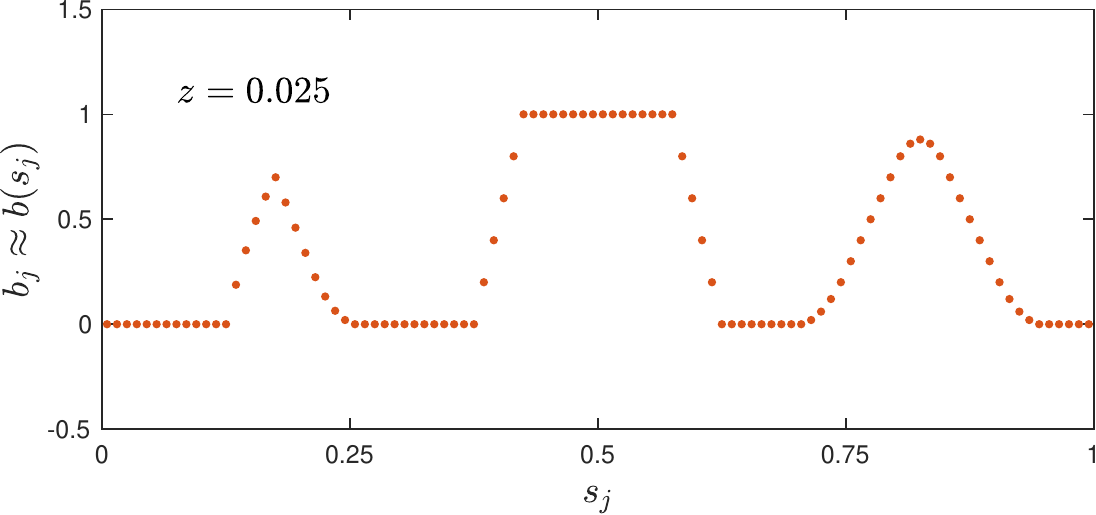}
\end{center}

\vspace*{-1.5em}
\caption{\label{ME:fig:fw}
The original signal $f(t)$, sampled at $t=t_k$ for $k=1,\ldots,100=n$ (top),
and the blurred values from $\Bb = \BA\Bf$, using the simple blurring 
kernel~(\ref{ME:avgker}) with $z=0.025$.
Compare these discretized versions to the continuous (exact) versions in
Figure~\ref{ME:fig:f1} and the top of Figure~\ref{ME:fig:b1}.
}
\end{figure}

Figure~\ref{ME:fig:fw} shows how this blurring operation works 
for the hat-function kernel with $z=0.025$ and $n=100$.
To the eye, the approximate blurring operation computed as $\Bb=\BA\Bf$
looks quite accurate. (Compare the bottom plot in Figure~\ref{ME:fig:fw}
to the top plot in Figure~\ref{ME:fig:b1}.)

We have seen how the calculus operation of blurring can be approximated 
by a basic linear algebra problem:  blurring a signal reduces to simple 
matrix-vector multiplication.
You can imagine applications where you would like to apply such blurring,
e.g., obscuring the face of a bystander in a photograph of a crime scene.
However, in many more situations we actually want to do the \emph{opposite}
of blurring:  We acquire a blurry signal (from a camera, or a telescope)
and we want to sharpen it up.  We want to \emph{deblur}.
Linear algebra immediately suggests a way to do this: 
\emph{invert} the blurring matrix.  
We will explore this possibility in the next section.

\begin{reflections}
\item The most aggressive blurring you could imagine would simply 
take $h(s,t) = 1$ for all $s,t \in [0,1]$.  
Describe the matrix $\BA$ for this choice of blurring function.
What would $\Bb = \BA\Bf$ look like?  
If you are given $\Bf$, could you solve for the unknown $\Bb$?

\item The matrix $\BA$ is said to be \emph{banded} when $a_{j,k} = 0$
provided $|j-k|$ is sufficiently large.  
Which of the three kernels studied in this section lead to a banded blurring
matrix?  How does the ``width'' of the band (the number of nonzeros per row)
depend on $z$?
\end{reflections}

\noindent
Now would be a good time to explore Exercises~\ref{ME:ex:fw}--\ref{ME:ex:Atime} 
starting on page~\pageref{ME:ex:fw}.

\section{The inverse problem: \emph{deblurring}}

In the last section we saw the blurring operation described by $\BA\Bf = \Bb$,
and we used knowledge of the blurring operation $\BA$ and the original signal
$\Bf$ to compute the blurry $\Bb$.
This set-up is called the \emph{forward problem}.

But now let us turn the tables.
Suppose we acquire (that is, we measure) some blurry signal $\Bb$,
and want to \emph{recover} the unknown sharp signal, $\Bf$. 
We want to \emph{deblur} $\Bb$.
Linear algebra gives an immediate way to accomplish this goal:
If $\BA\Bf = \Bb$, then we should have
\begin{equation} \label{ME:Ainv}
 \Bf = \BAs^{-1}\Bb,
\end{equation}
provided $\BA\in\R^{n\times n}$ is an invertible matrix.  
Since~(\ref{ME:Ainv}) involves the inverse of the blurring operation $\BA$,
we call this deblurring process an \emph{inverse problem}.
What could go wrong with the clean formula~(\ref{ME:Ainv})?

\medskip
We must first acknowledge that the $\BA\Bf=\Bb$ model is only a rough description
of the true process of blurring.  Here are a few ways the model falls short.
\begin{itemize} \setlength{\itemsep}{0pt} \setlength{\parskip}{0pt}
\item The integral model in~(\ref{ME:conv}) is an imperfect description of 
the blurring process.
\item The kernel $h(s,t)$ and the blurring parameter $z$ are only estimates 
of the properties of the instrument (camera, environmental conditions, etc.).
\item The discretization process that turned~(\ref{ME:conv}) into the 
linear algebra problem~(\ref{ME:blurmat}) introduced an additional approximation.
\item Measurements of the vector $\Bb$ (e.g., with a camera, telescope, or microscope)
will inevitably be imperfect, adding a layer of (hopefully random) errors 
that we will call ``noise.''
\end{itemize}
Humbled by this list of shortcomings, we might wonder if our model is any good at all!
Indeed, simple computational experiments do little to build our confidence.
Let us focus on the last flaw on our list, measurement errors in $\Bb$.  
Return to the example shown in Figure~\ref{ME:fig:fw},
but suppose that instead of exactly measuring $\Bb$, we acquire some noisy version, $\bnoise$.
Let the $j$th entry of this vector be
\[ (\bnoise)_j = b_j + e_j,\]
where $b_j$ is the $j$th entry of $\Bb$ and 
$e_j$ is a normal random variable of mean zero and standard deviation $\eps \|\Bb\|$,
for some small value of $\eps$.  We are using the Euclidean vector norm defined by
\[ \|\Bb\| = \sqrt{\Bb^T\Bb} = \Big(\sum_{j=1}^n b_j^2 \Big)^{1/2},\]
so that the magnitude of the noise that pollutes $\bnoise$ is calibrated 
to match the magnitude of $\Bb$.

\begin{figure}[b!]
\begin{center}
 \includegraphics[height=1.37in]{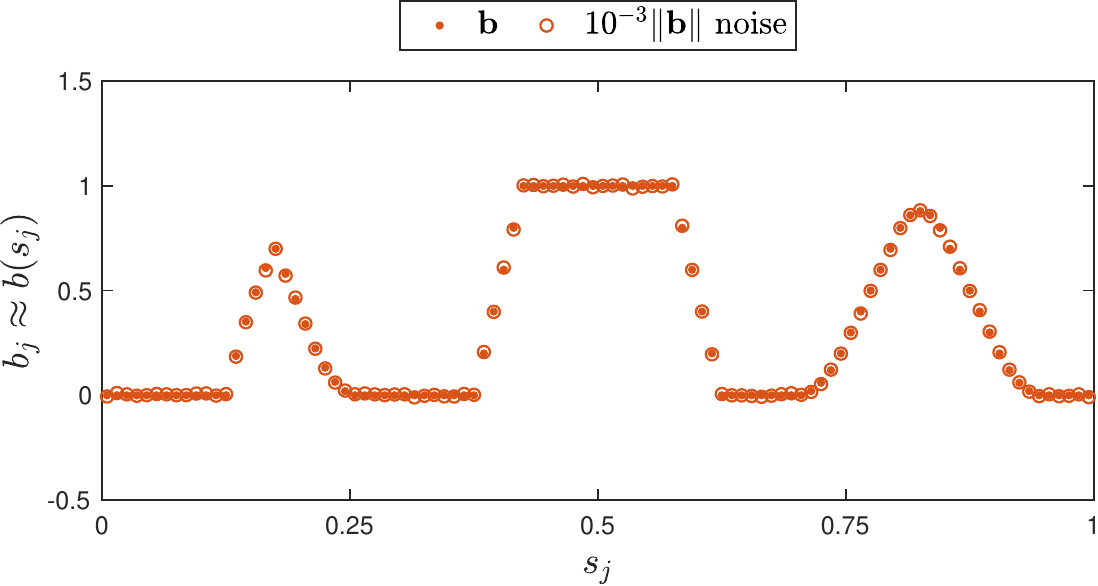}\quad
 \includegraphics[height=1.37in]{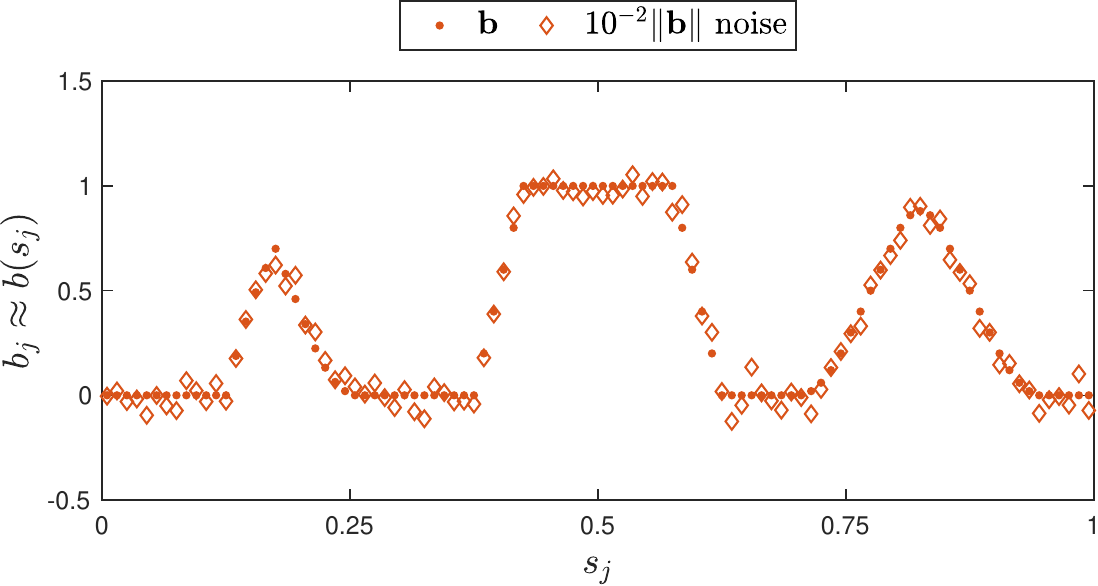}
 
\medskip
 \includegraphics[height=1.37in]{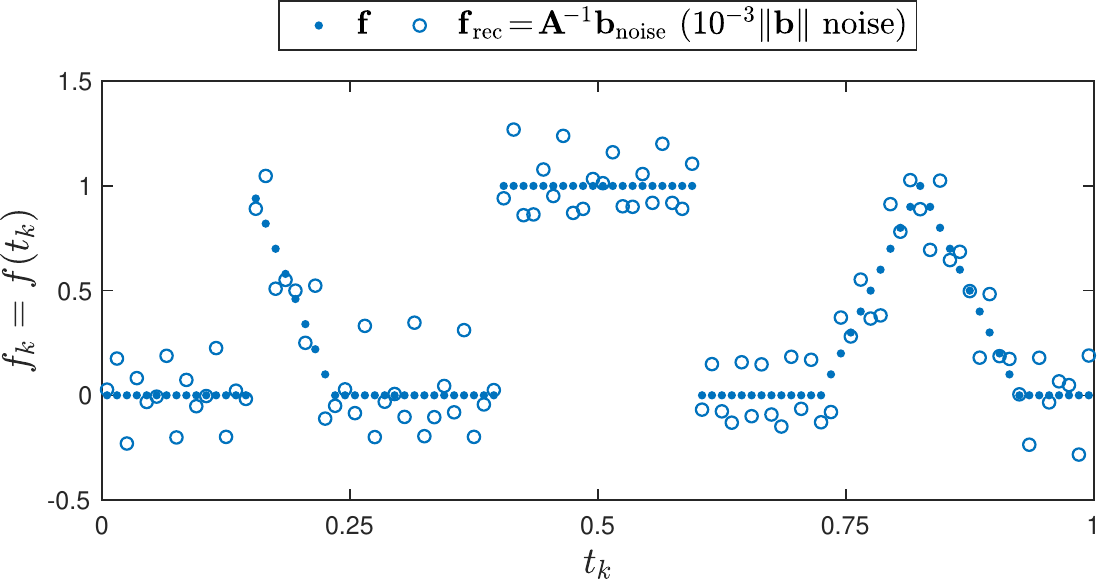}\quad
 \includegraphics[height=1.37in]{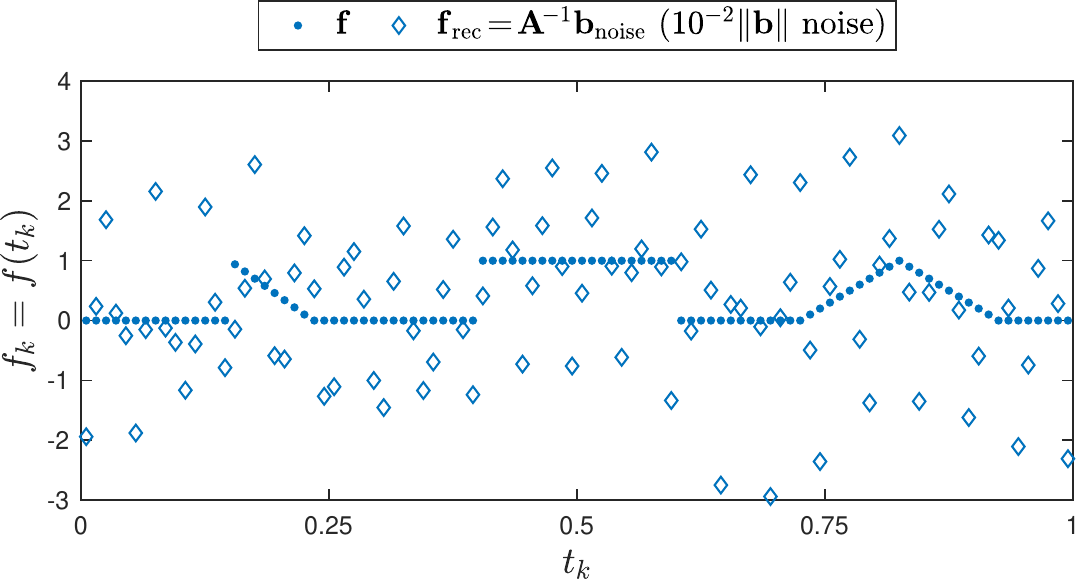}
\end{center}

\vspace*{-17pt}
\caption{\label{ME:fig:bw}
An attempt to recover $\Bf$ from noisy measurements
($n=100$, averaging kernel~(\ref{ME:avgker}) with $z=0.025$).
The true $\Bb = \BA\Bf$ is polluted with some random noise 
(normally distributed, mean~zero, standard deviation $10^{-3}\|\Bb\|$ (left)
or $10^{-2}\|\Bb\|$ (right)), giving $\bnoise$.
As seen in the top figures, 
$\bnoise$ differs very little from $\Bb$ for both noise levels.
We then recover the unblurred signal by computing $\frec = \BAs^{-1}\bnoise$.
As seen in the bottom figures, a whiff of noise in $\Bb$ causes
a large change to the recovered $\Bf$.
The values in $\frec$ ($\circ$ and $\diamond$) should fall
on top of the true values in $\Bf$ (\raisebox{2pt}{\tiny$\bullet$}), 
but most are far from correct.
}
\end{figure}

\begin{figure}[t!]
\begin{center}
 \includegraphics[height=1.37in]{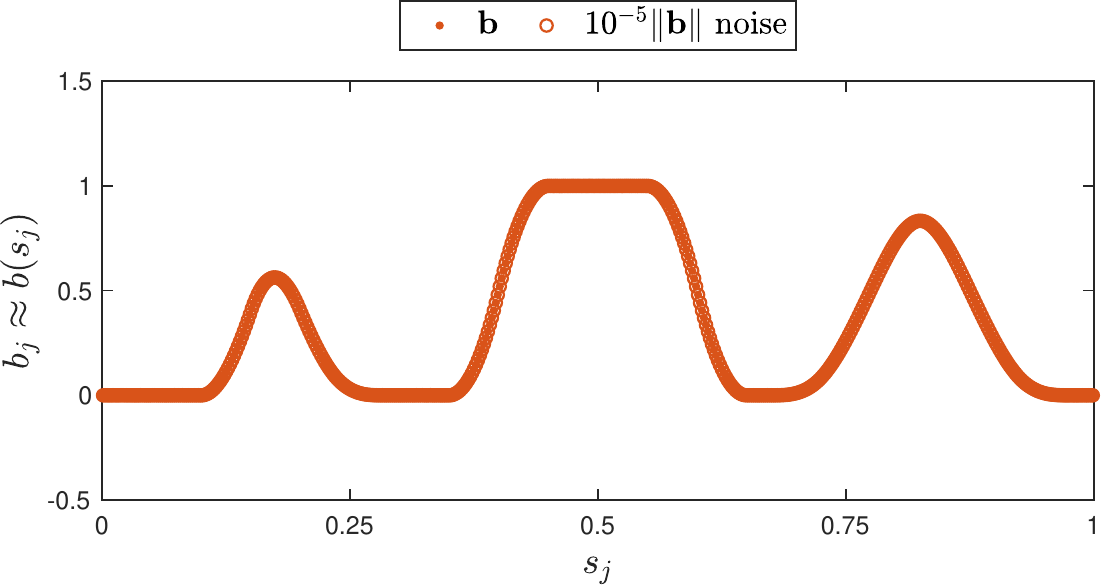}\quad
 \includegraphics[height=1.37in]{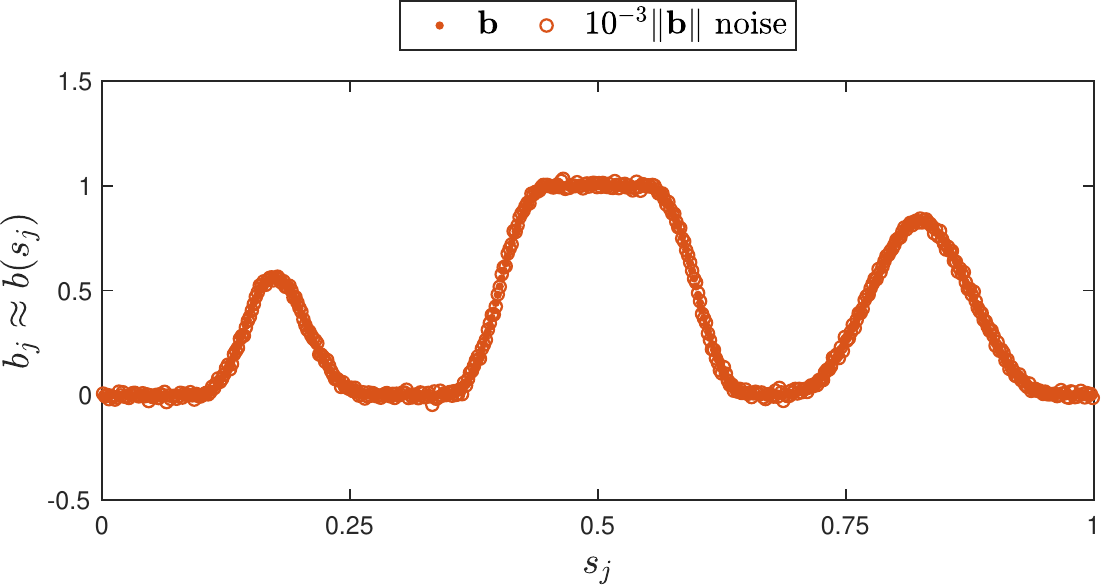}
 
\medskip
 \includegraphics[height=1.35in]{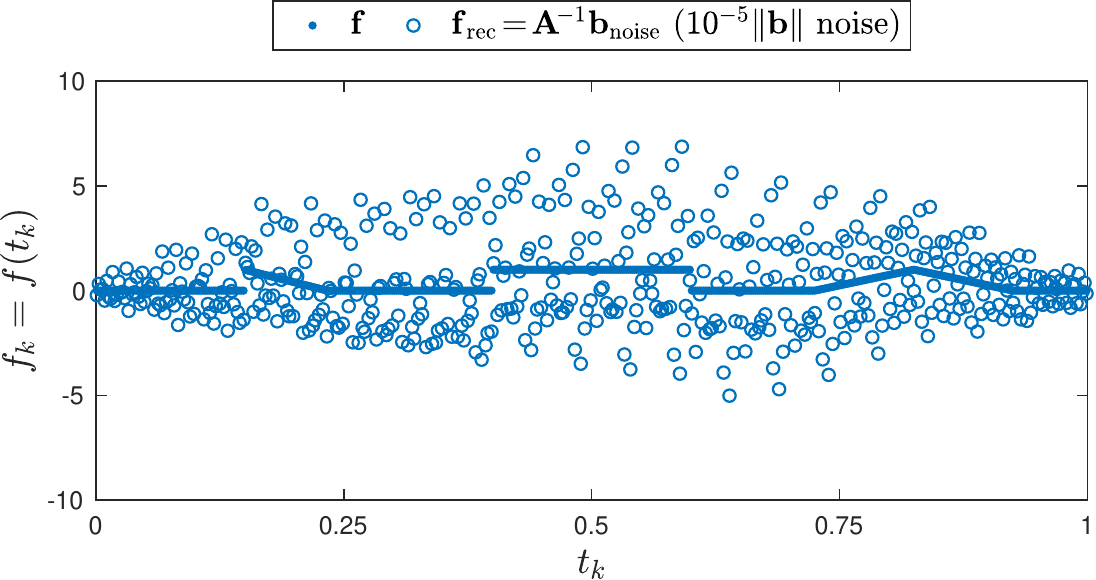}\quad
 \includegraphics[height=1.35in]{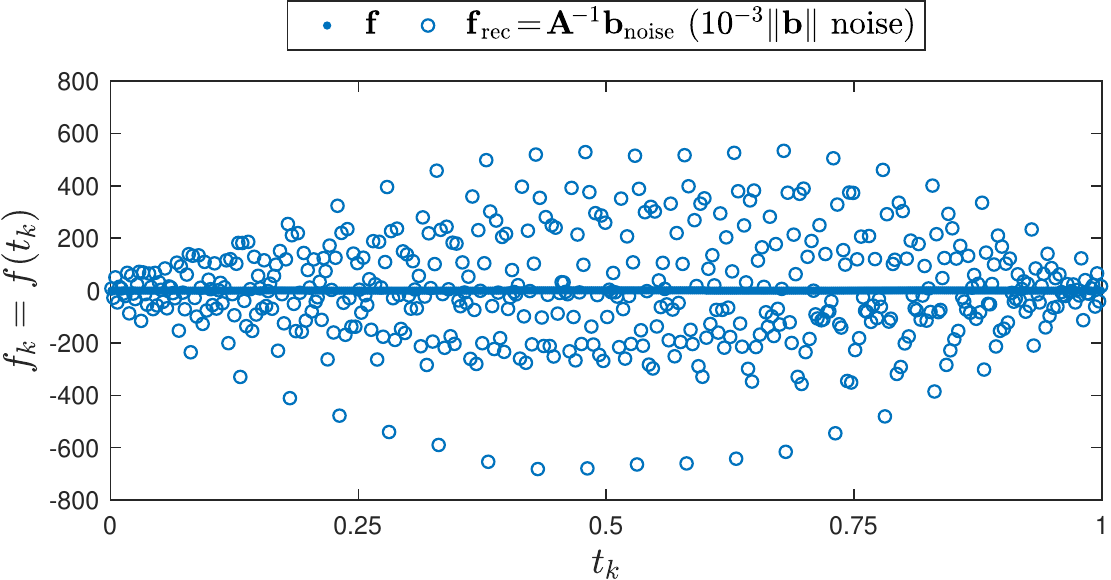}
\end{center}

\vspace*{-17pt}
\caption{\label{ME:fig:bw2}
Repetition of Figure~\ref{ME:fig:bw}, but now using $n=500$
with the hat-function kernel~(\ref{ME:hatker}) with $z=0.05$
and noise levels $10^{-5}\|\Bb\|$ (left) and $10^{-3}\|\Bb\|$ (right).
The very subtle level of noise in $\bnoise$ on the top-left is essentially 
unnoticeable to the eye, yet it causes profound errors in $\frec$ in the
bottom-left plot.  When the noise increases by two orders of magnitude
(top-right), so too do the errors in $\frec$.
}
\end{figure}

\medskip
We might instinctively expect that an error of relative size $\eps$ in $\bnoise$ 
will cause an error of similar  error in $\frec = \BAs^{-1}\bnoise$.
Figure~\ref{ME:fig:bw} shows that this is far from the case!
Errors with $\eps=10^{-3}$ and $\eps=10^{-2}$ in $\bnoise$ 
(just $0.1\%$ and $1.0\%$ relative noise level)
cause much larger errors in $\frec$.
The crisp image $\Bf$ is not recovered at all; 
the recovered values are bad for $\eps=10^{-3}$ and 
entirely useless for $10^{-2}$.  
(Indeed, the \emph{blurry, noisy} version $\bnoise$ gives us a better 
impression of the true $\Bf$ than the supposedly ``deblurred'' 
vector $\frec$ does!)
Lest one think there is something special about this particular example, 
in Figure~\ref{ME:fig:bw2} we repeat the experiment on a larger problem 
$(n=500)$ with the hat-function kernel~(\ref{ME:hatker}) with $z=0.05$.
Even the noise level $10^{-5}\|\Bb\|$ causes catastrophic errors in 
$\frec = \BAs^{-1}\bnoise$.

\begin{figure}[t!]
 \hspace*{95pt}\includegraphics[width=1.25in]{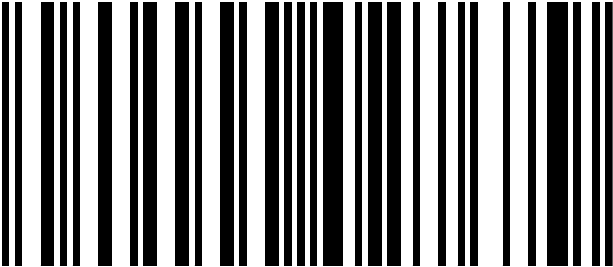}
 \hspace*{91pt} \raisebox{16pt}{\small \emph{the barcode}}\\[10pt]
\medskip
 \includegraphics[width=3.75in]{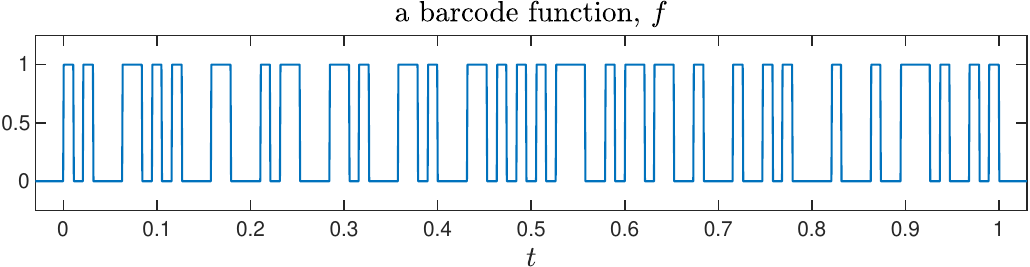} 
 \hspace*{8pt}\raisebox{35pt}{\small \emph{the original function}}\\
\medskip
 \includegraphics[width=3.75in]{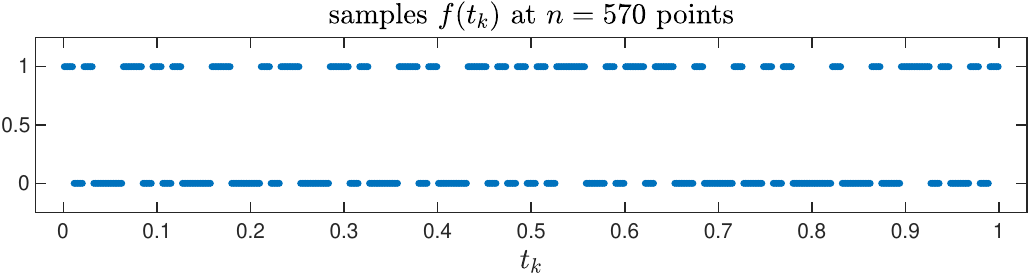}
 \hspace*{8pt}\raisebox{35pt}{\small \emph{sample $f$ at $t_k$ to get $\Bf$}}\\
\medskip
 \includegraphics[width=3.75in]{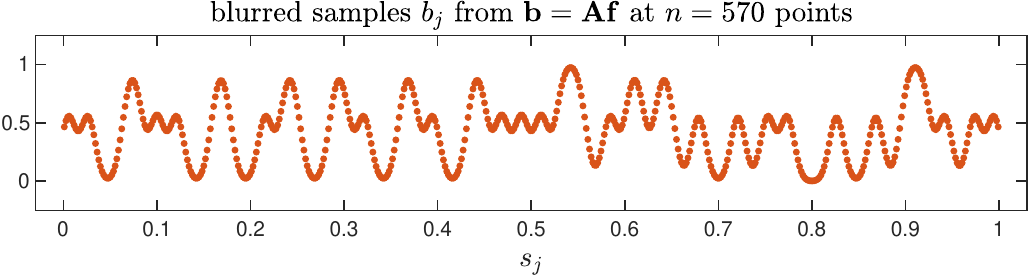}
 \hspace*{2pt}\raisebox{38pt}{\small \begin{tabular}{l}
 \emph{blur $\Bf$ to get $\Bb=\BA\Bf$}\\[3pt] Gaussian blur~(\ref{ME:gaussker})\\ with $z=0.01$\end{tabular}} \\
\medskip
 \includegraphics[width=3.75in]{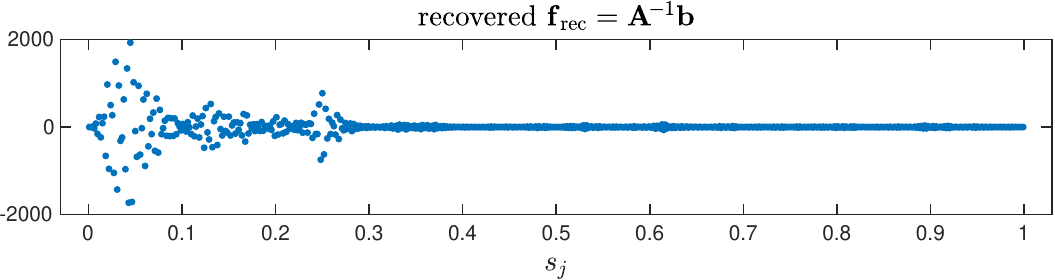}
 \hspace*{2pt}\raisebox{38pt}{\small \begin{tabular}{l} 
               \emph{undo the blur}\\ $\frec = \BAs^{-1}\Bb$ \\[3pt] (We \emph{should} get $\frec = \Bf$.)
               \end{tabular}}

\vspace*{-17pt}
\caption{\label{ME:fig:barcode1}
Deblurring our first barcode.  
A Universal Product Code (UPC) is shown at the top.  
We translate this barcode into a function $f$ that takes values 0~(white bar) and 1~(black bar).  
When sampled on a grid of $n=570$ points, 
the vector $\Bf\in\R^{570}$ only contains the values zero and one.
Blur this vector using the Gaussian kernel~(\ref{ME:gaussker}) 
with $z=0.01$ to get $\Bb=\BA\Bf$, shown in red.
Attempt to undo the blur by computing $\frec=\BAs^{-1}\Bb$.
Instead of recovering the true $\Bf$ (zeros and ones), we get much larger values!
Though we have not injected any extra noise, the inversion process is 
sufficiently sensitive to small numerical errors to make the result useless.}
\end{figure}

It is tempting to conclude that we need better measurements, that the problem can be
overcome if we can only obtain greater precision.
However, the quest for more precision is a fool's errand.
The difficulty exposed by this simple example is chronic, becoming increasingly problematic
for more realistic signals and blurring kernels.  
Figure~\ref{ME:fig:barcode1} shows such an example, reflecting the 
kinds of barcodes that we will focus upon later in this manuscript.
In this case we start with a vector $\Bf\in\R^{570}$ containing only the values
zero and one.  
Using the Gaussian kernel~(\ref{ME:gaussker}) with $z=0.01$, 
we blur this function to obtain $\Bb=\BA\Bf$.  
Now \emph{without injecting any intentional noise}, 
we simply try to undo 
this matrix vector product: $\frec = \BAs^{-1}\Bb$.
The small numerical rounding errors that occur when we compute with double-precision
floating point arithmetic are sufficient to make the inversion process go haywire:
while we hope for $\frec$ to be a vector containing only zeros and ones,
$\frec$ contains much larger values -- \emph{indeed, values over 100 times too large}.

In the face of such failure we might be tempted to give up
(or look for a bug in our code).  
Instead, we will use this as an opportunity to gain insight into the subtlety of 
matrix inversion, and to find a more robust approach to solving this kind of 
fragile inverse problem.


\begin{reflections}
\item In this section we suggest four potential sources of error 
       in the modeling process.  
       Which of these do you think is most significant, and which least significant?
       Can you think of any other sources of error that we have not mentioned?

\item  Can you imagine two vectors $\Bf_1 \not\approx \Bf_2$ that would be very similar
       when blurred?  That is, $\BA\Bf_1 \approx \BA\Bf_2$?  
       (Another way of thinking about this question:  Is there a vector $\Bz$
        that gets blurred to almost nothing: $\BA\Bz \approx \Bzero$?  In that case,
        set $\Bf_2 = \Bf_1 + \Bz$.)
\end{reflections}

\noindent
Now would be a good time to explore Exercises~\ref{ME:ex:bw1}--\ref{ME:ex:bign}
starting on page~\pageref{ME:ex:bw1}.

\section{What is wrong with \boldmath $\BAs^{-1}$?} \label{ME:sec:Ainv}

The image deblurring failure has an intuitive explanation.
Think about blurring a photograph.  Sharp lines are smoothed; fine distinguishing 
features often wash out.  Two similar images that have distinct details could look
nearly identical once they are blurred.  We are asking the inversion process to 
recover those fine details that are obscured by the blurring, a challenging task.

\begin{figure}[h!]
\medskip
\begin{center}
\includegraphics[width=1.75in]{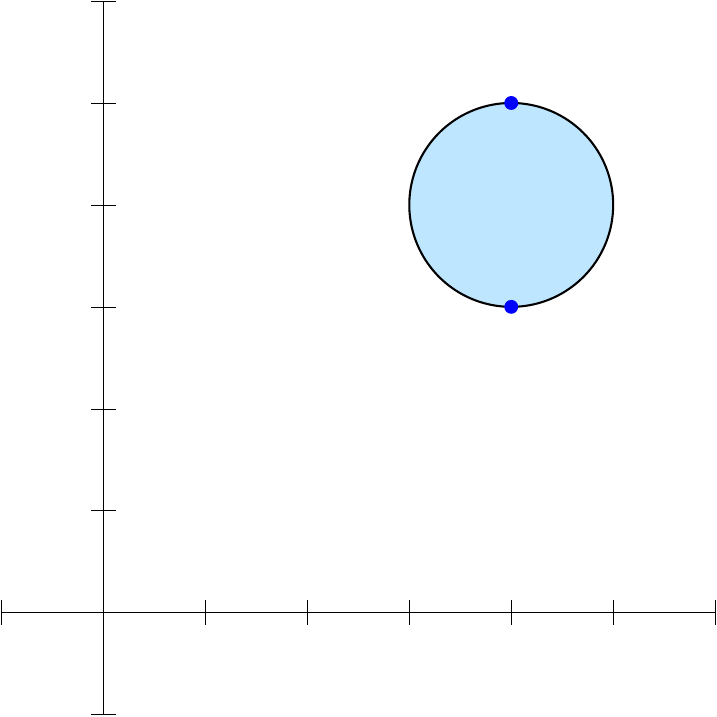}
\begin{picture}(0,0)
\put(-98,10){\sf \tiny 0.25}
\put(-80,10){\sf \tiny 0.50}
\put(-62,10){\sf \tiny 0.75}
\put(-50.5,1){\footnotesize $\alpha$}
\put(-134,79){\footnotesize $\beta$}
\put(-44,10){\sf \tiny 1.00}
\put(-26,10){\sf \tiny 1.25}
\put( -8,10){\sf \tiny 1.50}
\put(-124,34.5){\sf \tiny 0.25}
\put(-124,52.5){\sf \tiny 0.50}
\put(-124,70.5){\sf \tiny 0.75}
\put(-124,88.5){\sf \tiny 1.00}
\put(-124,106.5){\sf \tiny 1.25}
\put(-124,124.5){\sf \tiny 1.50}
\put(-42,111.5){\footnotesize $\Bf_1$}
\put(-43,64.5){\footnotesize $\Bf_2$}
\end{picture}
\hspace*{4em}
\includegraphics[width=1.75in]{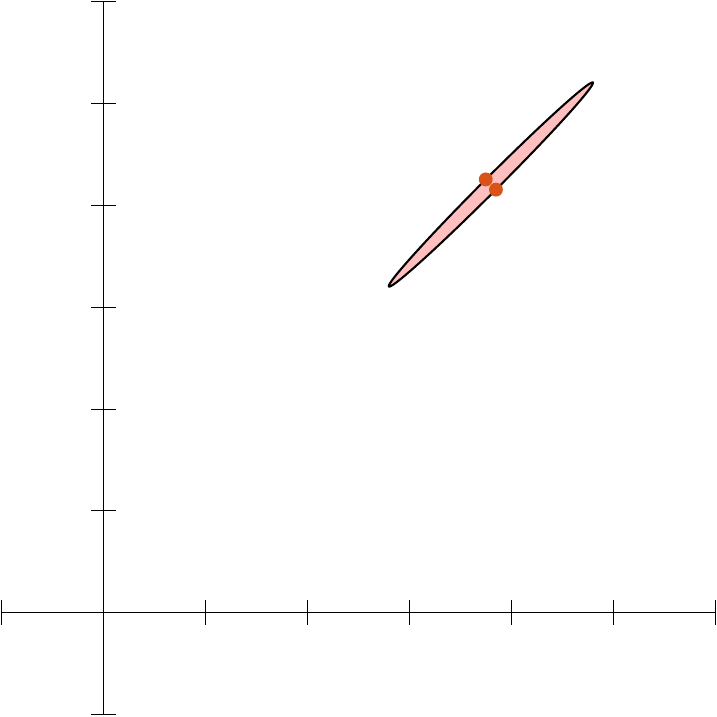}
\begin{picture}(0,0)
\put(-98,10){\sf \tiny 0.25}
\put(-80,10){\sf \tiny 0.50}
\put(-62,10){\sf \tiny 0.75}
\put(-44,10){\sf \tiny 1.00}
\put(-26,10){\sf \tiny 1.25}
\put( -8,10){\sf \tiny 1.50}
\put(-124,34.5){\sf \tiny 0.25}
\put(-124,52.5){\sf \tiny 0.50}
\put(-124,70.5){\sf \tiny 0.75}
\put(-124,88.5){\sf \tiny 1.00}
\put(-124,106.5){\sf \tiny 1.25}
\put(-124,124.5){\sf \tiny 1.50}
\put(-76,97.5){\footnotesize $\Bb_1 \!=\! \BA\Bf_1$}
\put(-42,86){\footnotesize $\BA\Bf_2\! =\! \Bb_2$}
\end{picture}
\end{center}

\vspace*{-17pt}
\caption{\label{ME:fig:invsketch}
On the left, the domain $\Bf\in\R^2$, with a ball of $\Bf$ vectors highlighted in blue.
On the right, the image $\BA\Bf \in \R^2$, with the red ellipse showing the image of the
blue disk on the left.  Note that the two vectors $\Bf_1$ and $\Bf_2$ are distant,
but their images $\BA\Bf_1$ and $\BA\Bf_2$ are quite close.
}
\end{figure}
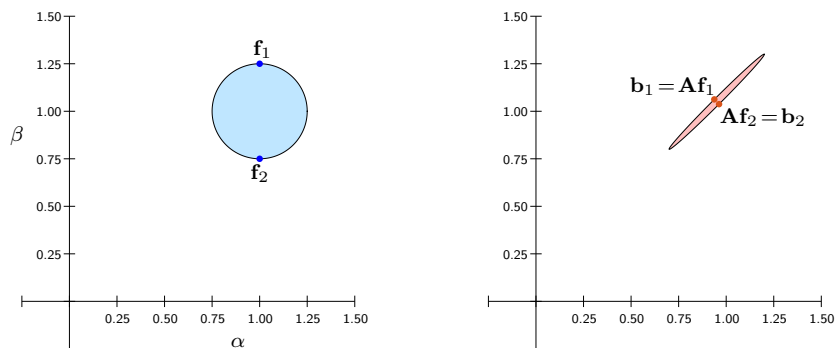

While this explanation might make intuitive sense, we can better understand the
issue through the language of linear algebra, and we do not need anything as 
complicated as a blurring matrix to see what is going on.
Instead, consider the simple $2\times 2$ matrix
\[ \BA = \left[\begin{array}{cc} 1 & -0.05 \\ 1 & \phantom{-}0.05 \end{array}\right].\]
This matrix has two linearly independent columns and two linearly independent rows;
the rank of the matrix is two, and it is invertible: $\det(\BA) = 0.10 \ne 0$.  
Indeed, you can compute
\[ \BAs^{-1} = \left[\begin{array}{rr} 0.5 & 0.5 \\ -10 & \phantom{-}10 \end{array}\right].\]
Notice that $\BA$ is an invertible matrix, but the rows of $\BA$ are \emph{nearly}
linearly dependent. (Think about plotting those two vectors -- they point
in nearly the same direction.)  
Even though the largest entry in $\BA$ is~1, $\BAs^{-1}$
has entries that are 10~times as large.  
Figure~\ref{ME:fig:invsketch} illustrates the implications.  
Consider the blue patch on the left, which contains all vectors of the form
\[ \Bf = \left[\begin{array}{c} \alpha \\ \beta \end{array} \right] \]
that fall within a disk of radius~0.25 centered at $[1,1]^T$; i.e.,
\[ |\alpha-1|^2 + |\beta-1|^2 \le 0.25.\]
The matrix $\BA$ maps these vectors to the red patch on the right, a narrow ellipse.
Two representative vectors $\Bf_1$ and $\Bf_2$ are far apart (in the plot on the left), 
but they are mapped quite close together by $\BA$, 
so that $\Bb_1 = \BA\Bf_1$ and $\Bb_2 = \BA\Bf_2$ 
appear nearby (in the plot on the right).

Now imagine we make a small change to 
\[ \Bb_1 = \BA\Bf_1 
         = \left[\begin{array}{cc} 1 & -0.05 \\ 1 & \phantom{-}0.05 \end{array}\right]
           \left[\begin{array}{c} 1 \\1.25 \end{array}\right]
         = \left[\begin{array}{c} 0.9375 \\1.0625 \end{array}\right]\]
to obtain
\[ \bnoise = \Bb_1 + \Be 
         = \left[\begin{array}{c} 0.9375 \\1.0625 \end{array}\right]
           + \left[\begin{array}{r} 0.0125 \\ -0.0125 \end{array}\right]
         = \left[\begin{array}{c} 0.95 \\1.05 \end{array}\right].\]
Then when we ``recover'' $\Bf_1$ from $\bnoise$, we obtain
\[ \frec = \BAs^{-1} \bnoise 
         = \left[\begin{array}{rr} 0.5 & 0.5 \\ -10 & \phantom{-}10 \end{array}\right]
           \left[\begin{array}{c} 0.95 \\1.05 \end{array}\right]
         = \left[\begin{array}{c} 1.00 \\1.00 \end{array}\right] \ne 
           \left[\begin{array}{c} 1 \\1.25 \end{array}\right] = \Bf_1.\]
Changing the entries of $\Bb_1$ by just $\pm 0.0125$ changed the recovered signal 
$\frec$ by $0.25$ in the second entry. 
As the next exercise shows, even in this simple $2\times2$ matrix example,
we can adjust the entries of $\BA$ to make the discrepancy between $\frec$ and
$\Bf_1$ \emph{arbitrarily large}.

\begin{reflection}
\item The theme of this section is that 
      multiplication by $\BA^{-1}$ can map close-together
      vectors $\Bb_1$ and $\Bb_2$ to far-away vectors 
      $\Bf_1 = \BA^{-1} \Bb_1$  and $\Bf_2 = \BA^{-1} \Bb_2$.

      Can you apply this same concept to functions of a single real 
      variable, $g(x):\R\to\R$?  That is, can you find a function $g$
      such that $g(x+h)$ is very different from $g(x)$ 
      for small values of $|h|>0$? 
      How does this relate to the derivative $g'(x)$?
\end{reflection}

\noindent
Now would be a good time to explore Exercise~\ref{ME:ex:2x2} 
on page~\pageref{ME:ex:2x2}.

\section{A Remedy from Regularization} \label{ME:sec:reg}

We have sought to compute $\Bf$ as the unique solution to the linear system
\[ \BA \Bf = \Bb,\]
which we simply write as
\[ \Bf = \BAs^{-1}\Bb.\]
When we try to compute the ``exact solution'' $\Bf$ and get some nonsensical
values (like we saw in Figure~\ref{ME:fig:barcode1}), we realize that 
satisfying the linear system $\BA\Bf=\Bb$ is only one of our objectives.  
We also want $\Bf$ to have realistic, meaningful entries.
For example, when solving the barcode problem, we might ask for
every entry of $\Bf$ to be either zero or one.
(For reasons not immediately obvious, this turns out to be
quite a difficult condition to impose, computationally.)
A simpler requirement is to simply ask that $\|\Bf\|$ not be too large.

So, we want $\BA\Bf \approx \Bb$ while keeping $\|\Bf\|$ as small as possible.
Let us introduce a parameter $\lambda>0$, and define an objective function that
balances these two requirements: 
\begin{equation} \label{ME:phi}
\phi(\Bf) = \|\Bb - \BA\Bf\|^2 + \lambda^2 \|\Bf\|^2.    
\end{equation}
Our goal is to find the $\Bf\in\R^n$ that minimizes $\phi(\Bf)$.

\medskip
The act of penalizing the misfit $\|\Bb-\BA\Bf\|$ with some term that
controls the size of the solution (such as $\|\Bf\|$) is known as \emph{regularization}.
More specifically, mathematicians see the optimization of 
the objective function~(\ref{ME:phi}) as an example of \emph{Tikhonov regularization},
while statisticians call this approach \emph{ridge regression}.

\medskip
What role does the parameter $\lambda$ play?  Consider the extreme choices.
\begin{itemize} \setlength{\itemsep}{0pt} \setlength{\parskip}{0pt}
\item By setting $\lambda=0$, we entirely neglect $\|\Bf\|$, placing all our emphasis on 
      making $\BA\Bf=\Bb$.  (Indeed, this was our initial strategy, above:  take $\Bf = \BAs^{-1}\Bb$.)
\item By taking $\lambda\to\infty$, we place increasing importance on controlling $\|\Bf\|$,
while neglecting $\BA\Bf=\Bb$.  In the limit $\lambda\to\infty$, we would simply obtain $\Bf=\Bzero$:
a solution that indeed makes $\|\Bf\|$ small, but is not interesting as a solution to $\BA\Bf\approx \Bb$.
\end{itemize}

\medskip
How can we minimize $\phi(\Bf)$?  
We will approach the question from the perspective of multivariable calculus.  
(One can derive the same solution using only properties of linear algebra; 
see, e.g., \cite[chap.~7]{Emb3606}.)
To minimize a smooth function $\phi(\Bf)$, 
we should take the gradient (with respect to $\Bf$), 
set this gradient equal to zero, and solve for $\Bf$.
Recalling that we can expand norms using inner products (e.g., $\|\By\|^2 = \By^T\By$) 
and distribute inner products (e.g., $(\Bx+\By)^T\Bz = \Bx^T\Bz + \By^T\Bz$),
we can write:
\begin{eqnarray}
\phi(\Bf) &=& \|\Bb - \BA\Bf\|^2 + \lambda^2 @ \|\Bf\|^2 \nonumber \\[5pt]
          &=& (\Bb-\BA\Bf)^T(\Bb-\BA\Bf) + \lambda^2 @@\Bf^T\Bf \label{ME:eq:phi1}\\[5pt]
          &=& \Bb^T\Bb - 2 @@\Bf^T\BAs^T\Bb + \Bf^T\BAs^T\!\BA\Bf + \lambda^2 @@\Bf^T\Bf. \label{ME:eq:phi2}
\end{eqnarray}
Compute the gradient of $\phi$ with respect to the variable $\Bf$, recalling that $\lambda$ is a constant:
\begin{eqnarray} 
   \nabla \phi(\Bf) &=& \Bzero - 2 \BAs^T\Bb + 2@\BAs^T\!\BA\Bf + 2 \lambda^2 @@\Bf \nonumber \\[5pt]
                    &=& 2 (\BAs^T\BA\Bf + \lambda^2@@ \Bf - \BAs^T\Bb).  \label{ME:gradient}
\end{eqnarray}
To find critical points of $\phi$, set $\nabla \phi(\Bf) = \Bzero$ and solve for $\Bf$.  
Thus we find that $\BAs^T\BA\Bf + \lambda^2@@ \Bf = \BAs^T\Bb$, i.e.,
\begin{equation} \label{ME:regne}
   \big(\BAs^T\!\BA + \lambda^2\BI)\Bf = \BAs^T\Bb.
\end{equation}
We denote the solution for this fixed value of $\lambda$ as
\begin{equation} \label{ME:regsol}
   \Bf_\lambda := (\BAs^T\!\BA + \lambda^2\BI)^{-1}\BAs^T\Bb.
\end{equation}
Notice that if $\BA$ is invertible and we take $\lambda=0$, we get
\[ \Bf_{@0} \ =\  (\BAs^T\!\BA)^{-1}\BAs^T\Bb \ =\  \BAs^{-1} \BAs^{-T}\BAs^T \Bb \ =\ \BAs^{-1}\Bb,\]
and so this approach recovers the original solution~(\ref{ME:Ainv}).
When $\lambda>0$, the $\lambda^2\BI$ term in~(\ref{ME:regsol}) 
increases all the diagonal entries of $\BAs^T\!\BA+\lambda^2\BI$, making the rows and columns increasingly
distinct from one another.  
This enhancement adds robustness to the inversion process -- although when $\lambda$ 
is too large the $\lambda^2@\BI$ term will dominate the $\BAs^T\BA$ term, 
resulting in an $\Bf$ that makes $\|\Bb-\BA\Bf\|$ unacceptably large.

Perhaps it is tempting to implement the formula~(\ref{ME:regsol}) by 
computing the inverse of the $n\times n$ matrix $\BAs^T\BA+\lambda^2@\BI$,
and then multiplying the result against $\Bb$.
Instead, it is again more efficient 
to solve a linear system of equations using Gaussian elimination.

\begin{lstlisting}[language=Python]
lam  = 1e-5;                                      % regularization parameter
I    = np.identity(n)                             % n-by-n identity matrix
flam = np.linalg.solve(A.T@A+(lam**2)*I,(A.T)@b)  % solve linear system
\end{lstlisting}

For subtle numerical reasons, one actually prefers to set this problem
up a little bit differently.  As will be derived in an exercise at the
end of this section, one can pose the objective function~(\ref{ME:phi})
as a standard \emph{linear least squares problem} of the form
\begin{equation} \label{ME:LS}
 \min_{\Bf \in \R^n} 
 \phi(\Bf) = 
 \min_{\Bf \in \R^n} 
   \left\| \left[\begin{array}{c} \Bb \\ \Bzero\end{array}\right]
         - \left[\begin{array}{c} \BA \\ \lambda@\BI \end{array}\right]
             \Bf\; \right\|.
\end{equation}
We will denote the \emph{augmented} matrix and vector in this equation as
\[ \BA_\lambda = \left[\begin{array}{c} \BA \\ \lambda@\BI \end{array}\right] \in \R^{2n\times n}, 
   \qquad
   \Bb_\lambda = \left[\begin{array}{c} \Bb \\ \Bzero \end{array}\right] \in \R^{2n}.
\] 
In Python, one can solve the least squares problem~(\ref{ME:LS}) 
using the \verb|np.linalg.lstsq| command.
(Other languages provide similar functionality;
e.g., in MATLAB, \verb|flam = Alam\blam|.)

\begin{lstlisting}[language=Python]
lam  = 1e-5;                                      % regularization parameter
I    = np.identity(n)                             % n-by-n identity matrix
Alam = np.block([[A],[lam*I]])
blam = np.block([b,np.zeros(n)])
flam = np.linalg.lstsq(Alam,blam,rcond=None)[0]   % solve least squares problem
\end{lstlisting}

\begin{figure}[t!]
\begin{center}
 \includegraphics[height=0.97in]{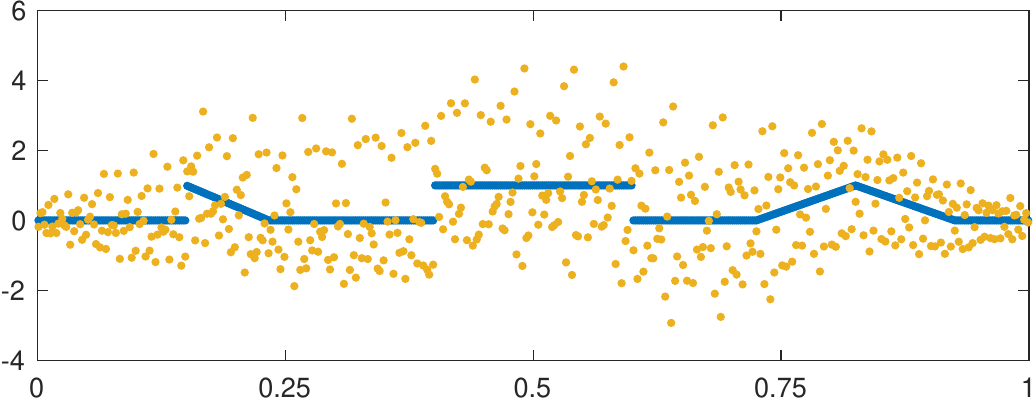}
 \begin{picture}(0,0)\put(-45.5,59){\footnotesize$\lambda=10^{-5}$}\end{picture}\quad
  \includegraphics[height=0.97in]{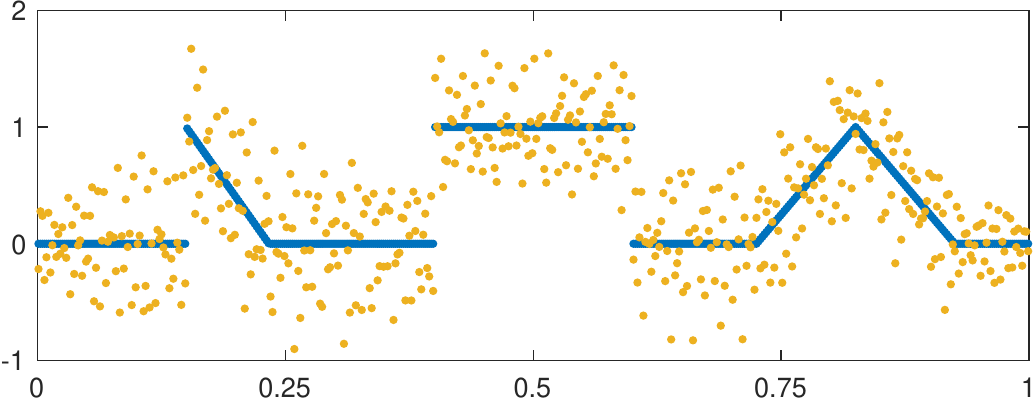}
 \begin{picture}(0,0)\put(-55,58){\footnotesize$\lambda=5\cdot10^{-5}$}\end{picture}

\medskip
 \includegraphics[height=0.97in]{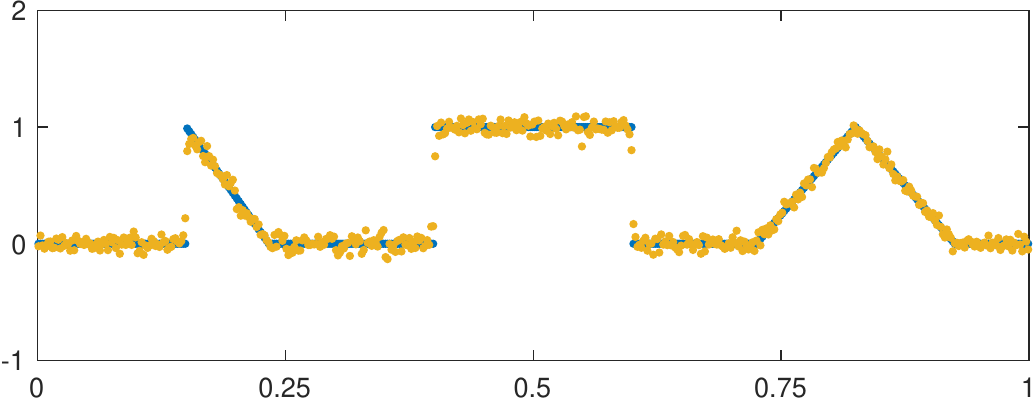}
 \begin{picture}(0,0)\put(-45,58){\footnotesize$\lambda=10^{-3}$}\end{picture}\quad
 \includegraphics[height=0.97in]{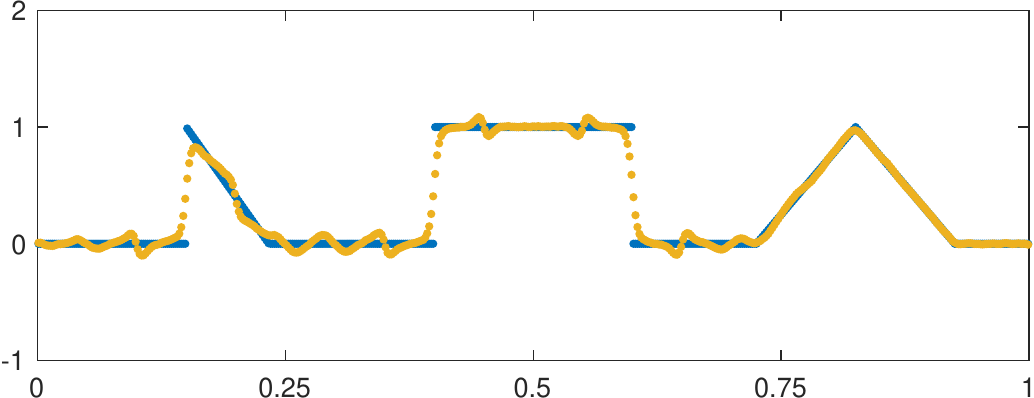}
 \begin{picture}(0,0)\put(-45,58){\footnotesize$\lambda=10^{-2}$}\end{picture}

\medskip
 \includegraphics[height=0.97in]{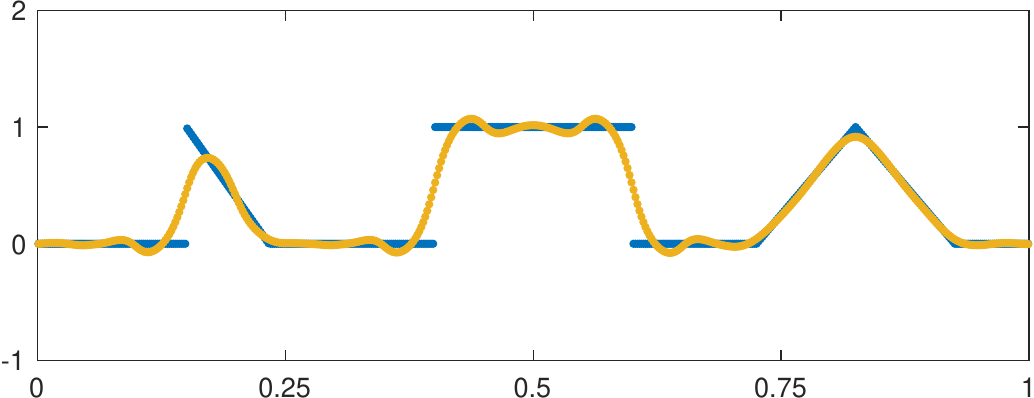}
 \begin{picture}(0,0)\put(-45,58){\footnotesize$\lambda=10^{-1}$}\end{picture}\quad
 \includegraphics[height=0.97in]{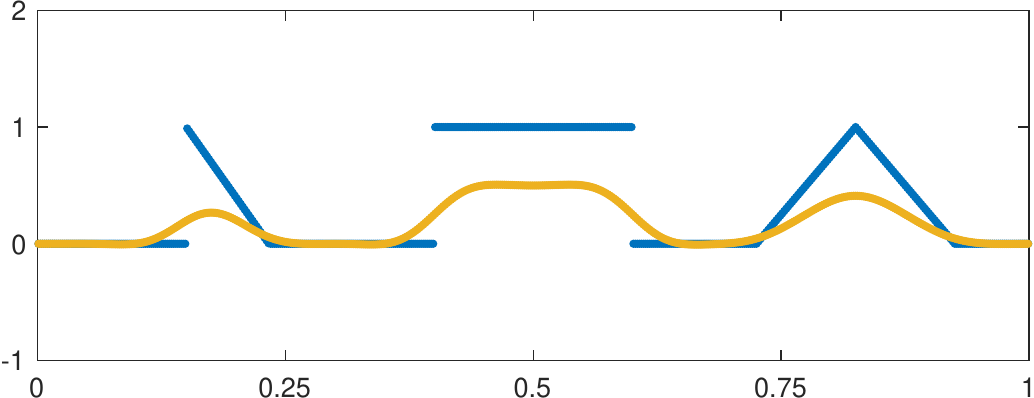}
 \begin{picture}(0,0)\put(-45,58){\footnotesize$\lambda=10^{0}$}\end{picture}
\end{center}

\vspace*{-1.5em}
\caption{\label{ME:fig:bw_reg}
Six attempts to recover $\Bf$ from $\bnoise$ polluted with random noise
(normally distributed, mean~zero, standard deviation $10^{-5}\|\Bb\|$),
using grid of $n=500$ points with the hat-function blurring kernel~(\ref{ME:hatker})
using $z=0.05$.
The blue points show the true vector $\Bf$ that we are trying to recover;
the yellow points show $\Bf_\lambda$ for the various $\lambda$ values.
When $\lambda$ is too small (top left), the regularization term
$\lambda^2\|\Bf\|^2$ has little effect, and the resulting $\Bf_\lambda$ is poor; 
when $\lambda$ is too large (bottom right),
$\Bf_\lambda$ is excessively suppressed.  The intermediate values
$\lambda=10^{-2}$ and $\lambda=10^{-1}$ do a better (but not perfect) job of recovering 
the true signal.}
\end{figure}
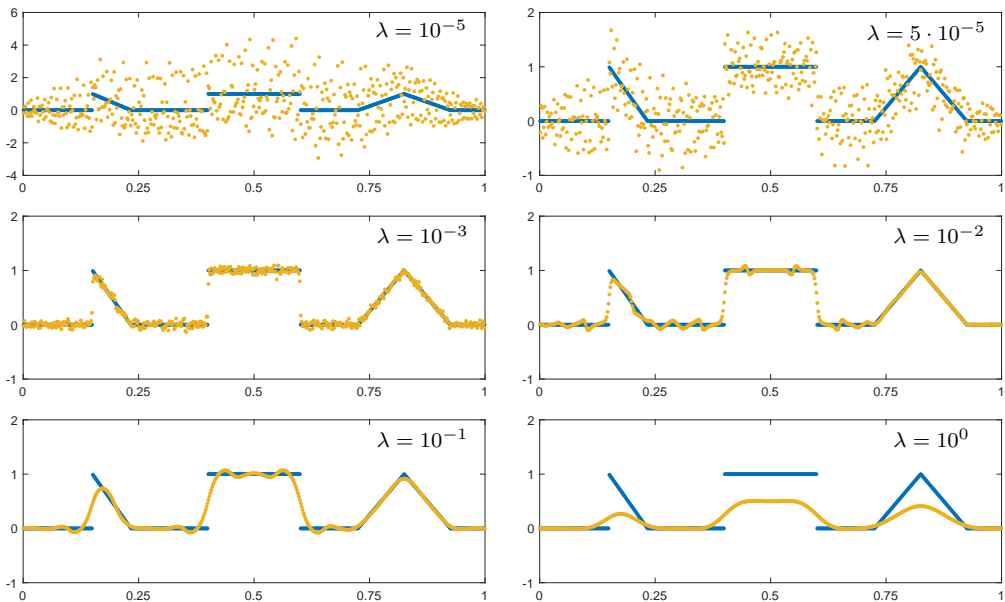

Figure~\ref{ME:fig:bw_reg} revisits the example on the left side of 
Figure~\ref{ME:fig:bw2}, where a whiff of noise (standard deviation 
$10^{-5}\|\Bb\|$) was enough to entirely break our initial recovery attempt,
$\frec = \BAs^{-1}\bnoise$.  Now we attempt to recover the signal using
regularization, $\Bf_\lambda = (\BAs^T\BA+\lambda^2\BI)^{-1} \bnoise$
for six different values of the regularization parameter $\lambda$.
When $\lambda$ is too small, regularization has little effect on $\Bf_\lambda$;
when $\lambda$ is too large, the solution is suppressed.  
Between these extremes the regularization gives an improved solution over the 
$\Bf = \BAs^{-1}\Bb$ solution we saw in Figure~\ref{ME:fig:bw2}.

\medskip

\begin{figure}[t!]
 \includegraphics[width=3.75in]{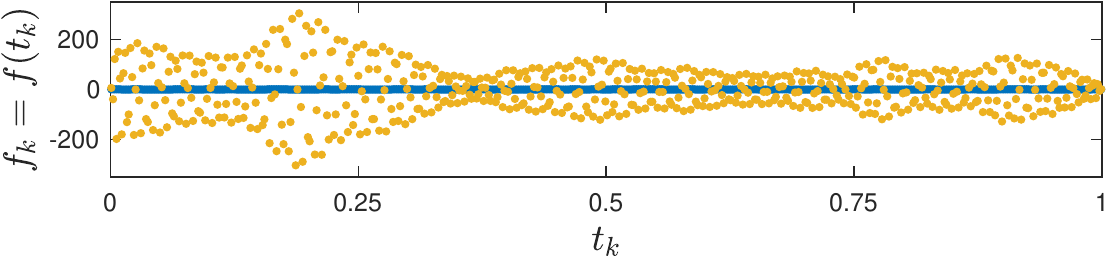} 
 \begin{picture}(0,0)
    \put(10,25){\includegraphics[width=1.25in]{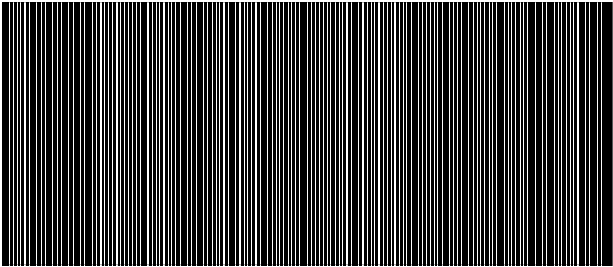}} 
    \put(37,10){\small \emph{$\lambda=10^{-10}$}}
 \end{picture}

\medskip
 \includegraphics[width=3.75in]{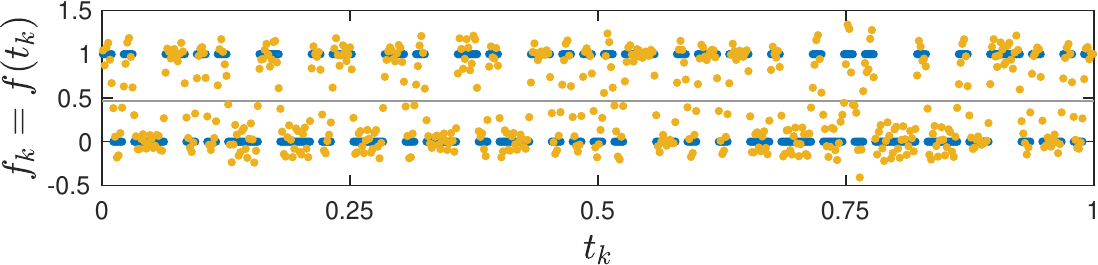} 
 \begin{picture}(0,0)
    \put(10,25){\includegraphics[width=1.25in]{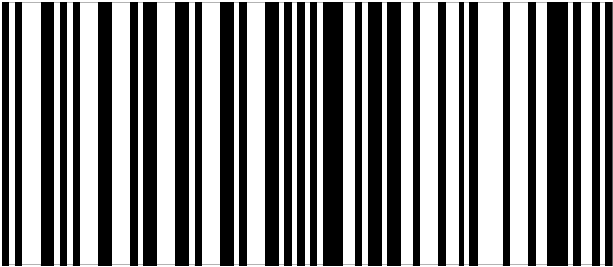}} 
    \put(37,10){\small \emph{$\lambda=10^{-7}$}}
 \end{picture}

\medskip
 \includegraphics[width=3.75in]{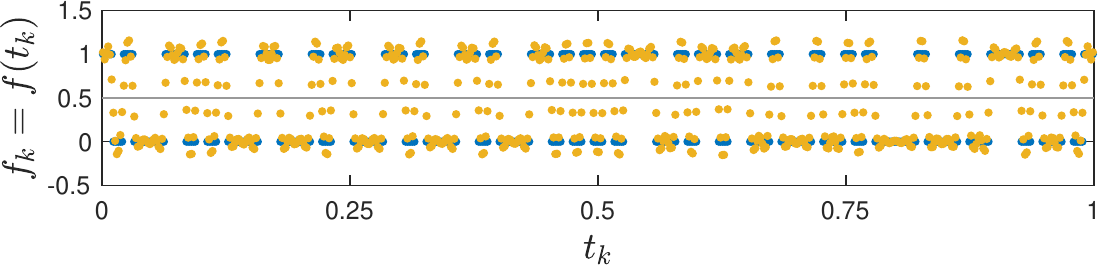} 
 \begin{picture}(0,0)
    \put(10,25){\includegraphics[width=1.25in]{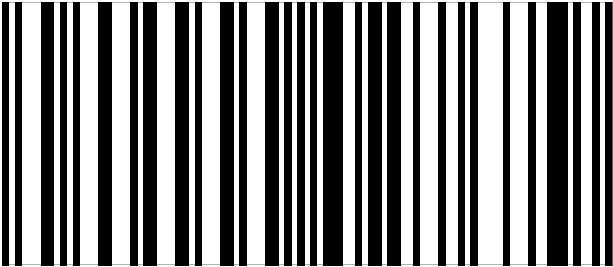}} 
    \put(37,10){\small \emph{$\lambda=10^{-5}$}}
 \end{picture}

\medskip
 \includegraphics[width=3.75in]{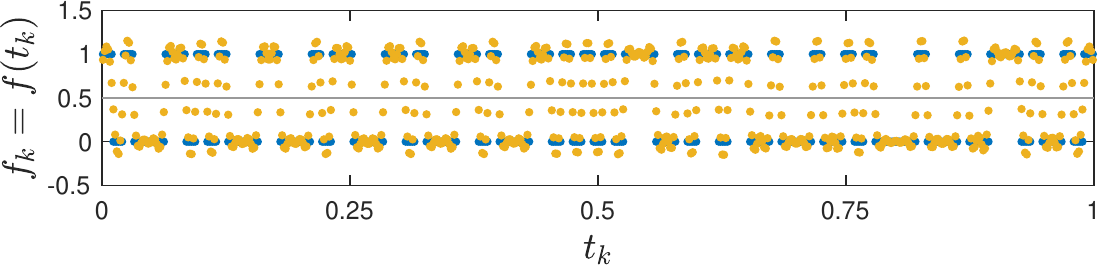} 
 \begin{picture}(0,0)
    \put(10,25){\includegraphics[width=1.25in]{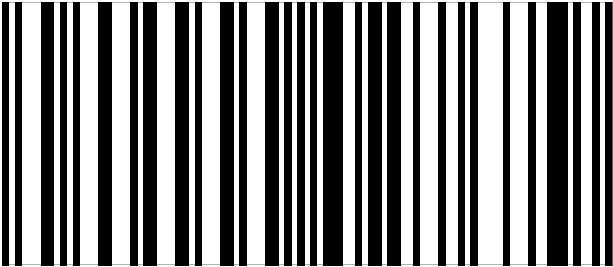}} 
    \put(37,10){\small \emph{$\lambda=10^{-3}$}}
 \end{picture}

\medskip
 \includegraphics[width=3.75in]{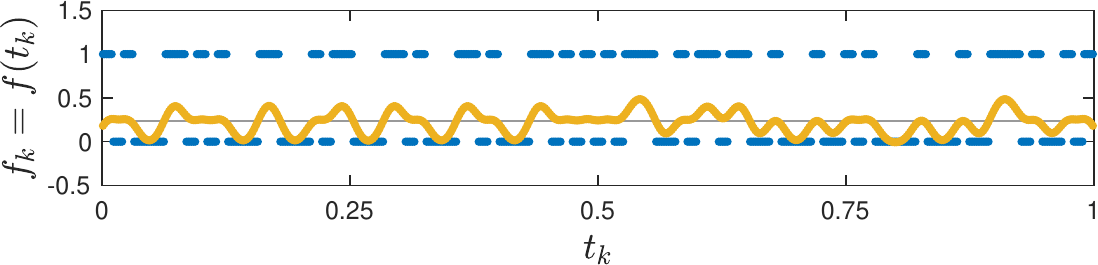} 
 \begin{picture}(0,0)
    \put(10,25){\includegraphics[width=1.25in]{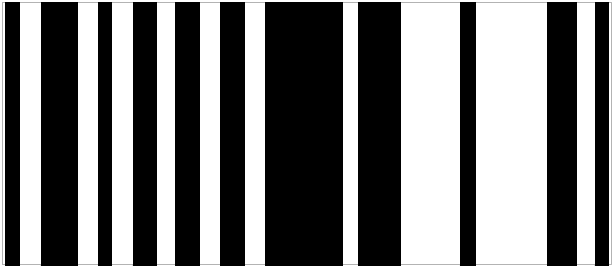}} 
    \put(37,10){\small \emph{$\lambda=10^{0}$}}
 \end{picture}

\vspace*{-0.5em}
\caption{\label{ME:fig:barcode1_reg}
Barcode deblurring with the help of regularization,
again with $n=570$ and the Gaussian kernel~(\ref{ME:gaussker}) with $z=0.01$.
We inject some random noise (size $10^{-8}\|\Bb\|$) to get $\bnoise$,
and use regularization to recover $\flam$ for five values of $\lambda$.
The yellow points show $\flam$; the blue points show the true values of $\Bf$.
In each case, we translate $\Bf_\lambda$ into a barcode by letting each entry 
above the midpoint of the extreme values of $\Bf_\lambda$ correspond to a black bar.
(This threshold is shown as the gray horizontal line visible on the bottom four plots.)
The value $\lambda=10^{-10}$ is too small to tame the unruly solution, 
and the barcode looks nothing like the correct version seen atop Figure~\ref{ME:fig:barcode1}.
The intermediate values $\lambda=10^{-7}$, $\lambda=10^{-5}$, and $\lambda=10^{-3}$ 
all do a reasonable job of recovering the desired barcode.
The large value $\lambda=10^0$ puts too much emphasis on minimizing $\|\Bf\|$
and not enough on making $\BA\Bf \approx \Bb$, 
resulting in a poor solution.}
\end{figure}
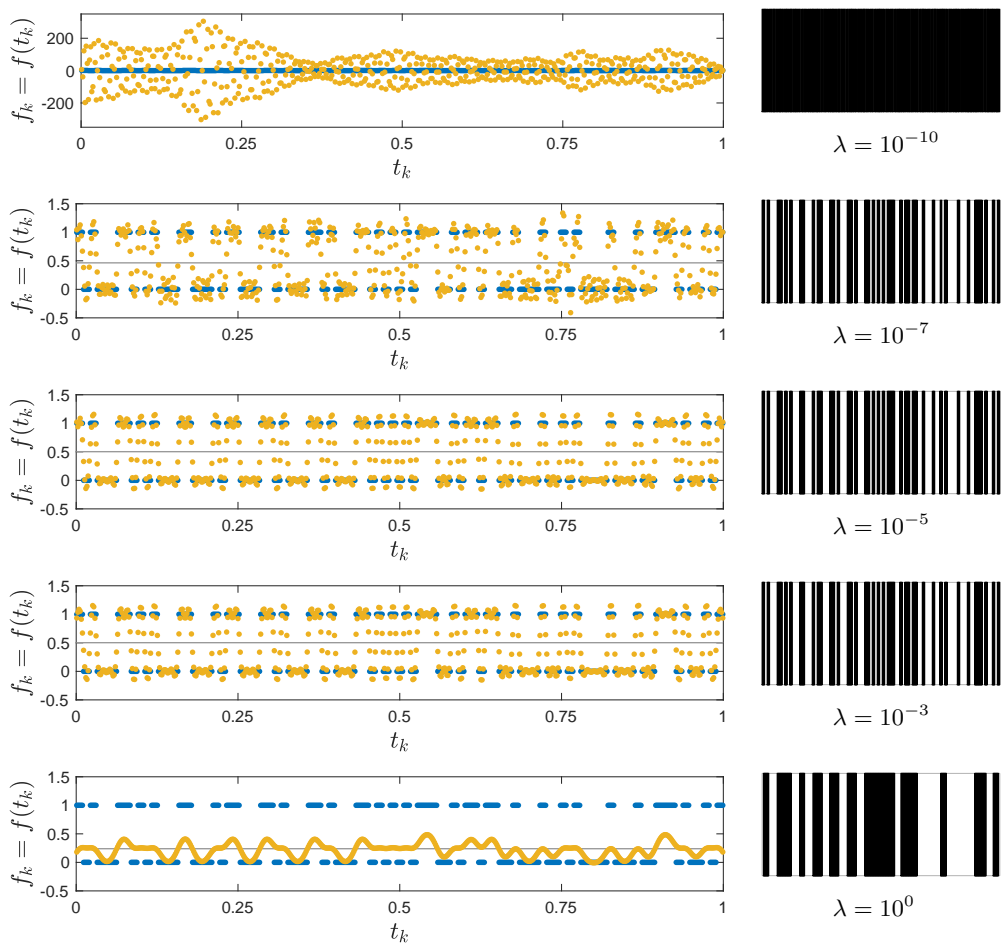

Figure~\ref{ME:fig:barcode1_reg} shows the regularized solutions for the 
barcode deblurring problem introduced in Figure~\ref{ME:fig:barcode1}, 
but now using $\bnoise$ with noise level $10^{-8}\|\Bb\|$.
Recall that in this case the true solution $\Bf\in\R^{570}$ is a binary vector
(all entries are either zero or one),
and so it suffices to recover a solution $\Bf_\lambda$ for which we can decide if
$f_j = 0$ or $f_j = 1$.  
For example, for each $\lambda$ we could set a threshold at the average of the extreme values
computed in the vector $\flam$,
\[ \tau_\lambda := \frac{1}{2}
                  \Big({\displaystyle \min_{1\le j\le n} (\flam)_j\ +\ \max_{1\le j\le n} (\flam)_j}\Big),\]
and then draw a vertical black bar at $t_j\in[0,1]$ if $(\flam)_j \ge \tau_\lambda$ for $j=1,\ldots, n$.
The barcodes on the right of Figure~\ref{ME:fig:barcode1_reg} are the result.
Again we see that taking $\lambda$ too small does not sufficiently regularize the solution,
and the resulting barcode is nonsense.
Taking $\lambda$ too large pushes the values in $\Bf_\lambda$ too close to zero, 
also giving a ridiculous solution.
The best choices for $\lambda$ lie in between these extremes.
While none of three values $\lambda = 10^{-7}$, $10^{-5}$, and $10^{-3}$ recover 
the function $f$ perfectly, using the $\tau_\lambda$ threshold allows us to recover
the bar code quite reliably in all three cases.
(The $10^{-7}$ result is off in just three of the $n=570$ positions; 
the $10^{-5}$ and $10^{-3}$ results match perfectly.)

When we compare the two deblurring examples shown in Figures~\ref{ME:fig:bw_reg} 
and~\ref{ME:fig:barcode1_reg}, we notice that
several of the ``good'' values of $\lambda$ for the barcode example 
(Figure~\ref{ME:fig:barcode1_reg}) would would have been much too small 
for the problem shown in Figure~\ref{ME:fig:bw_reg}.
Clearly there is some nuance to the choice of $\lambda$, and we should not expect one
magical value of $\lambda$ to work for all problems we encounter.
What guidance can we find for selecting this crucial parameter?
We will explore this question in the next section.

\begin{reflections}
\item Explain, in careful detail, how to move from~(\ref{ME:eq:phi1}) 
      to~(\ref{ME:eq:phi2}) in the derivation of the formula for $\phi(\Bf)$.

\item When solving a modeling problem, it is often easy to fixate on 
      minimizing an objective function, like $\|\Bb - \BA\Bf\|$,
      without thinking about what that means for the motivating application.
      For the barcode problem, which goal is more important:  getting an
      exact solution $\Bf$ to $\BA\Bf=\Bb$, or getting an approximate
      solution $\Bf$ that is good enough for you to recover the barcode 
      exactly?
\end{reflections}

\noindent
Now would be a good time to explore Exercises~\ref{ME:ex:xy}--\ref{ME:ex:reg2x2} 
starting on page~\pageref{ME:ex:xy}.

\begin{figure}[b!]
\begin{center}
\includegraphics[width=3.5in]{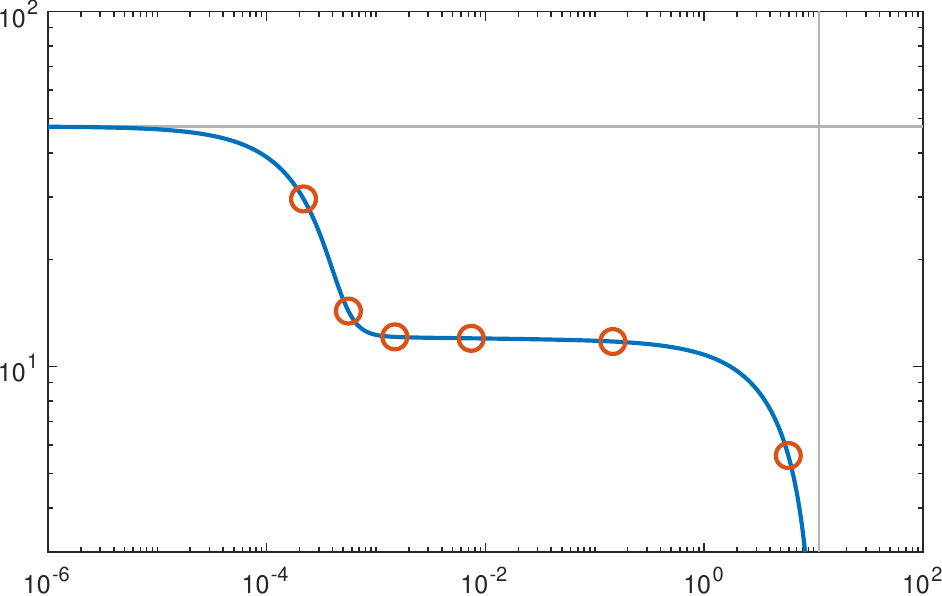}
\begin{picture}(0,0)
\put(-273,100){{$\|\Bf_\lambda\|$}}
\put(-160,-10){{$\|\bnoise - \BA\Bf_\lambda\|$}}
\put(-53,163){\rotatebox{0}{\footnotesize{$\|\bnoise\|$}}}
\put(-5,124){\rotatebox{0}{\footnotesize{$\|\frec\|$}}}
\put(-213,103){\rotatebox{0}{\footnotesize{$\lambda=10^{-5}$}}}
\put(-161,80){\rotatebox{30}{\footnotesize{$\lambda=5\cdot10^{-5}$}}}
\put(-151,74.5){\rotatebox{30}{\footnotesize{$\lambda=10^{-3}$}}}
\put(-132,74.25){\rotatebox{30}{\footnotesize{$\lambda=10^{-2}$}}}
\put(-94,73.5){\rotatebox{30}{\footnotesize{$\lambda=10^{-1}$}}}
\put(-77,36){\rotatebox{0}{\footnotesize{$\lambda=10^{0}$}}}
\end{picture}
\end{center}
\caption{\label{ME:fig:bw_L}
An ``L~curve'' corresponding to the example in Figure~\ref{ME:fig:bw_reg}.
The goal of regularization is to choose a value of $\lambda$ that balances the
goals of minimizing $\|\bnoise-\BA\Bf_\lambda\|$ while controlling $\|\Bf_\lambda\|$.
Appealing values of $\lambda$ typically fall near the central bend in this curve 
(near $\lambda=10^{-3}$ on this plot).
}
\end{figure}
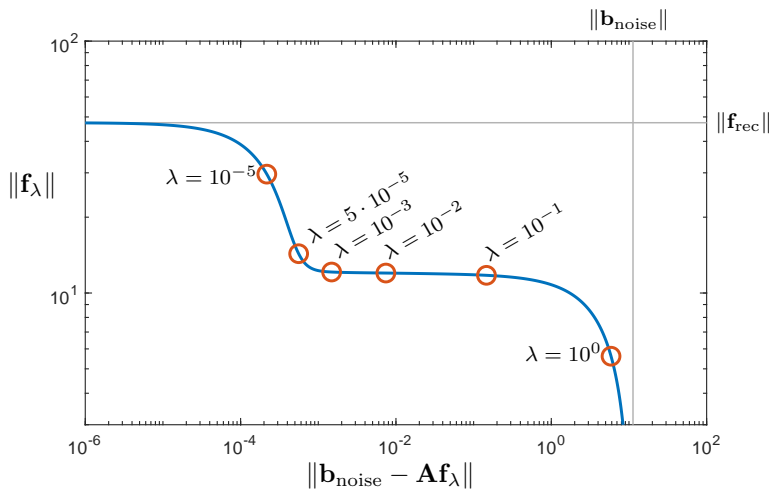

\section{Choosing the regularization parameter}

As suggested by Figures~\ref{ME:fig:bw_reg} and~\ref{ME:fig:barcode1_reg}, 
selection of the regularization parameter can be somewhat subtle.  
The optimal choice will depend on the application and the data.
Are you willing to accept a solution that exhibits a bit of chatter, 
but captures edges well 
(as for $\lambda=10^{-3}$ in Figure~\ref{ME:fig:bw_reg}), 
or do you want to suppress the noise at the cost of an overly smooth solution 
(as for $\lambda=10^{-2}$ or $\lambda=10^{-3}$ in Figure~\ref{ME:fig:bw_reg})?
While various sophisticated approaches exist to guide the selection of $\lambda$
(see~\cite[sect.~6.4.1]{GV12} and \cite{Han10} for some details), 
in this manuscript we will only explore trial-and-error investigations.

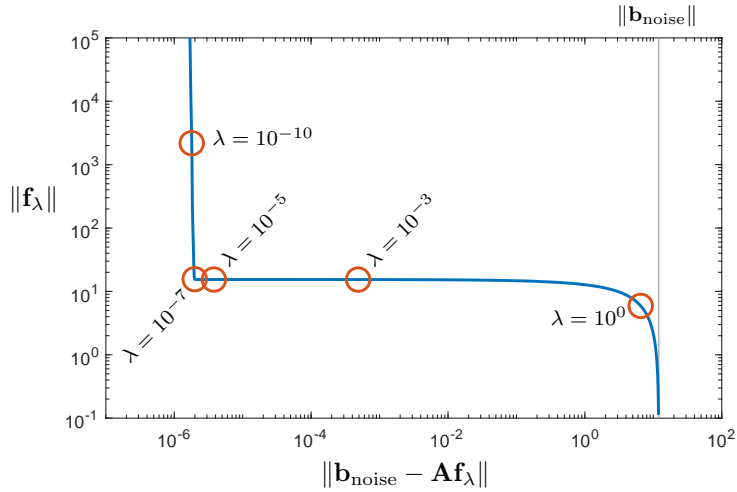
\begin{figure}[t!]
\begin{center}
\vspace*{5pt}
\includegraphics[width=3.5in]{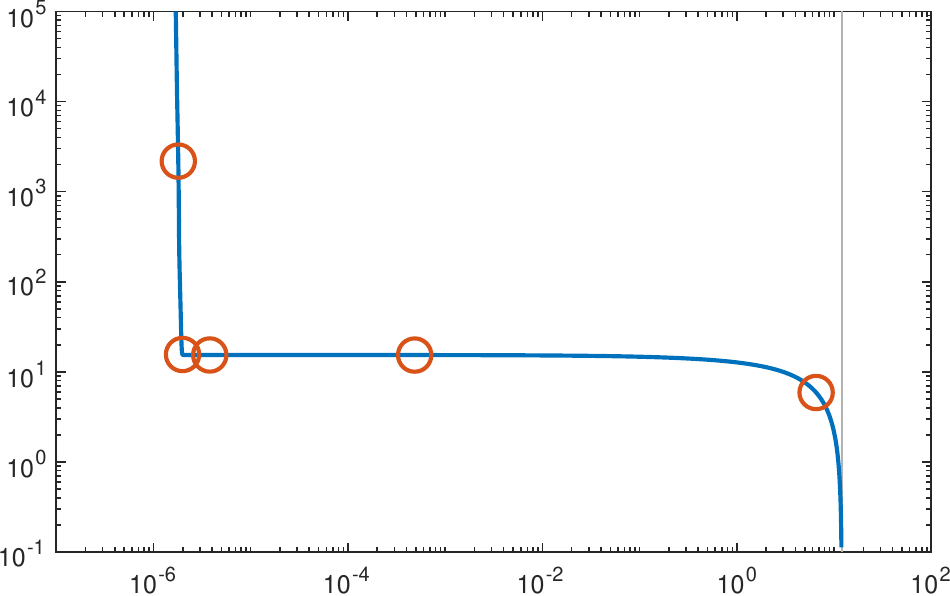}
\begin{picture}(0,0)
\put(-278,92){{$\|\Bf_\lambda\|$}}
\put(-160,-12){{$\|\bnoise - \BA\Bf_\lambda\|$}}
\put(-48,162){\rotatebox{0}{\footnotesize{$\|\bnoise\|$}}}
\put(-201,115){\rotatebox{0}{\footnotesize{$\lambda=10^{-10}$}}}
\put(-237,32){\rotatebox{45}{\footnotesize{$\lambda=10^{-7}$}}}
\put(-199.5,68.5){\rotatebox{45}{\footnotesize{$\lambda=10^{-5}$}}}
\put(-144.5,68.5){\rotatebox{45}{\footnotesize{$\lambda=10^{-3}$}}}
\put(-73,47){\rotatebox{0}{\footnotesize{$\lambda=10^{0}$}}}
\end{picture}
\end{center}
\vspace*{5pt}
\caption{\label{ME:fig:barcode1_L}
An ``L~curve'' corresponding to the example in Figure~\ref{ME:fig:barcode1_reg}.
The three intermediate values of $\lambda$ that lead to a reasonable recovery of the barcode
($\lambda = 10^{-7}, 10^{-5}, 10^{-3}$) all fall on the plateau of the L~shape, 
with $10^{-7}$ very close to the corner.
This plot makes clear that $\lambda=10^{-10}$ fails to control $\|\Bf_\lambda\|$,
while $\lambda=10^0$ permits $\|\bnoise-\BA\flam\|$ to be too large.}
\end{figure}

How does the balance between $\|\Bb-\BA\Bf_\lambda\|$ and $\|\Bf_\lambda\|$ tilt, as 
$\lambda$ increases from $\lambda=0$ to $\lambda\to\infty$?  
To see how this relationship evolves, 
you can solve the regularization problem for a wide range of $\lambda$ values, 
and then create a plot of the ordered pairs $(\|\Bb-\BA\Bf_\lambda\|, \|\Bf_\lambda\|)$.  
Because these quantities can vary over several orders of magnitude, 
it is helpful to draw this as a log-log plot
(i.e., plotting the pairs $(\log_{10} \|\Bb-\BA\Bf_\lambda\|, \log_{10} \|\Bf_\lambda\|)$.  
Figure~\ref{ME:fig:bw_L} shows a representative plot generated for 
the problem shown in Figure~\ref{ME:fig:bw_reg}, while 
Figure~\ref{ME:fig:barcode1_L} shows the analogous plot for the
barcode regularization problem in Figure~\ref{ME:fig:barcode1_reg}.

Such plots are called ``L~curves,''  named for the prominent central bend 
that these plots often (but not always) exhibit.  
As $\lambda$ increases from $0\to\infty$, the plot moves from the top-left to the bottom-right.
The ``L'' shows the transition from a vertical drop (where small increases in $\lambda$ 
give big, desirable reductions in $\|\Bf_\lambda\|$ without much increase in the misfit 
$\|\Bb-\BA\Bf_\lambda\|$) toward a horizontal plateau (where further increases in $\lambda$ 
give modest reductions in $\|\Bf_\lambda\|$ but undesirably large increases 
in $\|\Bb-\BA\Bf_\lambda\|$).

When you first encounter an L~curve, the plot can be a little tricky to understand.  
First off, appreciate that you cannot directly read off the values for $\lambda$ 
from the horizontal or vertical axes.  
You need to plot some test values of $(\|\Bb-\BA\Bf_\lambda\|, \|\Bf_\lambda\|)$
to appreciate how $\lambda$ changes along the curve.  
(The red circles in Figures~\ref{ME:fig:bw_L} and~\ref{ME:fig:barcode1_L}
show the points on the L~curve corresponding to the $\lambda$ values used 
to compute the solutions shown in Figures~\ref{ME:fig:bw_reg} and~\ref{ME:fig:barcode1_reg}.)
Now look for several landmark features on the L~curve.
\begin{itemize} \setlength{\itemsep}{0pt} \setlength{\parskip}{0pt}
\item
As $\lambda\to0$, the optimization places diminishing emphasis on $\|\Bf_\lambda\|$;
it primarily minimizes $\|\Bb-\BA\Bf_\lambda\|$.  
For $\lambda=0$ the minimum is attained by $\frec = \BAs^{-1}\Bb$,
which gives $\Bb-\BA\frec = \Bzero$.
Look for the horizontal asymptote at the level $\|\frec\|$ at the top-left
of the L~curve, visible in Figure~\ref{ME:fig:bw_L}.

\medskip
\item
As $\lambda\to\infty$, the optimization completely focuses on minimizing 
$\|\Bf_\lambda\|$, and so $\Bf_\lambda\to\Bzero$ as $\lambda\to\infty$.
Thus, as $\lambda$ gets large $\|\Bb-\BA\Bf_\lambda\| \to \|\Bb\|$.
This limit is seen on the bottom-right of Figures~\ref{ME:fig:bw_L} 
and~\ref{ME:fig:barcode1_L}.
\end{itemize}

Sitting between these extremes, the corner of the L is a ``sweet spot'' 
where the goals of minimizing both $\|\Bf_\lambda\|$ 
and $\|\Bb-\BA\Bf_\lambda\|$ are balanced.
Due to the logarithmic nature of the axes, 
choosing $\lambda$ even a bit before the corner 
(i.e., too small) can yield a poor solution due to the large value of $\|\flam\|$,
whereas choosing $\lambda$ a bit after the corner
increases $\|\Bb-\BA\flam\|$, but often not in a problematic way.
Consider the barcode example from Figure~\ref{ME:fig:barcode1_L}, 
where we know the entries of the true solution $\Bf$ must be~0 or~1.  
In this case, a factor of~10 difference in $\|\flam\|$ 
(say, $10^1$ versus $10^2$) has a significant impact on the 
quality of the recovered solution; however, a factor of~10 difference 
in the small values of $\|\Bb-\BA\flam\|$ (say, from $10^{-5}$ to $10^{-4}$) 
will hardly be noticeable in plots of the solution.

The following code will produce an L~curve like the one shown in
Figure~\ref{ME:fig:bw_L} (presuming you have \verb|A| and \verb|bnoise|
defined already). Note that code may take a few seconds to run; 
while it is straightforward, this approach is not a very expedient 
way to draw L curves.
(When solving regularization problems for many values of $\lambda$,
it can be more efficient to use the singular value decomposition,
as outlined in Section~\ref{ME:sec:SVD}, to solve these systems,
instead of solving independent least squares problems at each step.

\newpage

\begin{lstlisting}[language=Python]
lamvec = np.logspace(-7,.5,100);        # 100 (log-spaced) lambda values

I    = np.identity(n)                   # identity matrix
blam = np.block([bnoise,np.zeros(n)])   # b_lambda vector
fnorm = [];                             # vector for ||f_lam|| for each lambda
rnorm = [];                             # vector for ||b_noise-A*f_lam||

for lam in lamvec:                      # loop over all lambda values
    Alam = np.block([[A],[lam*I]])                  # form A_lam
    flam = np.linalg.lstsq(Alam,blam,rcond=None)[0] # solve for f_lam
    fnorm.append(np.linalg.norm(flam))              # store ||f_lam||
    rnorm.append(np.linalg.norm(bnoise-A@flam))     # store ||b_noise-A*f_lam||

plt.loglog(rnorm,fnorm,'-')             # produce the L curve plot
plt.xlim([1e-6, 1e2])
plt.ylim([3, 1e2])
plt.xlabel('$||b_{noise}-A f_\lambda||$',fontsize=14)
plt.ylabel('$||f_\lambda||$',fontsize=14)
\end{lstlisting}

\begin{reflection}
\item Explain why choosing $\lambda$ so that the point 
      $(\|\bnoise-\BA\flam\|,\|\flam\|)$ occurs a bit \emph{after} 
      the corner of the L-curve is often less problematic 
      than choosing $\lambda$ so that 
      $(\|\bnoise-\BA\flam\|,\|\flam\|)$ 
      occurs a bit \emph{before} the corner.
\end{reflection}

\noindent
Now would be a good time to explore Exercise~\ref{ME:ex:Lcurve} 
on page~\pageref{ME:ex:Lcurve}.

\section{UPC barcodes} \label{ME:sec:UPC}

We have shown an example of barcode deblurring in Figures~\ref{ME:fig:barcode1}
and~\ref{ME:fig:barcode1_L}.
Now it is time for you to try your hand at decoding several such images to
discover the product encoded by a blurry ``scan'' of a UPC symbol.
To decode these symbols, you need to know a little bit about
the UPC barcode system itself.  
We have digested the summary here from the Universal Product Code 
Wikipedia page~\cite{UPCwiki}.
Those wishing to dive deeper can consult the authoritative 
\emph{GS1 General Standard Specification}~\cite[pp.~253--256]{GS1}.
(What we describe here is more specifically referred to as the ``UPC-A'' format.)
Barcodes are themselves a fascinating technology.
For some information about their history, we recommend the podcast~\cite{Mar14} 
and the collection~\cite{Hab01};
additional information on the UPC decoding problem can be found in, e.g., 
\cite{ISW13,SG22,Wit04}.

\begin{table}[t!]
\caption{\label{ME:tbl:upc} Key for decoding UPC symbols, adapted from~\cite{UPCwiki}:
{\sf b} and {\sf w} denote black and white bars.}
\begin{center}
\begin{tabular}{lll}
\multicolumn{1}{c}{\textsl{colors}} & 
\multicolumn{1}{c}{\textsl{number of bars}} & 
\multicolumn{1}{c}{\textsl{description}} \\[2pt]
$\sf b@w@b$   & three bars, each of width~1 & \textsl{start code}  \\
$\sf w@b@w@b$  & four bars of total width 7 & first digit \\
$\sf w@b@w@b$  & four bars of total width 7 & second digit \\
$\sf w@b@w@b$  & four bars of total width 7 & third digit \\
$\sf w@b@w@b$  & four bars of total width 7 & fourth digit \\
$\sf w@b@w@b$  & four bars of total width 7 & fifth digit \\
$\sf w@b@w@b$  & four bars of total width 7 & sixth digit \\
$\sf w@b@w@b@w$ & five bars, each of width~1        & \textsl{middle code}\\
$\sf b@w@b@w$  & four bars of total width 7 & seventh digit \\
$\sf b@w@b@w$  & four bars of total width 7 & eighth digit \\
$\sf b@w@b@w$  & four bars of total width 7 & ninth digit \\
$\sf b@w@b@w$  & four bars of total width 7 & tenth digit \\
$\sf b@w@b@w$  & four bars of total width 7 & eleventh digit \\
$\sf b@w@b@w$  & four bars of total width 7 & twelfth digit \\
$\sf b@w@b$    & three bars, each of width~1 & \textsl{end code} 
\end{tabular}
\end{center}
\end{table}

All UPC barcodes contain 59 alternating black and white bars of varying widths;
together, these bars encode 12 decimal digits.
Each bar, be it black ({\sf b}) or white ({\sf w}), has one of four widths: 
either 1, 2, 3, or 4 units wide.  
Table~\ref{ME:tbl:upc} shows how the bars are grouped to define the digits that
make up the 12-digit code.
Each of these digits is encoded with four bars (two white, two black); 
the widths of these bars will differ, but the \emph{sum} of the width of
the four bars is always~7, so that all UPCs have the same width
(95~units wide).  Each digit (0--9) corresponds to a particular pattern
of bar widths, according to Table~\ref{ME:tbl:upc_digits}.  
(Two different mirror-image patterns are used for each digit.)

\begin{table}[h!]
\caption{\label{ME:tbl:upc_digits} Correspondence of digits to patterns in UPC codes,
adapted from~\cite{UPCwiki}.}
\begin{center} 
\hspace*{-3em}\begin{tabular}{ccc}
\textsl{digit} & \textsl{pattern 1} & \textsl{pattern 2} \\
0 & 3--2--1--1 &  1--1--2--3 \\
1 & 2--2--2--1 &  1--2--2--2 \\
2 & 2--1--2--2 &  2--2--1--2 \\
3 & 1--4--1--1 &  1--1--4--1 \\
4 & 1--1--3--2 &  2--3--1--1 \\
5 & 1--2--3--1 &  1--3--2--1 \\
6 & 1--1--1--4 &  4--1--1--1 \\
7 & 1--3--1--2 &  2--1--3--1 \\
8 & 1--2--1--3 &  3--1--2--1 \\
9 & 3--1--1--2 &  2--1--1--3 
\end{tabular}\end{center}
\end{table}

For example, if the four bars (read left-to-right) have the pattern \\
\centerline{width~3 -- width~2 -- width~1 -- width~1}\\
or \\
\centerline{width~1 -- width~1 -- width~2 -- width~3,}\\
the corresponding UPC digit is~0.

Let us step through the example shown in Figure~\ref{ME:fig:barcode1}.
To check that you understand the system, try decoding the UPC below 
(for a can of Coke).  
To start you out, we give the code for the first
three digits.  In some cases it might be tricky to judge some of
the widths: you can appreciate the accuracy of optical scanners!

\begin{center}
\label{ME:cokeUPC}
\includegraphics[width=2.25in]{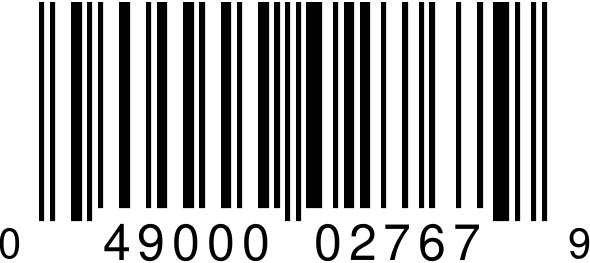}
\end{center}

\begin{center}
\begin{tabular}{r|ccc|cccc|cccc|cccc}
\textsl{color} & {\sf b} & {\sf w} & {\sf b} & {\sf w} & {\sf b} & {\sf w} & {\sf b}  & {\sf w} & {\sf b} & {\sf w} & {\sf b} & {\sf w} & {\sf b} & {\sf w} & {\sf b}\\
\textsl{width} &  1      &    1    &   1 & 3 & 2 & 1 & 1 & 1 & 1 & 3 & 2 & 3 & 1 & 1 & 2 \\
\textsl{meaning} &  \multicolumn{3}{c}{start code} &  \multicolumn{4}{|c}{0} & \multicolumn{4}{|c}{4} &  \multicolumn{4}{|c}{9}
\end{tabular}
\end{center}

We want to simulate the reading of this UPC by an optical scanner,
e.g., in a supermarket check-out line.  
The barcode is represented mathematically as
a function $f(t)$ that takes the values zero and one:
zero corresponds to white ({\sf w}) bars, 
one corresponds to black ({\sf b}) bars.  
The function corresponding to the Coke bar code is shown below.

\vspace*{1em}
\begin{center}\includegraphics[scale=.65]{IMAGES/barcode1_ftrue}\end{center}

We now need to explain how we discretize the interval $[0,1]$ into $n$ points.
Recall that a UPC contains a total of 59 bars that span 95~\emph{units} 
(i.e., the sum of the widths of all the black and white bars together is~95).  
We choose to use $n=570$, so that we have 6~discretization points \emph{per each unit} 
of the barcode: hence a bar of width~1 will correspond to 6~consecutive entries
in our vector $\Bf$.  Since all UPCs begin with the start code of three bars of
width one ({\sf bwb}), when we use $n=570$ the true vector $\Bf$ that discretizes 
the function $f$ will always begin like
\[\Bf = [\;\underbrace{1,\; 1,\; 1,\; 1,\; 1,\; 1}_{\sf b},\; 
           \underbrace{0,\; 0,\; 0,\; 0,\; 0,\; 0}_{\sf w},\;
           \underbrace{1,\; 1,\; 1,\; 1,\; 1,\; 1}_{\sf b},\; \ldots\;]^T.\]

Suppose the optical scanner can only acquire $\Bb\in\R^n$, 
samples of $f$ blurred by the matrix $\BA$.
Let us assume that this blur can be described by the Gaussian kernel~(\ref{ME:gaussker})
with parameter $z=0.01$.  
Applying the discrete blurring process~(\ref{ME:blurmat})
to the true barcode~$f$ for the Coke UPC gives the samples shown below.
From this blurred function, it would be difficult to determine the widths
of the bars, and hence to interpret the barcode.  
We shall try to improve the situation by solving the inverse problem 
$\BA\Bf = \Bb$ for the vector $\Bf$ that samples the function $f(t)$ at
the points $t_k = (k-1/2)/n$.

\begin{center}\includegraphics[scale=.65]{IMAGES/barcode1_b}\end{center}

However, we are unlikely to measure the exact blurred vector, $\Bb$.
Instead, we model some pollution from random noise.
The plot below shows $\Bb$ plus a vector with normally distributed
random entries having mean zero and standard deviation $5 \times 10^{-3} \|\Bb\|$.
(We will use a bit less noise in our initial computations; 
we use higher noise here to make it easier to see in the plot.)

\begin{center}\includegraphics[scale=.65]{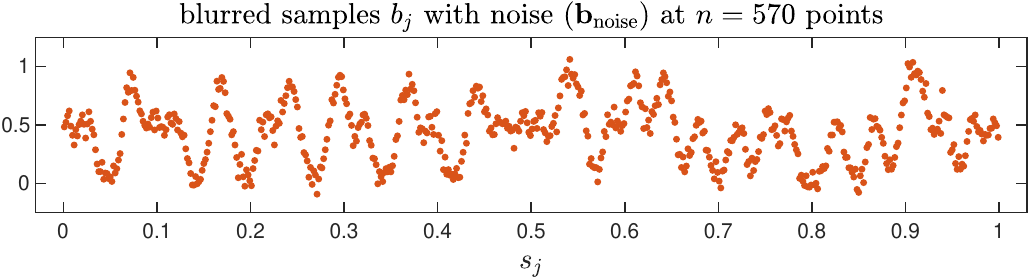}\end{center}

The Jupyter notebook \verb|coke_upc.ipynb| posted at 

\begin{center}
\verb|https://github.com/markembree/deblurring|
\end{center}

defines a function \verb|coke_upc| with the interface 
\begin{lstlisting}[language=Python]
A, b, bnoise, ftrue = coke_upc()
\end{lstlisting}
This code generates, for $n=570$: the blurring matrix 
$\BA\in\R^{570\times 570}$ for the kernel described above;
the blurred function sampled at 570~points, $\Bb\in\R^{570}$;
that same vector but with the addition of random noise of standard deviation 
$10^{-3} \|\Bb\|$, yielding $\Bb_{\rm noise}\in\R^{570}$;
the exact bar code solution $\Bf_{\rm true}$.  (We generated $\Bb$ 
as \verb|b = A@ftrue|,  then polluted it with random noise.
Since different noise is generated each time you call the routine, you
could get slightly different answers, but the qualitative results will
be the same.)

\medskip
\textsl{Some of the following exercises are intentionally open-ended.  
You will be evaluated on the thoroughness of your experiments, rather 
than for recovering a particular value for the barcodes.  
Include numerous plots and label what they show; explain (in words) what you discover
from each plot.  Your explanations are a key part of these exercises.}

\begin{reflections}
\item To demonstrate your understanding of the UPC barcode system,
      finish decoding the Coke UPC that we began on page~\pageref{ME:cokeUPC}.  
      Explain how the bars lead to the values {\sf 049000027679}.
      (Follow the lead given in the text for decoding the initial {\sf 049} section.)

\item Find a UPC barcode for some product that you have on hand.
      Decode the bars following the description in this section.
      (Most UPC barcodes include their numerical values directly underneath
      the bars, allowing you to check that you are correct.)

\item In recent years, two-dimensional ``QR codes'' have become a popular 
      way to transmit information.  Perform some background research on 
      how QR codes work, and their differences from one-dimensional UPC codes.
\end{reflections}

\noindent
Now would be a good time to explore Exercises~\ref{ME:ex:geupc}--\ref{ME:ex:mystery2} 
starting on page~\pageref{ME:ex:geupc}.

\section{Insight from the Singular Value Decomposition (SVD)} \label{ME:sec:SVD}

This optional concluding section offers some additional insight
for students who are familiar with the Singular Value Decomposition (SVD).\ \ 
For more information about the SVD oriented toward an undergraduate audience, 
see, e.g., \cite[chap.~5]{Emb3606}, \cite[chap.~7]{Str23}.
Let us begin here with a little background, to establish notation.

Let $\BA \in \Rmn$ be a matrix with $m$ rows and $n$ columns, and assume $m \ge n$.
The (economy-sized) SVD expresses $\BA$ as the product of three matrices,
\begin{equation} \label{ME:SVD}
   \BA = \BU\BSigma\BV^T,
\end{equation}
where 
\begin{itemize}
\item $\BU\in\R^{m\times n}$ has $n$ orthonormal columns, so $\BU^T\BU=\BI\in\R^{n\times n}$;
\item $\BSigma \in \R^{n\times n}$ is a diagonal matrix of \emph{singular values},
      sorted in decreasing order:
\[ \sigma_1\ge \sigma_2 \ge \cdots \ge \sigma_n\ge 0;\]
\item $\BV\in\R^{n \times n}$ has $n$ orthonormal columns, so $\BV^T\BV = \BI\in\R^{n\times n}$.
\end{itemize}
We can write the matrices out as
\[ \BU = \left[ \begin{array}{cccc} \rule[-8pt]{0pt}{20pt}\Bu_1 & \Bu_2 & \cdots & \Bu_n \end{array}\right] 
           \in \R^{m\times n},\qquad
   \BV = \left[ \begin{array}{cccc} \rule[-8pt]{0pt}{20pt}\Bv_1 & \Bv_2 & \cdots & \Bv_n \end{array}\right] 
           \in \R^{n\times n},
\]
and
\[ 
\BSigma = \left[ \begin{array}{cccc} \sigma_1 \\ & \sigma_2 \\ & & \ddots \\ & & & \sigma_n\end{array}\right] 
          \in \R^{n\times n}.
\]
We call the $\Bu_j$ \emph{left singular vectors} and the $\Bv_j$ \emph{right singular vectors}.
For each $j=1,\ldots, n$, 
\begin{equation} \label{ME:eq:Avsigu}
 \BA\Bv_j = \sigma_j@\Bu_j.
\end{equation}
The number of nonzero singular values is the $\emph{rank}$ of $\BA$, which we will denote by $r$.
In Python, you can compute the economy-sized SVD via the command
\begin{lstlisting}[language=Python]
U, S, Vt = np.linalg.svd(A, FullMatrices=False)
\end{lstlisting}
(Note that this command returns $\BV^T$ in the variable \verb|Vt|.
In contrast, the analogous command in MATLAB, 
\verb|[U,S,V] = svd(A,'econ')|, returns $\BV$ in the variable \verb|V|.)

Using properties of matrix-matrix multiplication, the common matrix form 
of the SVD~(\ref{ME:SVD}) can also be expressed in the \emph{dyadic form}
\begin{equation} \label{ME:dyadic}
 \BA = \sum_{j=1}^n \sigma_j^{} \Bu_j^{} \Bv_j^T.
\end{equation}
Take a moment to think about the elements of this formula:  $\sigma_j$ is just a number,
and thus acts like a scaling factor; since $\Bu_j \in \R^m$ and $\Bv_j\in\R^n$, their
\emph{outer product} $\Bu_j^{} \Bv_j^T$ is a rank-one matrix of dimension $m\times n$.
\begin{center}
\emph{The formula~(\ref{ME:dyadic}) thus expresses $\BA$ as the weighted sum of $n$
special rank-one matrices.}
\end{center}
This form of the SVD will be especially illuminating as we seek to understand
why the formula $\Bf = \BAs^{-1}\Bb$ proved to be such a problematic way to deblur signals,
why the addition of random noise made the deblurring process even harder, 
and why regularization was so helpful for taming unruly solutions.

\subsection{The singular values of blurring matrices often decay quickly}
\label{ME:ssec:blur_sv}

\begin{figure}[t!]
\begin{center}
\includegraphics[width=3in]{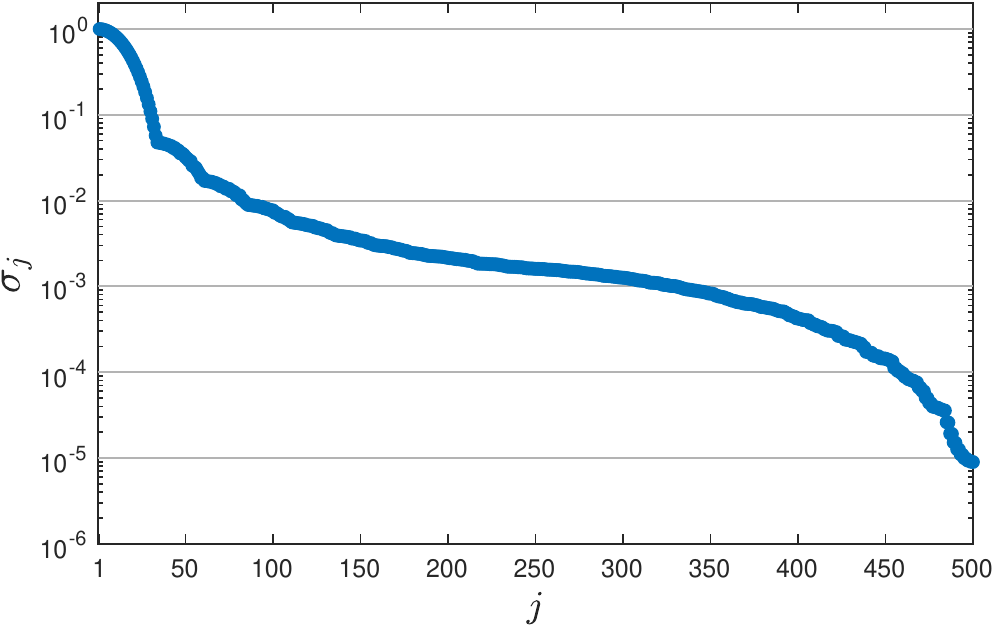}
\end{center}

\vspace*{-18pt}
\caption{\label{ME:fig:blur2_svd}
Singular values of the blurring matrix $\BA$ in~(\ref{ME:blurmat}), 
using the hat function kernel~(\ref{ME:hatker}) with $n=500$ and $z=0.05$, 
as used in Figure~\ref{ME:fig:bw2}.}
\end{figure}

Let us examine the singular values of the blurring matrix $\BA\in\R^{n\times n}$
that we introduced in equation~(\ref{ME:blurmat}).  
Figure~\ref{ME:fig:blur2_svd} shows the singular values for $\BA$ derived from
the hat-function kernel~(\ref{ME:hatker}) with dimension $n=500$ and blurring
factor $z=0.05$.  Note the logarithmic scale on the vertical axis:
the singular values drop by more than a factor $10^5$, 
that is, $\sigma_1/\sigma_{500} > 10^5$.
Such \emph{singular value decay} is typical for blurring matrices.  

Recall the equation $\BA\Bv_j = \sigma_j@\Bu_j$ from~(\ref{ME:eq:Avsigu}). 
The observation that $\BA$ has small singular values $\sigma_j$ 
should mean that there exist some unit vectors $\Bv_j$ that get blurred 
by $\BA$ to very small quantities, $\sigma_j\Bu_j$.  
Since we can view blurring as a kind of ``local averaging,'' we might ask:  
Are there vectors that have large entries, but very small local averages?  
Such vectors would presumably take on positive and negative values at nearby entries 
(so as to average (or blur) out to approximately zero).  
Alternatively, are there vectors that are not much diminished by blurring?  
Such vectors, corresponding to larger singular values, would vary quite slowly 
from entry to entry, so the local average at a point is essentially the same
as the value of the vector at that point.  Said another way, such vectors
should already be blurry to start with, not exhibiting any fine-scale detail.

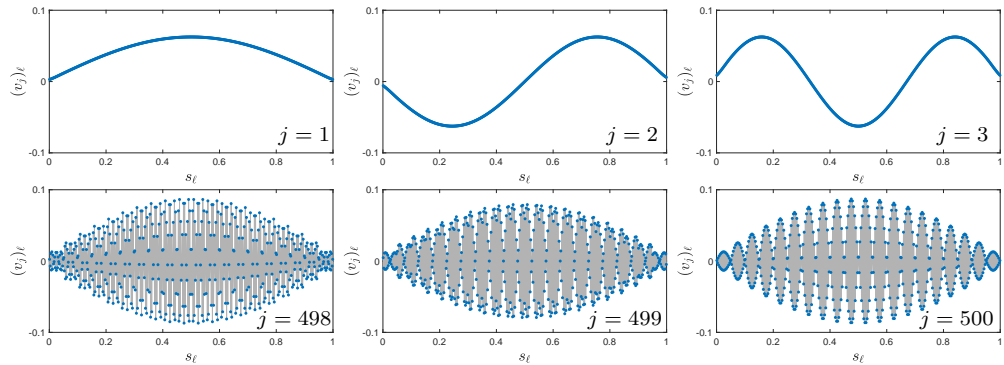
\begin{figure}[t!]
\vspace*{2pt}
\begin{center}
\includegraphics[width=1.7in]{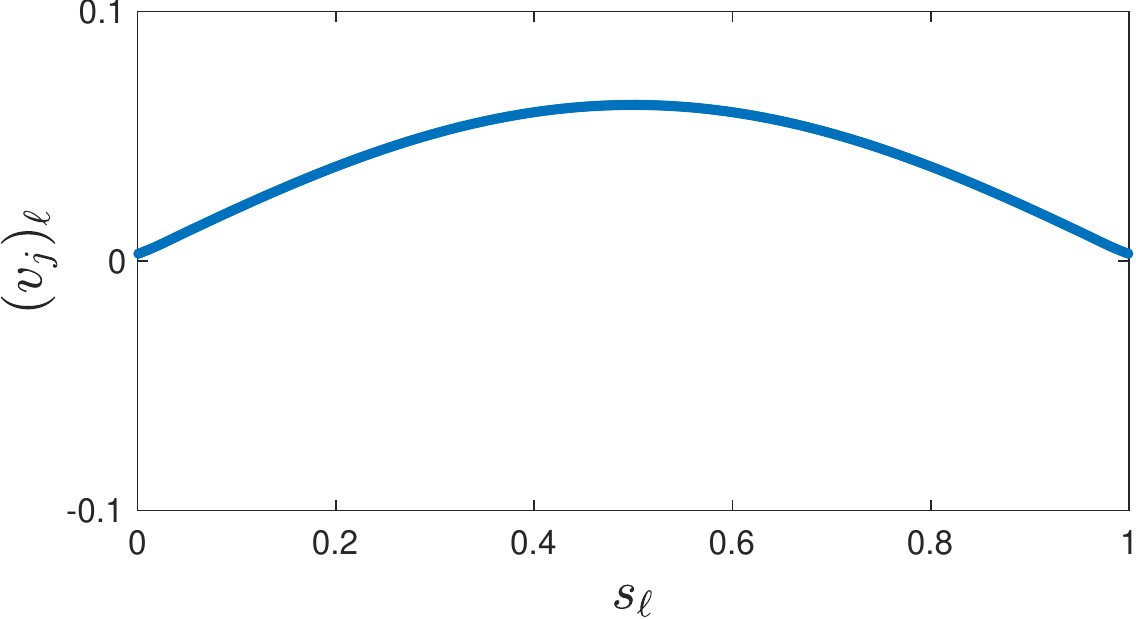}\ 
\includegraphics[width=1.7in]{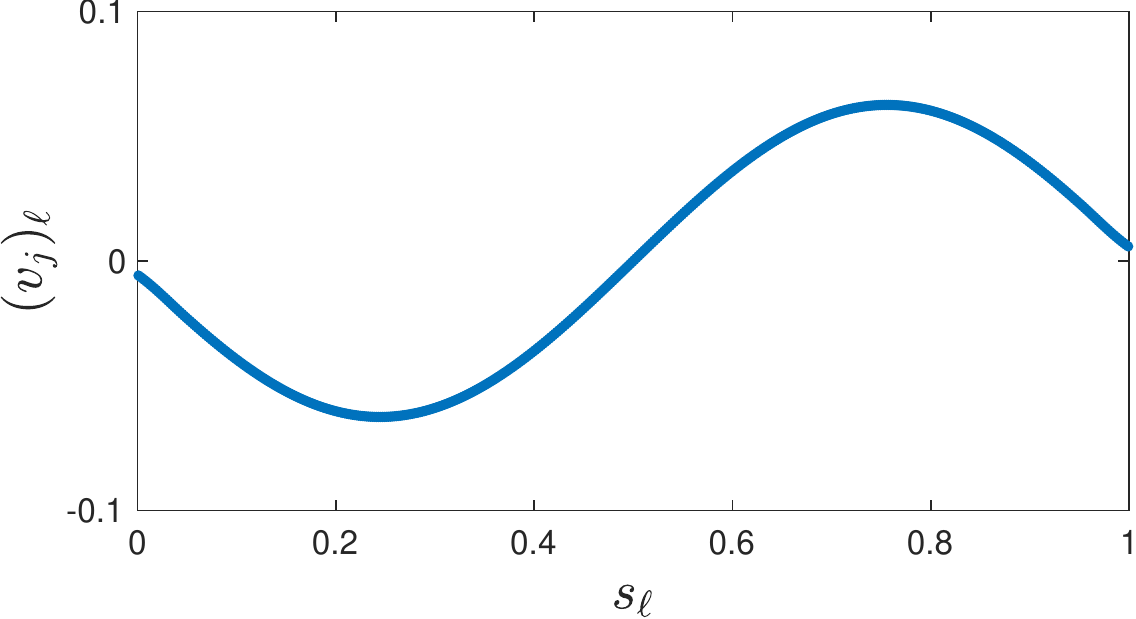}\ 
\includegraphics[width=1.7in]{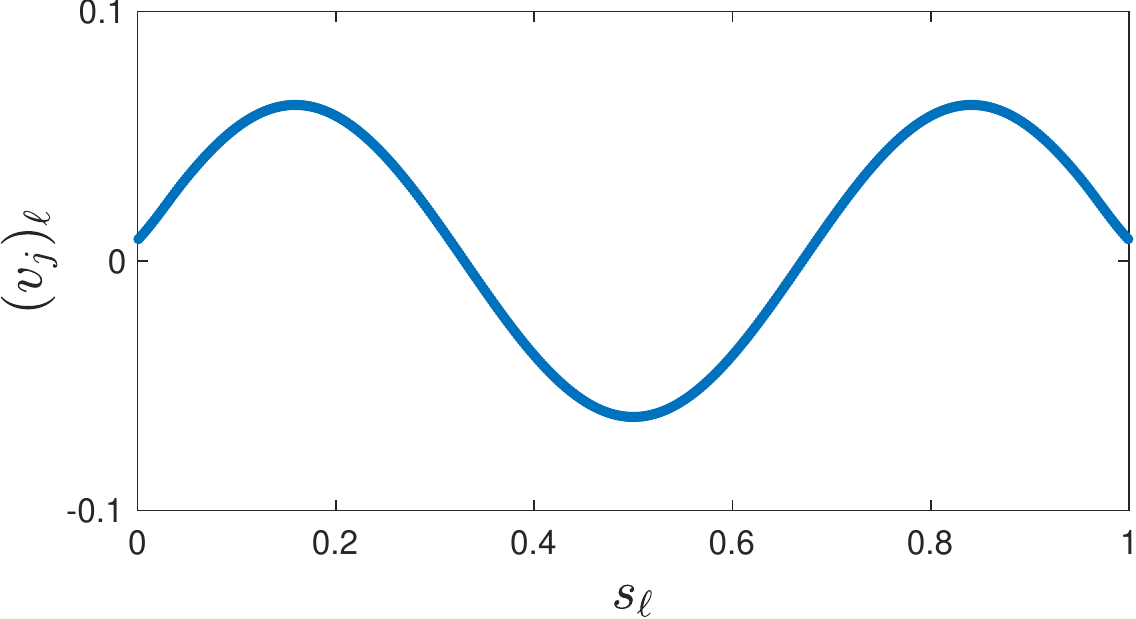}

\includegraphics[width=1.7in]{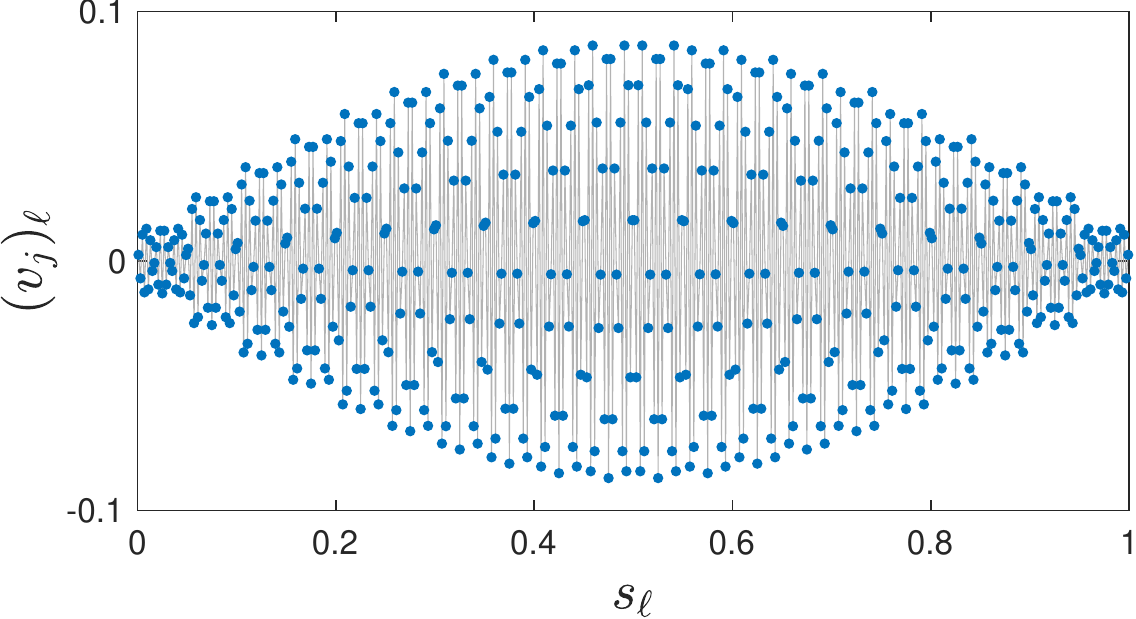}\ 
\includegraphics[width=1.7in]{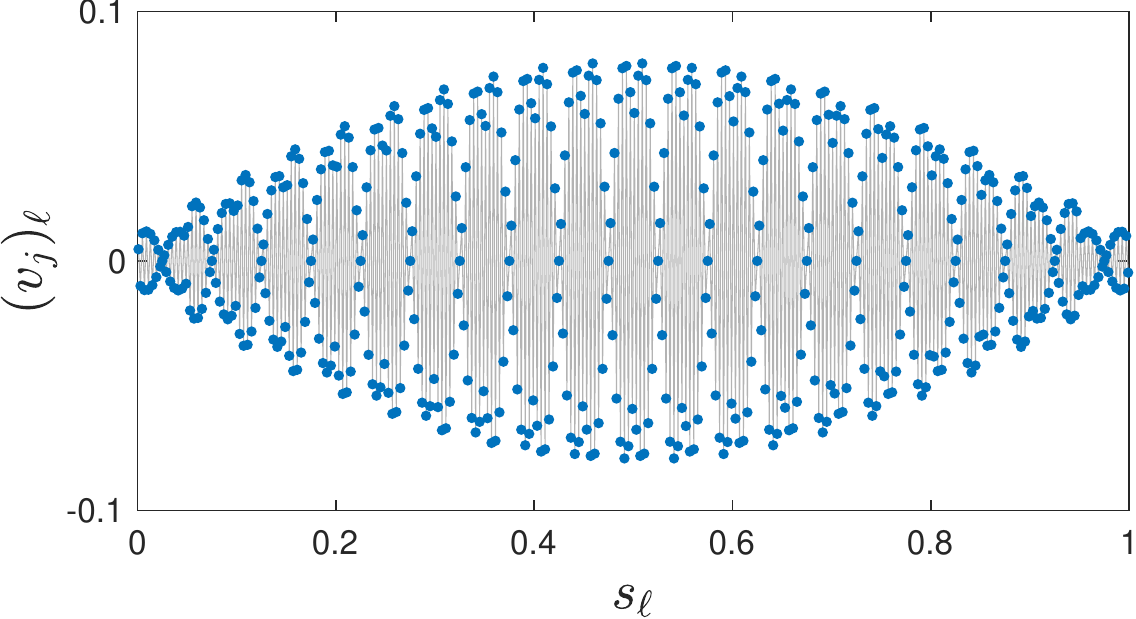}\ 
\includegraphics[width=1.7in]{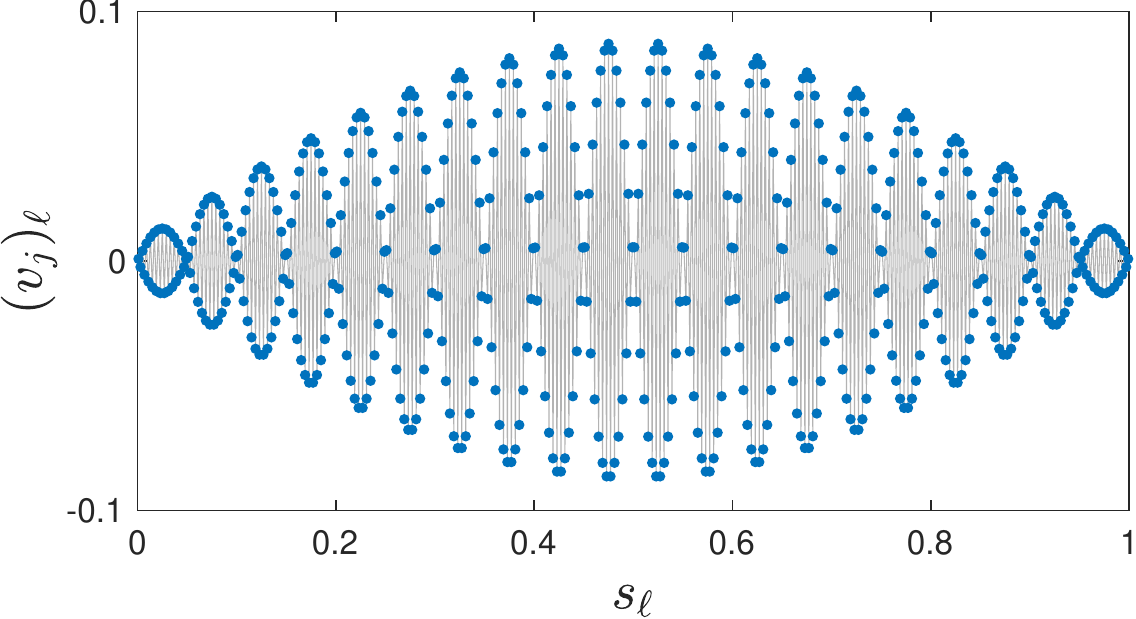}

\begin{picture}(0,0)
\put(-86,96){\footnotesize  $j=1$}
\put( 38,96){\footnotesize  $j=2$}
\put(163,96){\footnotesize  $j=3$}
\put(-94,26.5){\footnotesize  $j=498$}
\put(31,26.5){\footnotesize  $j=499$}
\put(156.5,26.5){\footnotesize  $j=500$}
\end{picture}
\end{center}

\vspace*{-30pt}
\caption{\label{ME:fig:blur2_singvecs}
Singular vectors $\Bv_j$ of the blurring matrix $\BA$ in~(\ref{ME:blurmat}),
using the hat function kernel~(\ref{ME:hatker}) with $n=500$ and $z=0.05$.
The singular vectors associated with the largest singular values ($j=1,2,3$)
vary quite smoothly from entry to entry, resembling sine waves; these smooth
vectors do not change much when blurred.  (We have plotted dots for each
entry of these vectors, but they appear to be fused into smooth, continuous
curves.)  In contrast, the singular vectors
associated with the smallest singular values ($j=498,499,500$)
oscillate quite a bit between entries (highlighted by the gray lines 
connecting the consecutive values); these vectors have local averages 
that are quite close to zero.}
\end{figure}

Figure~\ref{ME:fig:blur2_singvecs} shows the right singular vectors 
associated with the three largest singular values ($\Bv_1$, $\Bv_2$, 
and $\Bv_3$) and the three smallest singular value 
($\Bv_{498}$, $\Bv_{499}$, and $\Bv_{500}$) of $\BA$.\ \ 
As expected, the ``leading'' singular vectors ($\Bv_1$, $\Bv_2$, $\Bv_3$)
are very smooth, varying little from one entry to the next; 
they exhibit no fine-scale detail, and are not much affected by blurring.
In contrast, the ``trailing'' singular vectors ($\Bv_{498}$, $\Bv_{499}$,
$\Bv_{500}$ oscillate quite a bit between entries; they are dominated 
by fine-scale detail, and blurring will smooth them out almost to zero.

What are the implications of these singular values and vectors for the
deblurring operation?

\subsection{Inverting a matrix emphasizes the small singular values}

Suppose $\BA\in\R^{n\times n}$ is a square, invertible matrix. 
Note that $\BU, \BV \in \Rnn$ are square matrices with orthonormal columns,
$\BU^T\BU=\BI$ and $\BV^T\BV=\BI$, which implies that $\BU^T$ and $\BV^T$
function as the inverse of a matrix:  $\BU^{-1}=\BU^T$ and $\BV^{-1}=\BV^T$.\ \ 
Thus one can write
\[ \BAs^{-1} = (\BU\BSigma\BV^T)^{-1} = \BV \BSigma^{-1}\BU^T,\]
or in the dyadic form
\begin{equation} \label{ME:dyadic_inv}
   \BAs^{-1}  = \sum_{j=1}^n {1\over \sigma_j} \Bv_j^{} \Bu_j^T.
\end{equation}
This last equation reveals how the singular values of $\BAs^{-1}$ are the 
\emph{reciprocals} (or \emph{inverses}) of the singular values of $\BA$.\ \ 
Thus if $\BA$ has very \emph{small} singular values, 
$\BAs^{-1}$ will have some very \emph{large} singular values.  
This fact should give us some pause, in light of the observation we made in 
section~\ref{ME:ssec:blur_sv} about the rapidly decaying singular values of blurring matrices! 
(The idea that ``problematic'' matrices have large entries in their inverse also
arose in the simple $2\times 2$ example in section~\ref{ME:sec:Ainv}.)

Now apply the formula~(\ref{ME:dyadic_inv}) to solve the system $\BA\Bf = \Bb$
for the unknown vector $\Bf$.
We find that
\begin{equation} \label{ME:Ainvb}
 \Bf = \BAs^{-1}\Bb = \bigg(\sum_{j=1}^n {1\over \sigma_j} \Bv_j^{} \Bu_j^T\bigg) \Bb 
                      = \sum_{j=1}^n {\Bu_j^T\Bb \over \sigma_j} \Bv_j.
\end{equation}
Here we have distributed the $\Bb$ across the sum, and used the fact that 
$\Bu_j^T\Bb \in \R$ is just a scalar value and can thus be moved to the 
other side of $\Bv_j$.
Equation~(\ref{ME:Ainvb}) is crucial to understanding how the inversion process 
works -- or not.  It means that 
\begin{center}
\em
$\Bf$ can be expressed as a linear combination of the right singular vectors $\Bv_j$,\\
and the weight associated with $\Bv_j$ is given by $(\Bu_j^T\Bb)/\sigma_j$.
\end{center}

\subsection{Barcodes put more ``energy'' in larger singular values}

What can be said about the coefficients $(\Bu_j^T\Bb)/\sigma_j$ for the 
deblurring examples we have been studying?

First, we note that the blurring matrices $\BA\in\Rnn$ that we have been 
considering are \emph{symmetric}, and have matching right and left singular vectors: 
$\Bv_j = \Bu_j$.
Thus, the intuition you can draw from Figure~\ref{ME:fig:blur2_singvecs} about 
the $\Bv_j$ vectors also applies to the $\Bu_j$ vectors.

\begin{figure}[h!]
\begin{center}
\includegraphics[width=2.5in]{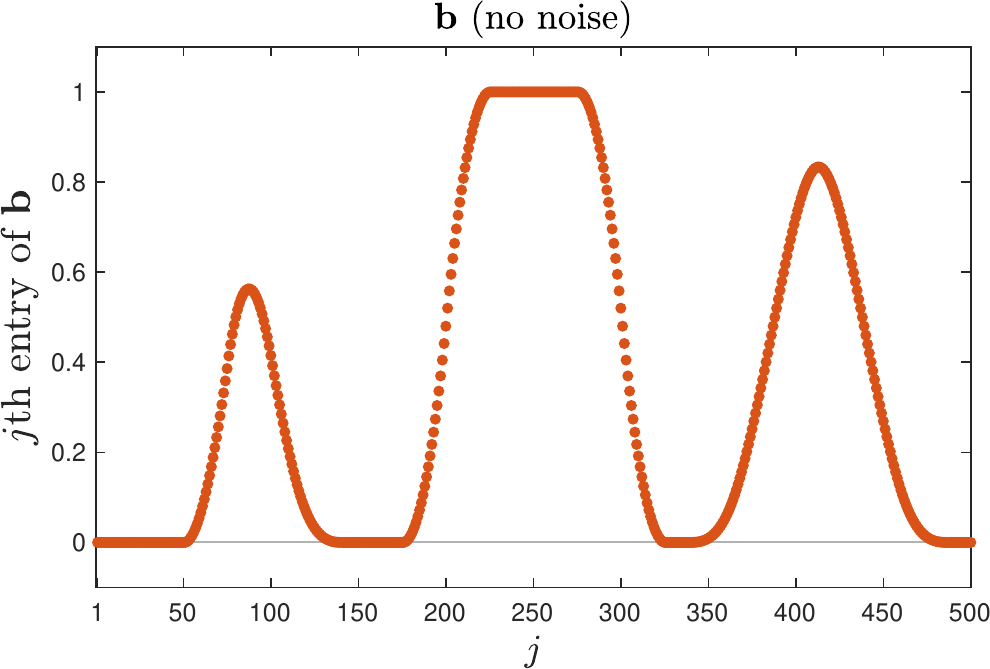}\quad
\includegraphics[width=2.5in]{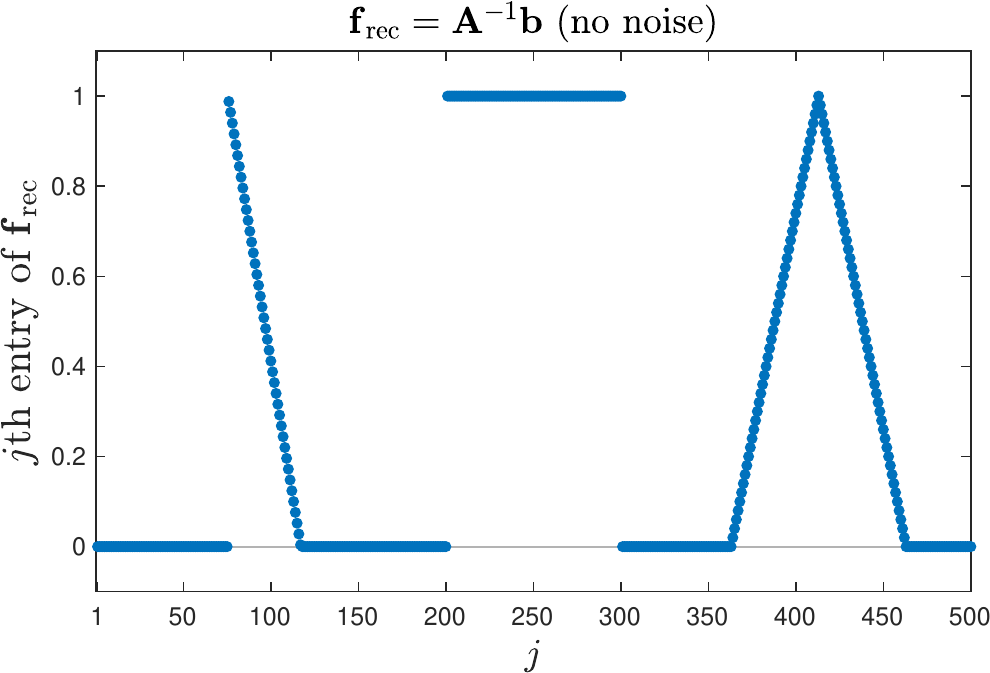}

\includegraphics[width=2.5in]{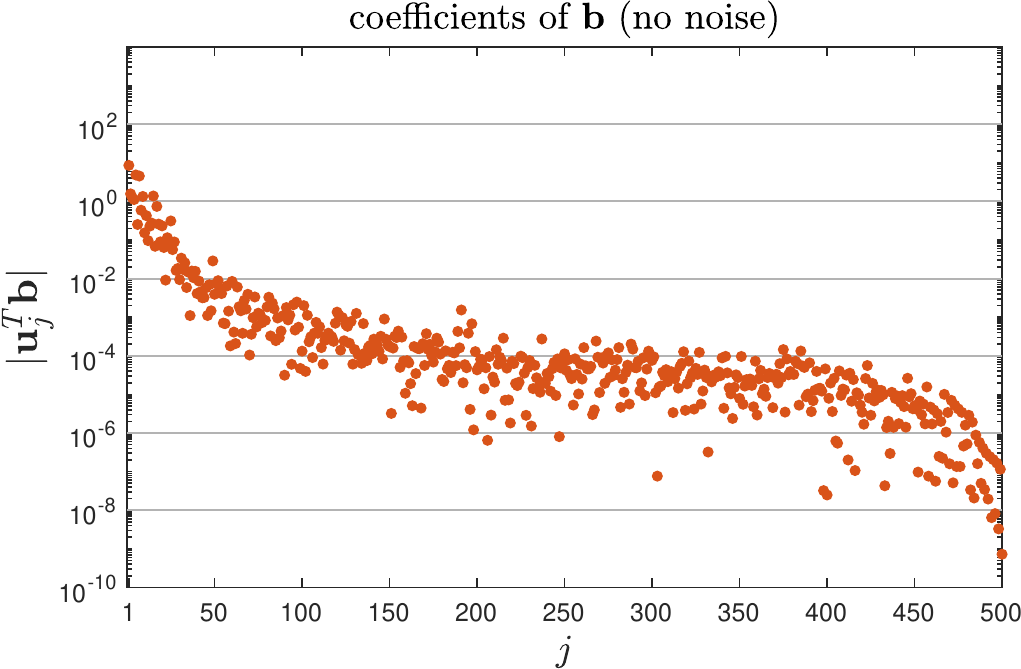}\quad
\includegraphics[width=2.5in]{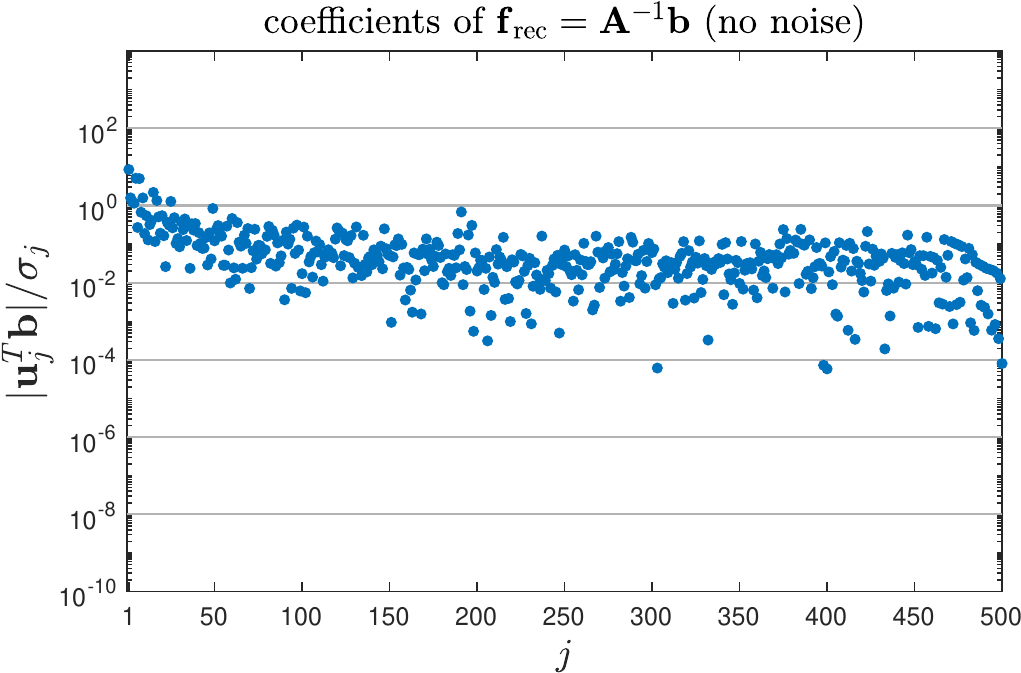}
\end{center}

\vspace*{-15pt}
\caption{\label{ME:fig:blur2_c0}
The vectors $\Bb$ (top left) and $\frec = \BA^{-1}\Bb$ (top right),
with no noise added,
along with the magnitude of the coefficients $\Bu_j^T\Bb$ for $\Bb$ 
(bottom left)
and $\Bu_j^T\Bb/\sigma_j$ for $\frec = \BAs^{-1}\Bb$ (bottom right).
Notice that $\Bb$ has very small components corresponding to the
smallest singular values (i.e., large $j$),  but these values are
lifted much larger in $\frec$ due to the division by $\sigma_j$.
The coefficients for the smoother vector $\Bb$ are smaller 
in magnitude than those for the rougher vector $\frec$, for large~$j$.
}
\end{figure}

Since $\BA$ is square, so too is $\BU\in\Rnn$; recall that the orthonormality 
of the singular vectors, $\BU^T\BU=\BI$, means that $\BU^{-1} = \BU^T$, 
and hence $\BU\BU^T=\BI$, too.  
Thus we can write
\begin{equation} \label{ME:expand_b}
 \Bb = (\BU\BU^T)@\Bb = \BU (\BU^T\Bb) = \sum_{j=1}^n (\Bu_j^T\Bb)@\Bu_j.
\end{equation}
We see that 
\begin{center}
\emph{the coefficient $\Bu_j^T\Bb$ reveals the component 
of $\Bb$ in the $\Bu_j$ direction.}
\end{center}
Since the blurred barcodes are relative smooth 
(see, for example, the plot of $\Bb$ in Figure~\ref{ME:fig:barcode1}), 
we intuitively expect that smooth $\Bu_j$ will make the 
greatest contribution to $\Bb$, while rough $\Bu_j$ will 
make much smaller contributions.
Extrapolating from the six plots of singular vectors shown in 
Figure~\ref{ME:fig:blur2_singvecs}, we expect $\Bu_j$ 
to become increasingly rough as $j$ increases, and thus 
\begin{center}
\emph{we expect $|\Bu_j^T\Bb|$ to exhibit, overall, a decreasing trend as $j$ increases.}
\end{center}
The left plot in Figure~\ref{ME:fig:blur2_c0} reveals that this is 
generally the case: for small $j$, $|\Bu_j^T\Bb|>1$,
while for many larger values of $j$, $|\Bu_j^T\Bb| < 10^{-6}$.
(Notice in particular the sharp drop off that begins around $j=450$.)

The right plot in Figure~\ref{ME:fig:blur2_c0} shows how these coefficients
$\Bu_j^T\Bb$ change when we divide them by $\sigma_j$ in the formula~(\ref{ME:Ainvb}) 
for $\Bf = \BAs^{-1}\Bb$.  Essentially, we are dividing the left plot of 
Figure~\ref{ME:fig:blur2_c0} by the plot of singular values in Figure~\ref{ME:fig:blur2_svd}.
In this case, dividing the coefficients $\Bu_j^T\Bb$ by the 
rapidly decaying singular values $\sigma_j$ elevates the values of 
$\Bu_j^T\Bb/\sigma_j$, most significantly for larger $j$.  
Indeed, aside from $j=1,\ldots, 50$, the value of $|\Bu_j^T\Bb|/\sigma_j$ 
is generally around $10^{-2}$.
Those larger terms (small $j$) will dominate the sum in~(\ref{ME:Ainvb}), 
in this case by a margin sufficient to make the inversion process stable:
$\Bf = \BAs^{-1}\Bb$ produces a reasonable solution in this case,
\emph{provided we are using the exact blurring vector $\Bb$}.  
What happens when we add a whiff of noise to $\Bb$?

\subsection{Random noise scatters ``energy'' across all singular values}

\begin{figure}[t!]
\begin{center}
\includegraphics[width=2.5in]{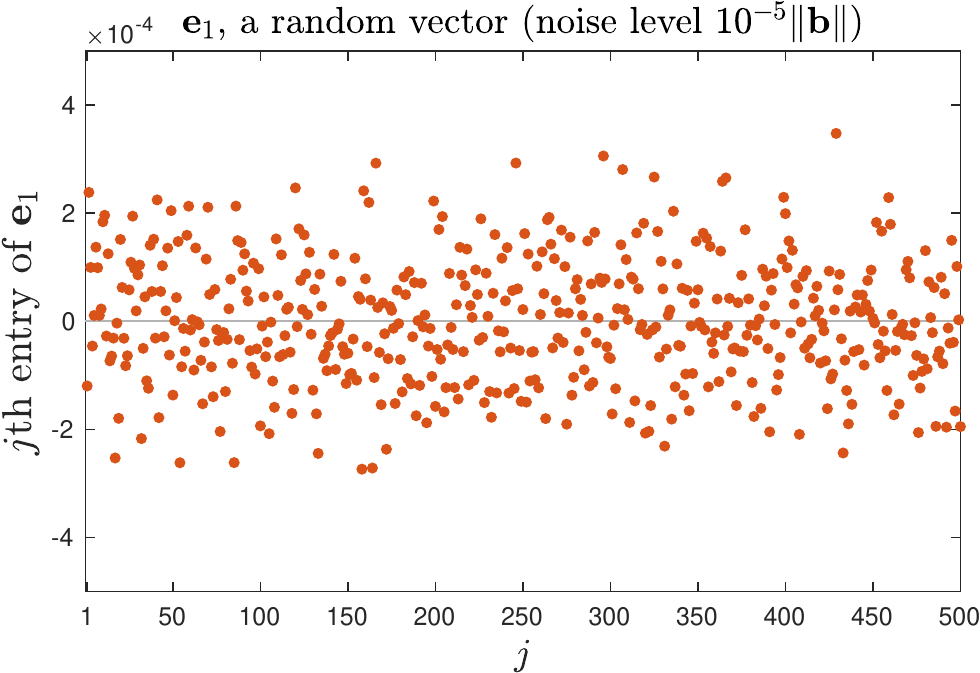}\quad
\includegraphics[width=2.5in]{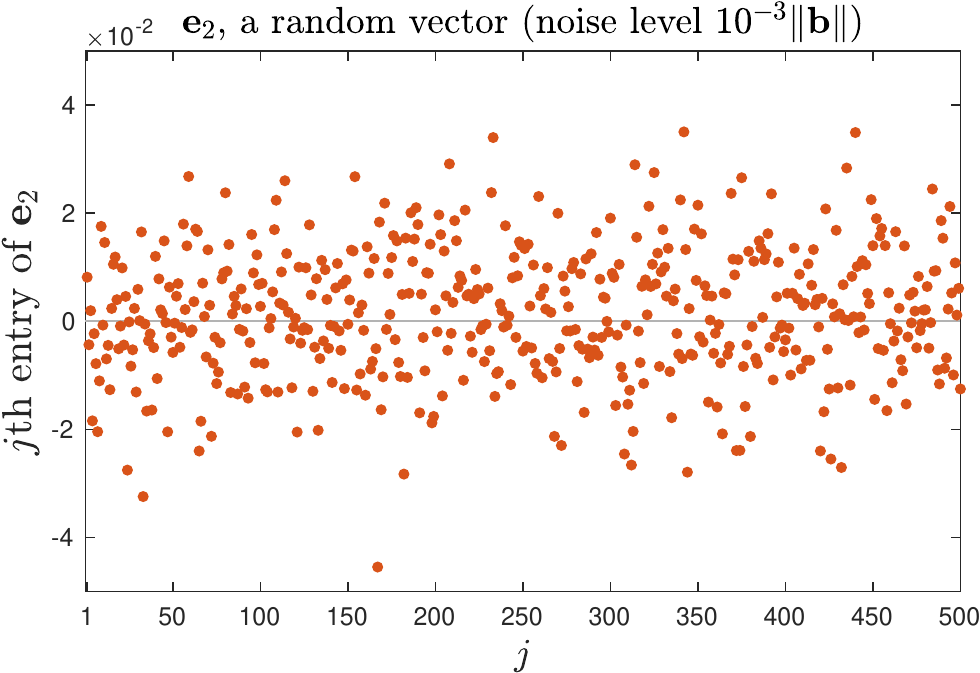}

\includegraphics[width=2.5in]{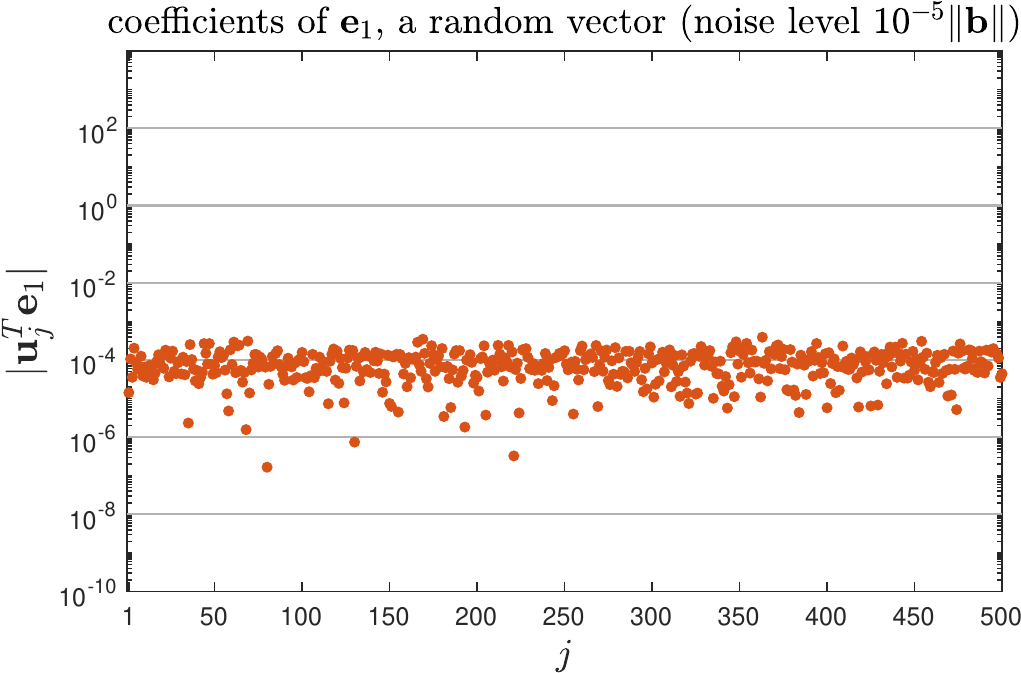}\quad
\includegraphics[width=2.5in]{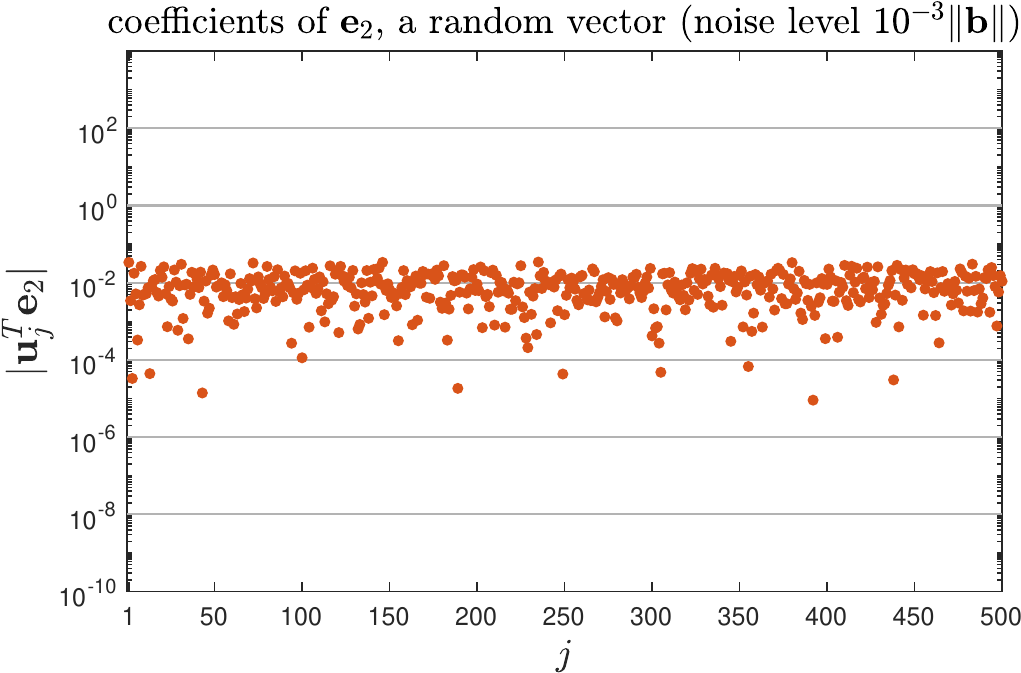}
\end{center}

\vspace*{-15pt}
\caption{\label{ME:fig:blur2_e}
The noise vectors $\Be_1$ (top left,  noise level $10^{-5}\|\Bb\|$)
and $\Be_2$ (top right,  noise level $10^{-3}\|\Bb\|$),
along with the magnitude of the coefficients $\Bu_j^T\Be_1$ 
for $\Be_1$ (bottom left)
and $\Bu_j^T\Be_2$ for $\Be_2$ (bottom right).
The noise vectors are not at all smooth.
As a result, in contrast to the coefficients for $\Bb$ 
(seen on the left in Figure~\ref{ME:fig:blur2_c0}),
these coefficients for $\Be_1$ and $\Be_2$ are fairly uniform across all $j$;
no drop-off is seen as $j$ increases.}
\end{figure}

Now let us consider two random ``noise'' vectors, $\Be_1$ and $\Be_2$.
Both vectors will have normally distributed random entries with 
mean~0 and standard deviation $10^{-5}\|\Bb\|$ (for $\Be_1$) and  
$10^{-3}\|\Bb\|$ (for $\Be_2$).  We can expand these noise vectors
in the same way that we expanded $\Bb$: for $k=1,2$,
\[ \Be_k = \sum_{j=1}^n (\Bu_j^T\Be_k)@\Bu_j.\]
These vectors $\Be_1$ and $\Be_2$ are not at all smooth:
they vary quite a bit (randomly) from entry to entry.  
There is no reason why the smooth singular vectors, like
 $\Bu_1$ and $\Bu_2$, should make any greater or lesser 
contribution to $\Be_k$ than the highly oscillatory vectors
like $\Bu_{n-1}$ and $\Bu_n$.
It is no surprise, then, that their coefficients $\Bu_j^T\Be_k$
look nothing like those seen for $\Bb$ in Figure~\ref{ME:fig:blur2_c0}.
Indeed, Figure~\ref{ME:fig:blur2_e} shows that $\Be_1$ (left) and 
$\Be_2$ (right) have roughly the same contributions from most
singular vectors $\Bu_j$; there is no drop-off as $j$ increases.
\begin{center}
  \emph{Since random vectors $\Be_k$ exhibit no inherent smoothness, \\
we expect $|\Bu_j^T \Be_k|$ to be roughly the same magnitude for all $j$.\\
(There should \emph{not} be any noticeable drop-off as $j$ increases.)}
\end{center}
The general \emph{order of magnitude} of $|\Bu_j^T\Be_k|$ is controlled
by the noise level:  the values in the left plot (noise level $10^{-5}\|\Bb\|$)
are about $10^2$ times smaller than those on the right (noise level $10^{-3}\|\Bb\|$).

\subsection{Random noise + small singular values = disaster}

\begin{figure}[b!]
\begin{center}
\includegraphics[width=2.5in]{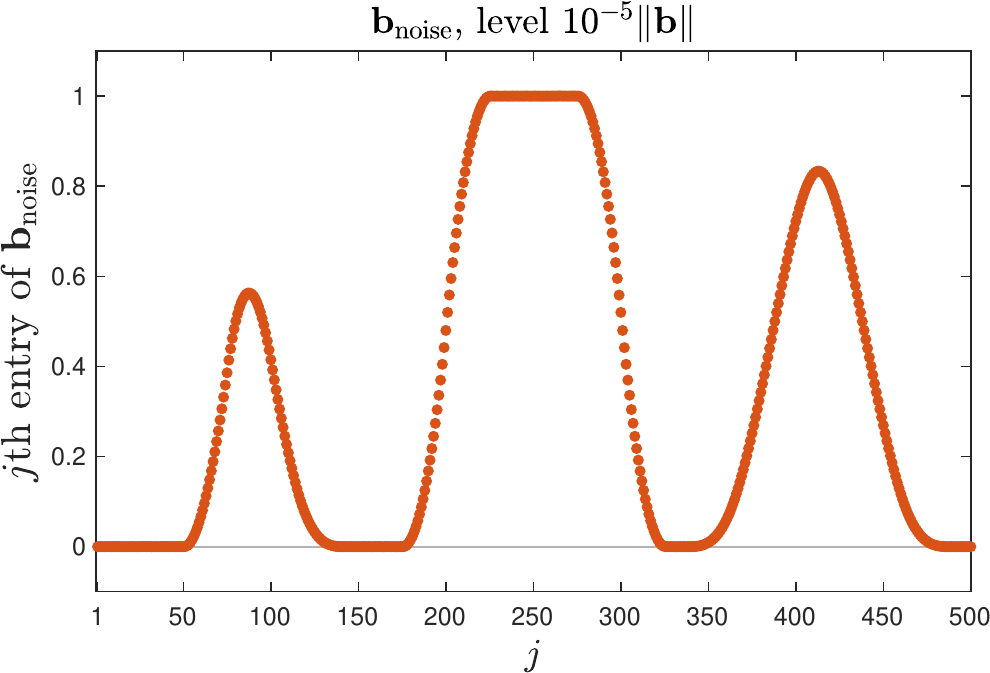}\quad
\includegraphics[width=2.5in]{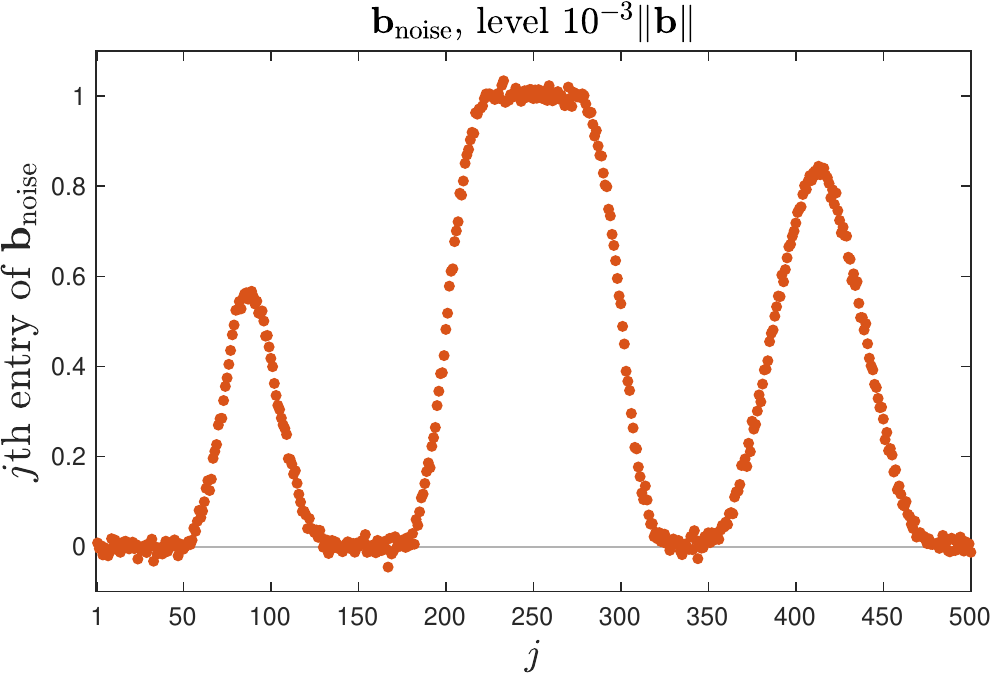}

\includegraphics[width=2.5in]{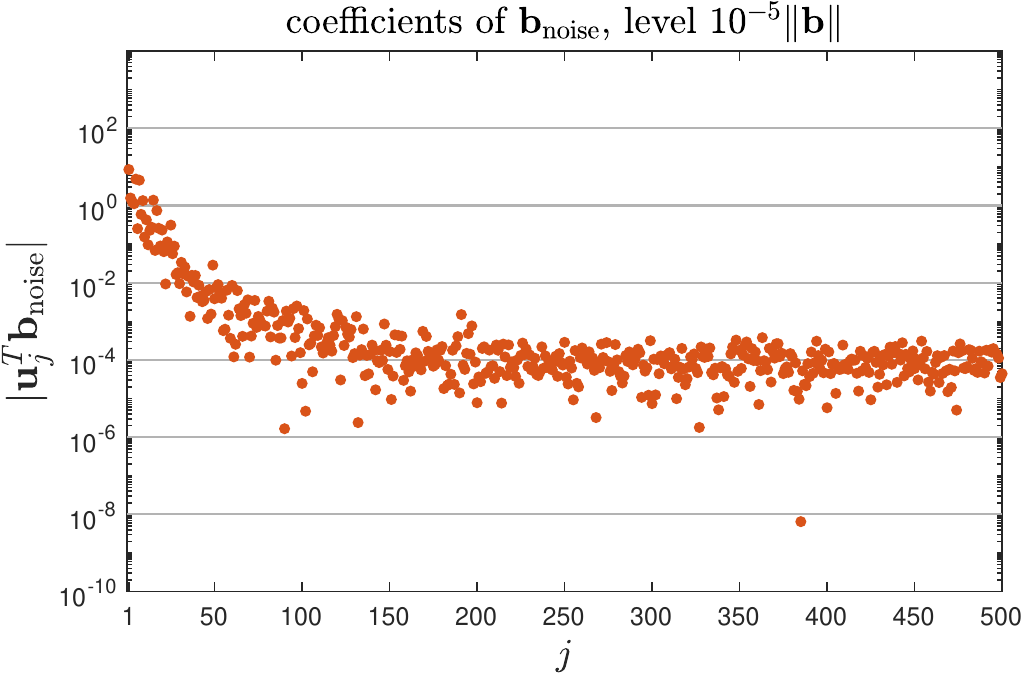}\quad
\includegraphics[width=2.5in]{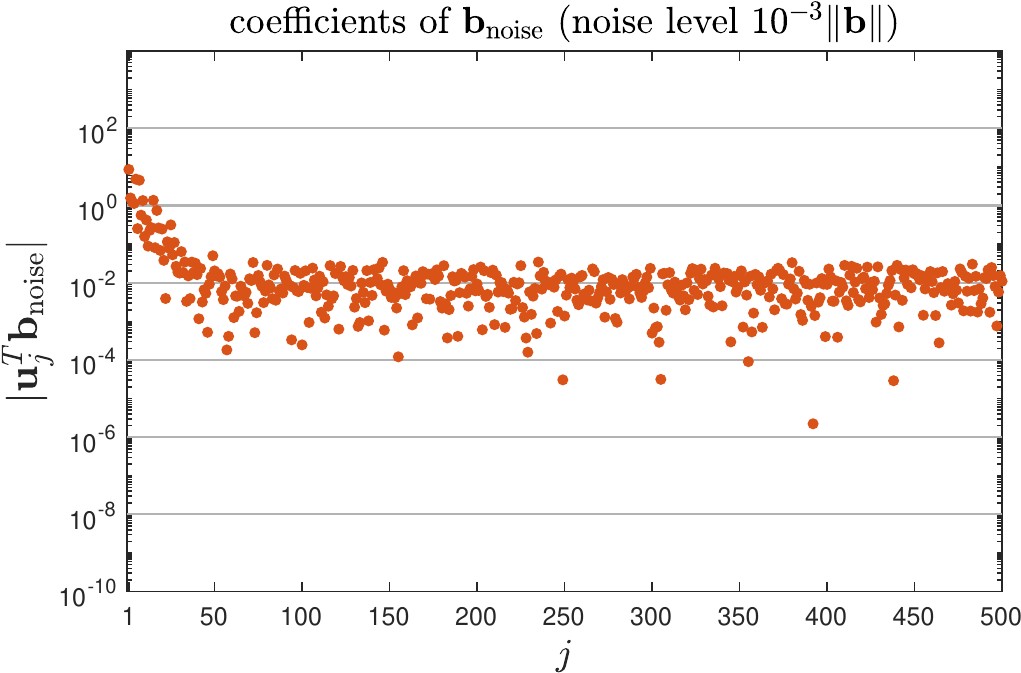}
\end{center}

\vspace*{-15pt}
\caption{\label{ME:fig:blur2_ca}
The polluted vectors $\bnoise$ for noise level $10^{-5}\|\Bb\|$
(top left) and $10^{-3}\|\Bb\|$ (top right),
along with the magnitude of the coefficients $\Bu_j^T\bnoise$ below.
Compare the left plots to those in Figure~\ref{ME:fig:blur2_c0}:
even though the $10^{-3}\|\Bb\|$ noise is hardly evident in 
the plot of $\bnoise$, it significantly elevates the 
coefficients $\Bu_j^T\bnoise$ in the smallest singular values (large $j$).
For the larger noise level on the right, these undesirable 
components are elevated even more.
}
\end{figure}

It is now time to combine the key observations from the last two subsections.
We simulate a noisy sample by constructing the vector 
$\bnoise = \Bb + \Be_k$ for one of the noise vectors $\Be_1$ or $\Be_2$.
When we expand $\bnoise$ in the singular vectors, were merely sum the 
coefficients from $\Bb$ and $\Be_k$:
\[ \bnoise = \sum_{j=1}^n \Big(\Bu_j^T(\Bb + \Be_k)\Big)@\Bu_j 
           = \sum_{j=1}^n (\Bu_j^T\Bb)@\Bu_j 
             + \sum_{j=1}^n (\Bu_j^T\Be_k)@\Bu_j.\]
We are now adding the noise plots from Figure~\ref{ME:fig:blur2_e}
to the coefficients $\Bu_j^T\Bb$ from the left plot of Figure~\ref{ME:fig:blur2_c0}.
The results are shown in Figure~\ref{ME:fig:blur2_ca}:
For small $j$, the coefficients for $\bnoise$ are dominated by the coefficients
$\Bu_j^T\Bb$ for $\Bb$, \emph{until these coefficients decay beneath the noise level
seen in Figure~\ref{ME:fig:blur2_e}}.  Beyond that point, the coefficients
$\Bu_j^T\bnoise$ are dominated by the noise contribution $\Bu_j^T\Be_k$ for $\Be_k$;
the rapid drop-off in $\Bu_j^T\Bb$ is completely swamped by the noise.

\begin{figure}[t!]
\begin{center}
\includegraphics[width=2.5in]{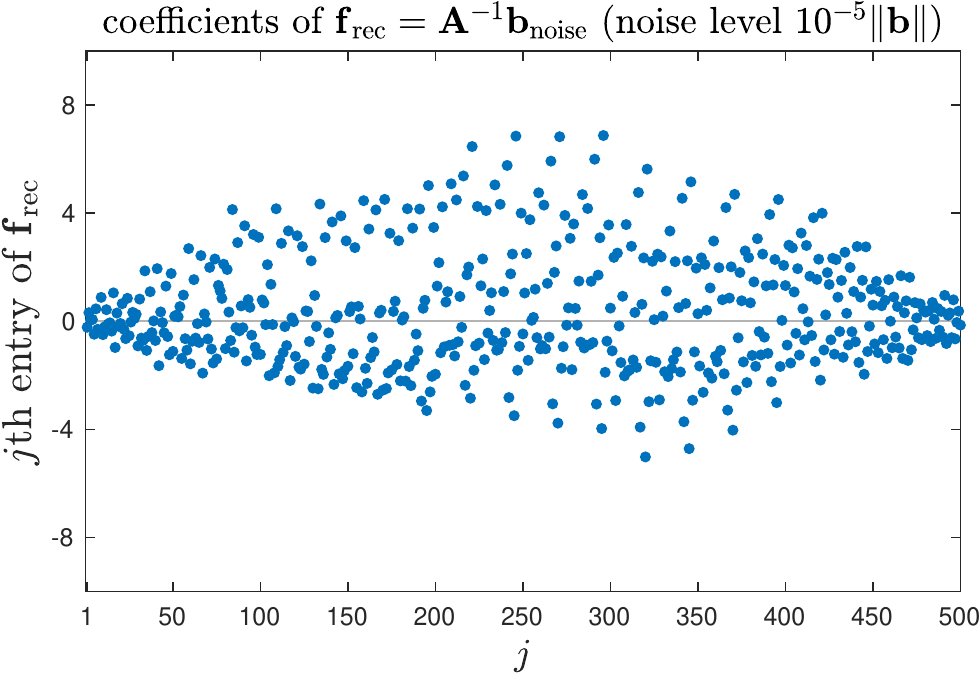}\quad
\includegraphics[width=2.5in]{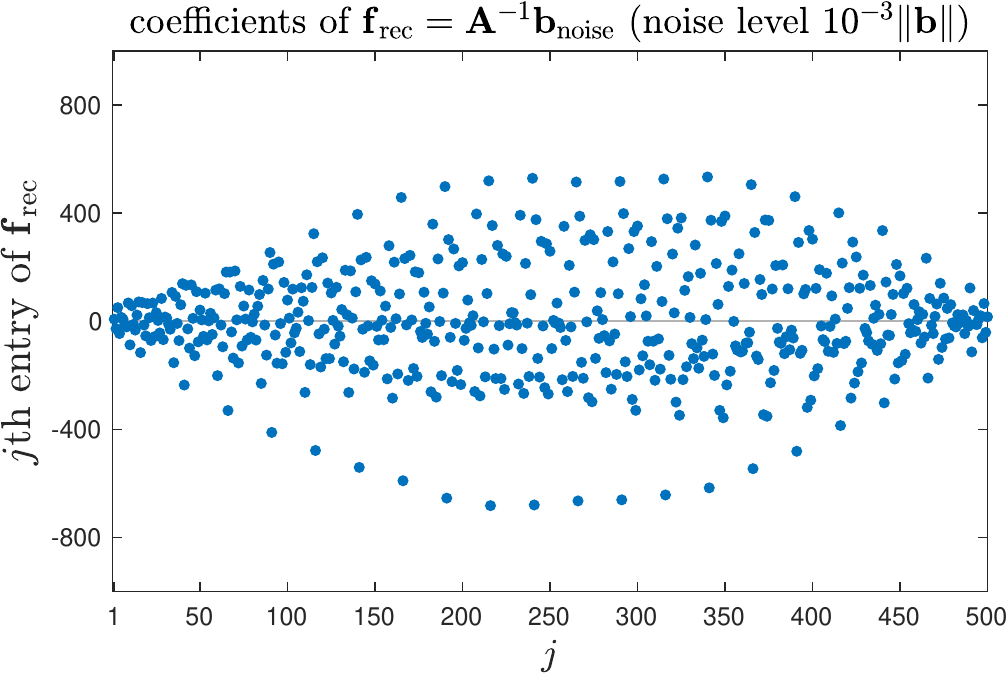}

\includegraphics[width=2.5in]{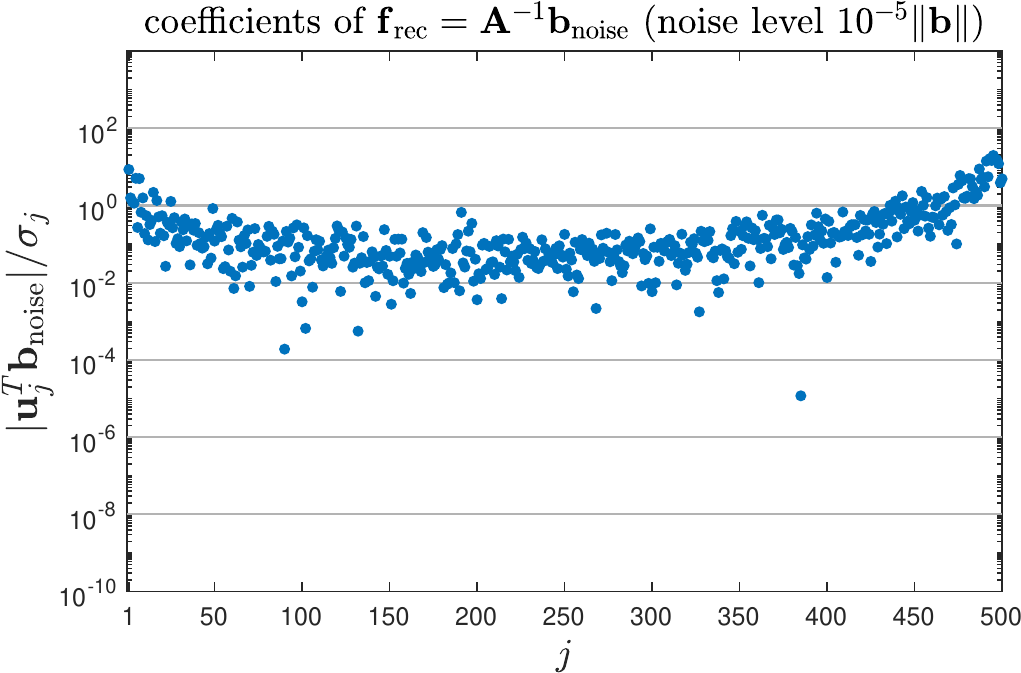}\quad
\includegraphics[width=2.5in]{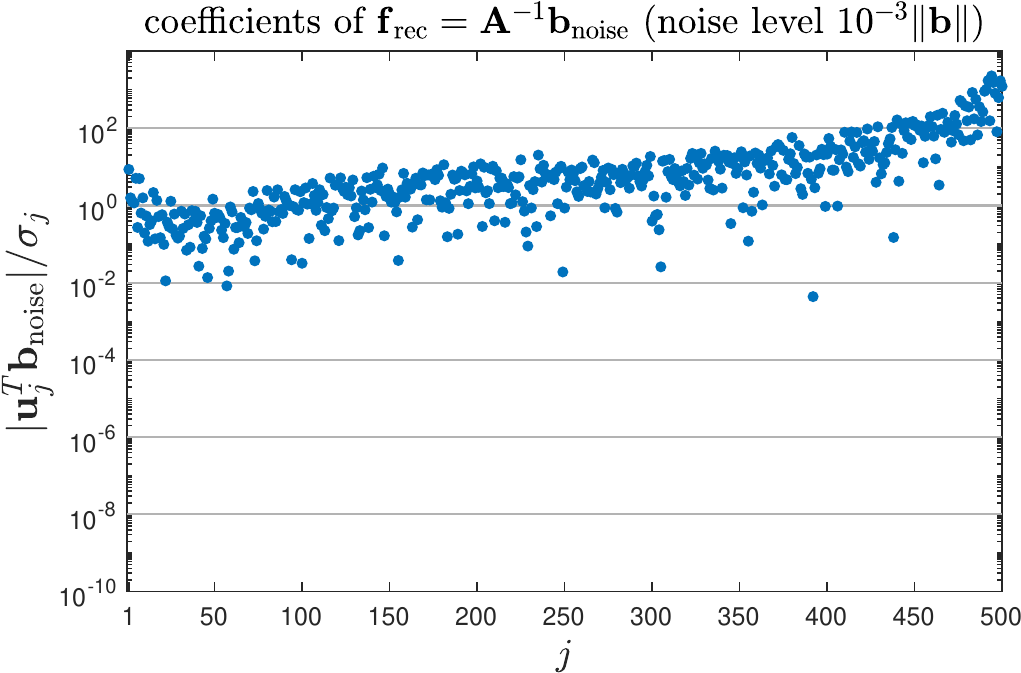}
\end{center}

\vspace*{-15pt}
\caption{\label{ME:fig:blur2_cb}
The polluted $\bnoise$ vectors shown in Figure~\ref{ME:fig:blur2_ca} 
lead to catastrophically poor recovered solutions $\frec$,
for noise level $10^{-5}\|\Bb\|$ (top left) 
and $10^{-3}\|\Bb\|$ (top right).
The bottom plots, which show the coefficients $(\Bu_j^T\bnoise)/\sigma_j$ for 
$\frec = \BAs^{-1} \bnoise$, make the problem apparent:
For large $j$, moderate values of $\Bu_j^T\bnoise$ are divided by small singular
values $\sigma_j$, causing these $(\Bu_j^T\bnoise)/\sigma_j$ terms to dominate.
The effect is disastrous: the contributions to $\frec$ from the trailing singular vectors
(e.g., bottom of Figure~\ref{ME:fig:blur2_singvecs}) totally dominate the desirable 
part of the solution associated with the leading singular vectors 
(e.g., top of Figure~\ref{ME:fig:blur2_singvecs}).
}
\end{figure}

This perspective exposes the problem that occurs when we compute 
\begin{equation} \label{ME:frec_sv}
 \frec = \BAs^{-1}\bnoise 
           = \sum_{j=1}^n {\Bu_j^T\bnoise \over \sigma_j}@\Bv_j.
\end{equation}
When we divide $\Bu_j^T\bnoise$ by $\sigma_j$, we can no longer count on the
rapidly decaying values of $\Bu_j^T\Bb$ to counteract the rapidly decaying 
singular values $\sigma_j$ (as seen on the right plot in Figure~\ref{ME:fig:blur2_c0}).
\begin{center}
\em 
For small $j$, $\displaystyle{ {\Bu_j^T\bnoise \over \sigma_j} \approx {\Bu_j^T\Bb \over \sigma_j}}$, which is accurate.\\[5pt]
For large $j$, $\displaystyle{ {\Bu_j^T\bnoise \over \sigma_j} \approx {\Bu_j^T\Be_k \over \sigma_j}}$, which is rubbish.
\end{center}
The major problem is, as we can see in Figure~\ref{ME:fig:blur2_cb},
\emph{the rubbish overwhelms the accurate data}.
As a result, we will get a nonsense result for $\frec = \BAs^{-1}\bnoise$. 
Take a careful look back on the plots of $\frec$ shown in 
Figure~\ref{ME:fig:bw2} for these two same noise vectors.
What shape does $\frec$ take?  Since the expansion~(\ref{ME:frec_sv})
for $\frec$ is dominated by $(\Bu_j^T\bnoise/\sigma_j)\Bv_j$ terms for large values of $j$, 
it should be no surprise that the computed $\frec$ is dominated by shapes that 
resemble the trailing singular vectors $\Bv_j$, 
like those pictured in the bottom half of Figure~\ref{ME:fig:blur2_singvecs}.
\begin{center}
\emph{When contributions from small singular values (like $\sigma_n$) dominate $\frec = \BAs^{-1}\bnoise$,\\
expect $\frec$ to resemble the shape of the corresponding singular vectors (like $\Bv_n$).}
\end{center}
\subsection{Regularization controls small singular values in the inversion}

Now the key problem has been revealed:  $\bnoise$ has nontrivial but unwanted 
components in the most highly oscillatory singular vectors of $\BA$ (i.e., large $j$), 
and when we compute $\frec = \BAs^{-1}\bnoise$ the division by the singular values 
magnifies those unwanted components to a dominant level.  
How does regularization fix this issue?

\begin{figure}[b!]
\begin{center}
\includegraphics[width=2.5in]{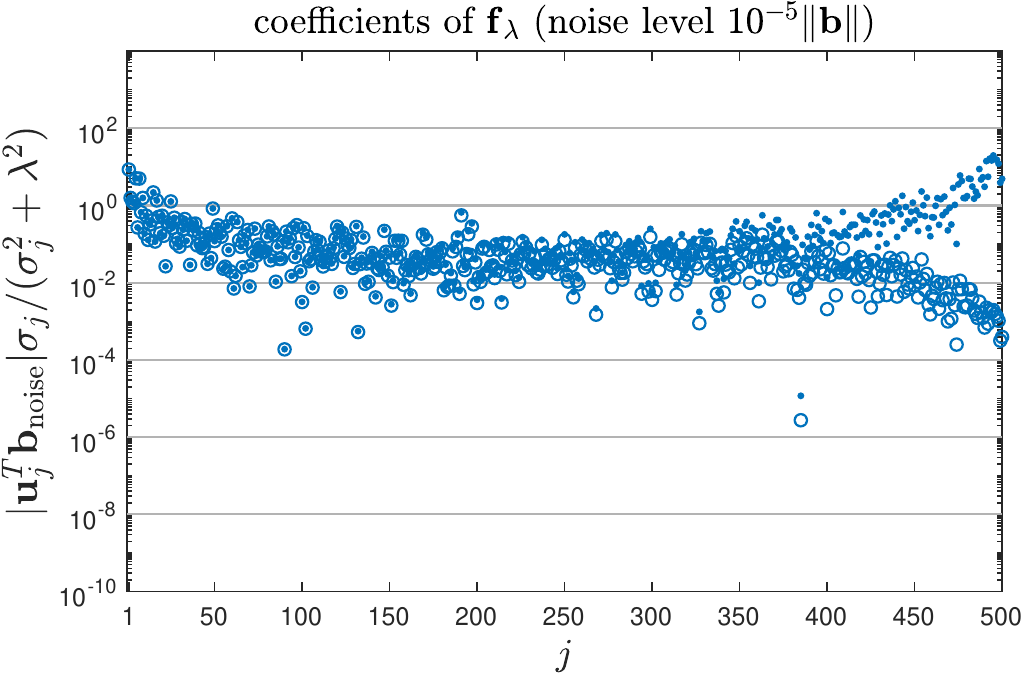}\quad
\includegraphics[width=2.5in]{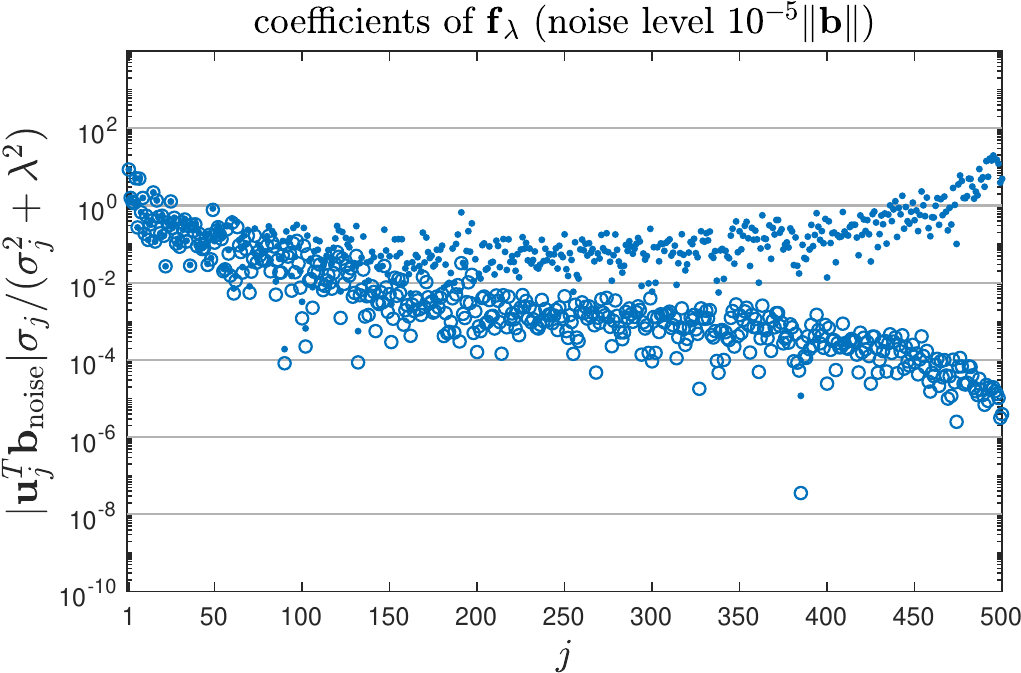}

\includegraphics[width=2.5in]{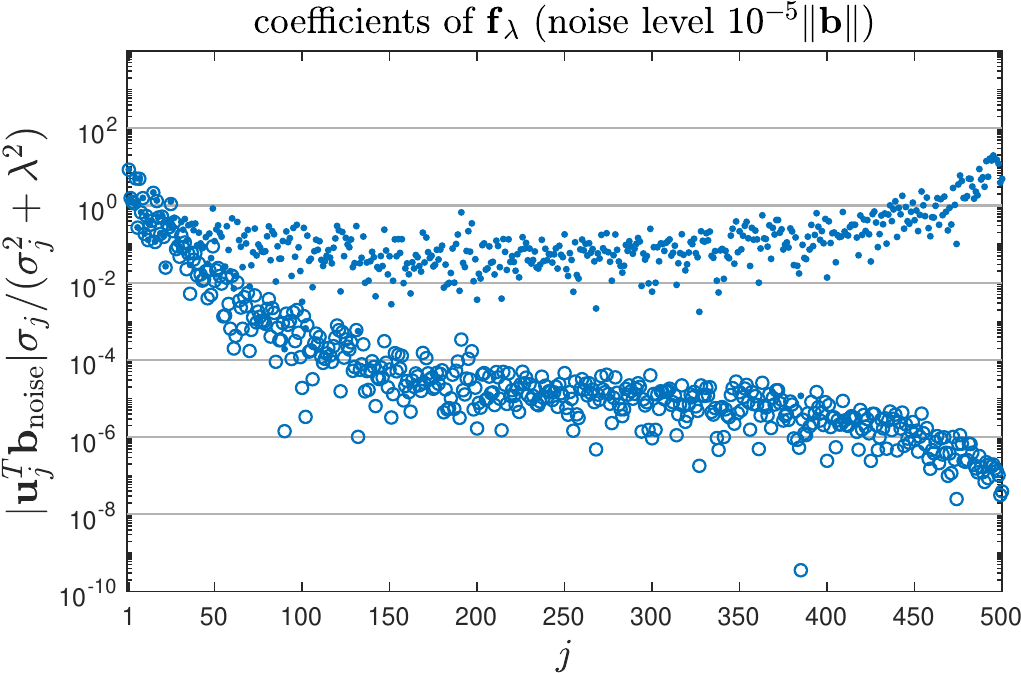}\quad
\includegraphics[width=2.5in]{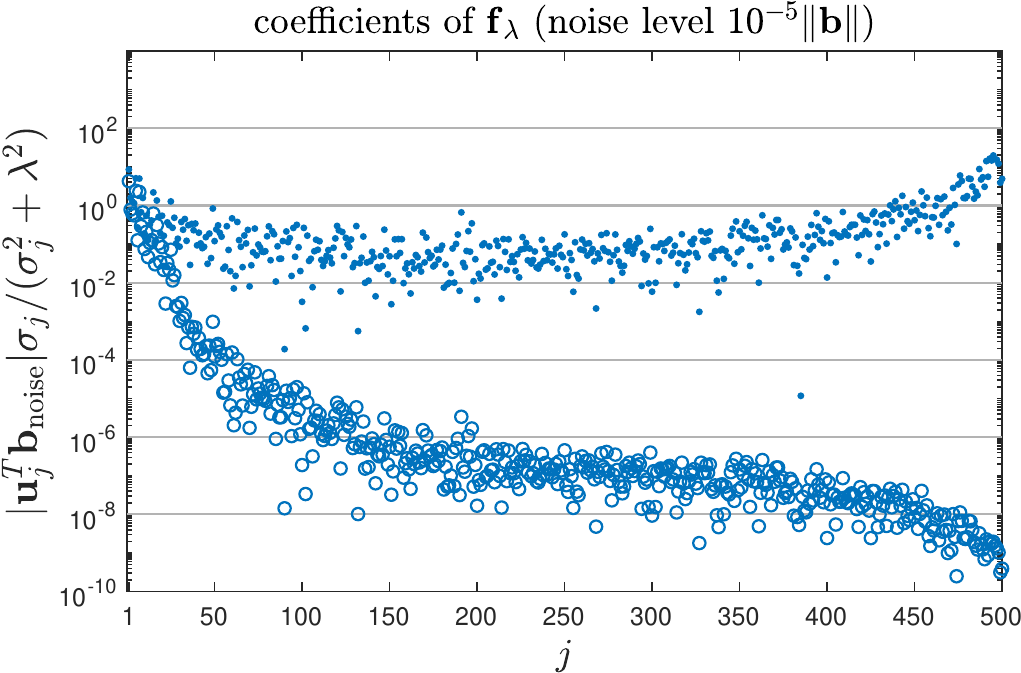}

\begin{picture}(0,0)
\put(-155,150.5){\footnotesize $\lambda=10^{-3}$}
\put(35,150.5){\footnotesize $\lambda=10^{-2}$}
\put(-155,31){\footnotesize $\lambda=10^{-1}$}
\put(35,31){\footnotesize $\lambda=10^{0}$}
\put(-100,97){\footnotesize \emph{no regularization}}
\put(-95,42){\footnotesize \emph{regularization}}
\end{picture}
\end{center}

\vspace*{-25pt}
\caption{\label{ME:fig:blur2_c_reg}
Magnitude of the coefficients of the regularized solution
 $\flam$ derived from $\bnoise$ with noise level $10^{-5}\|\Bb\|$,
for four different values of $\lambda$. 
(These plots correspond to the bottom four solutions in Figure~\ref{ME:fig:bw_reg}.)
The solid dots show $|\Bu_j^T\bnoise|/\sigma_j$ for $\frec$; the hollow dots show
the regularized coefficients $|\Bu_j^T\bnoise| \sigma_j/(\sigma_j^2+\lambda^2)$.
By adjusting $\lambda$ we can dial back the influence of the smaller singular values,
but taking $\lambda$ too large will damp out too many terms, thus giving
an overly smoothed solution.
}
\end{figure}

To answer that question, we first need to derive an expression for the 
regularized solution $\flam$ in the same style as the 
expansion~(\ref{ME:frec_sv}) of $\frec$ in the basis of $\{\Bv_j\}$ vectors.
Recall from~(\ref{ME:regsol}) that we can write
\[ \flam = (\BAs^T\BA + \lambda^2\BI)^{-1} \BAs^T\Bb.\]
We can use the SVD to reveal additional structure behind this formula.
Start at the matrix level with $\BA = \BU\BSigma\BV^T$, and recall that
orthonormality of the singular vectors ensures $\BU^T\BU=\BI$ and $\BV\BV^T=\BI$.\ \ 
Using these facts we find that
\begin{eqnarray*}
\BAs^T\BA + \lambda^2@\BI \;&=&\; (\BU\BSigma\BV^T)^T (\BU\BSigma\BV^T) + \lambda^2 \BI \\[5pt]
                          \;&=&\; \BV \BSigma^T \BU^T\BU\BSigma\BV^T + \lambda^2 \BV\BV^T 
                          \ =\ \BV\big(\BSigma^2+\lambda^2\BI\big)\BV^T,
\end{eqnarray*}
where we have used the fact that $\BSigma^T=\BSigma$.\ \ 
Invert this formula to get
\begin{eqnarray*}  
\big(\BAs^T\BA + \lambda^2@\BI\big)^{-1} 
              \;&=&\; \Big(\BV\big(\BSigma^2+\lambda^2\BI)\BV^T\Big)^{-1}  \\
              \;&=&\; (\BV^T)^{-1} \big(\BSigma^2+\lambda^2\BI\big)^{-1}\BV^{-1}  
              \ =\ \BV \big(\BSigma^2+\lambda^2\BI\big)^{-1}\BV^T.
\end{eqnarray*}
Multiply this last expression against $\BA^T\Bb = \BV\BSigma\BU^T\Bb$ to arrive at
\[ \flam = \BV \big(\BSigma^2 + \lambda^2@\BI\big)^{-1}\BSigma@ \BU^T\Bb,\]
which can be expressed in terms of the individual singular values and singular vectors,
akin to the formula~(\ref{ME:frec_sv}) for $\frec$, as:
\begin{equation}   \label{ME:reg_sv}
\flam = \sum_{j=1}^n \bigg((\Bu_j^T\Bb) {\sigma_j\over \sigma_j^2 + \lambda^2}\bigg) \Bv_j.
\end{equation}
Notice the key role that $\lambda>0$ plays in the formula~(\ref{ME:reg_sv}) when $\sigma_j$
is very small.  As $\sigma_j\to 0$, observe that $1/\sigma_j \to \infty$, 
(causing blow-up in the formula for $\frec$) 
while $\sigma_j / (\sigma_j^2+\lambda^2) \to 0$
(avoiding blow-up in the formula for $\flam$). 
At the same time, for large values $\sigma_j \gg \lambda$, 
we have $\sigma_j/(\sigma_j^2+\lambda^2) \approx 1/\sigma_j$, 
and so in this regime $\lambda$ does not interfere much with the solution.

In short, by tuning $\lambda$ correctly we can suppress components of $\flam$
associated with small singular values, while allowing the desirable components
to pass through with little change.  
Figure~\ref{ME:fig:blur2_c_reg} clearly illustrates this effect,
comparing the coefficients $(\Bu_j^T\bnoise)/\sigma_j$ for $\frec$ 
(solid dots) to the coefficients $(\Bu_j^T\bnoise) \sigma_j/(\sigma_j^2+\lambda^2)$
of $\flam$ (hollow dots), for the same $\lambda$ values used in the bottom four plots of $\flam$ in 
Figure~\ref{ME:fig:bw_reg} (noise level $10^{-5}\|\Bb\|$).
In all four of these examples, the `tail' of elevated coefficients
for $\frec$ seen in the left plot of Figure~\ref{ME:fig:blur2_cb}
is suppressed by regularization, to an increasing degree as $\lambda$ grows.
In Figure~\ref{ME:fig:blur2_c_reg} we saw that the three values 
$\lambda=10^{-3}$, $10^{-2}$, and $10^{-1}$ give reasonable 
solutions that increase in smoothness, consistent with the 
more aggressive suppression of $\Bv_j$ for large $j$ due to the 
$\sigma_j/(\sigma_j^2+\lambda^2)$ term.  
Taking $\lambda=10^0$ suppresses too many terms, leading to an 
unacceptably smoothed solution.

\medskip
This section has illustrated how the singular value decomposition can illuminate
the regularization process.  Indeed, it is a remarkable tool for analyzing the
solution of linear systems, least squares problems, low-rank approximation 
and dimension reduction, and many other rich applications.
In the same way, the simple deblurring example described in this manuscript
is but a taste of the rich and vital world of \emph{inverse problems},
with applications that range from image deblurring to radar, 
from landmine detection to ultrasound.
Many avenues remain for your exploration!

\begin{reflections}
\item It might seem like a leap to move from the product-of-matrices form
      of the SVD in~(\ref{ME:SVD}) and the dyadic form~(\ref{ME:dyadic}).
      Consider the following $2\times 2$ matrices
\[ \BU = \left[\begin{array}{cc} \Bu_1 & \Bu_2 \end{array}\right]
       = \left[\begin{array}{cc} u_{1,1} & u_{2,1} \\ u_{1,2} & u_{2,2} \end{array}\right], 
\qquad 
   \BV = \left[\begin{array}{cc} \Bv_1 & \Bv_2 \end{array}\right]
       = \left[\begin{array}{cc} v_{1,1} & v_{2,1} \\ v_{1,2} & v_{2,2} \end{array}\right], \]
where, e.g., $u_{j,k}$ refers to the $k$th entry of the vector $\Bu_j$.

By multiplying each expression out, entry-by-entry, 
show that $\BU\BV^T = \Bu_1^{}\Bv_1^T + \Bu_2^{}\Bv_2^T$.

\item Think back to the local averaging kernel~(\ref{ME:avgker}).
      Can you design a function $f(t)$ such the blurred function
      \[ b(s) = \int_0^1 h(s,t) f(t)\,\dop t \]
      is close to zero for \emph{all} $s\in[0,1]$?  
      How would such a function relate to the singular vectors of the
      discretization matrix $\BA$ for this blurring kernel?
\end{reflections}

\noindent
Now would be a good time to explore Exercises~\ref{ME:ex:comp_svs}--\ref{ME:ex:trunc_svd} 
starting on page~\pageref{ME:ex:comp_svs}.

\bigskip

\hrule

\begin{exercises}
\setcounter{section}{2} \setcounter{prob}{0}

\item \label{ME:ex:fw}
Experiment with the blurring operation as shown in Figure~\ref{ME:fig:fw}.
You can generate the function $f$ using the code below.
\begin{lstlisting}[language=Python]
def f_function(t=0.5):
    f1 = (t>=0.15)*np.maximum(1-12*(t-.15),0)     # a down ramp
    f2 = np.abs(t-0.5)<=.1                        # a step
    f3 = np.maximum(1-10*np.abs(t-0.825),0)       # a hat
    return f1+f2+f3 
\end{lstlisting}
For $n=100$, you can create a vector $\Bf$, blurring matrix $\BA$, and blurred
vector $\Bb = \BA\Bf$ using the code below.  Note that \verb|build_blur_A| 
uses the Gaussian blurring kernel, here called with parameter $z=0.025$.
\begin{lstlisting}[language=Python]
n = 100
z = 0.025
t = np.array([(k+.5)/n for k in range(n)])  # t_k points
f = f_function(t)                           # evaluate f(t_k)
A = build_blur_A(n,z)                       # build blurring matrix
b = A@f                                     # b=A*f (matrix times vector)
\end{lstlisting}
You can then plot the entries of $\Bf$ and $\Bb$ as dots 
using the following commands.
\begin{lstlisting}[language=Python]
import matplotlib.pyplot as plt

plt.plot(t,f_function(t),'.',markersize=5,label='f')
plt.plot(t,b,'.',markersize=5,label='b')
plt.legend()
\end{lstlisting}

Execute these commands to produce your own plot of $\Bf$ and $\Bb$.
(Your plot will differ slightly from the one shown in Figure~\ref{ME:fig:fw},
which uses the local averaging blurring kernel~(\ref{ME:avgker}), 
whereas the code above uses the Gaussian kernel~(\ref{ME:gaussker}).)

Now experiment with the value of $z$, the dimension $n$, and the choice of 
blurring kernel.  
(As you change these parameters, be sure to regenerate your $\BA$ matrix.) 
Produce plots to show how these parameters affect the blurring operation.

\medskip
\item \label{ME:ex:Atime}
The code in Figure~\ref{ME:fig:buildA} built the blurring matrix $\BA$
\emph{by row}.  What if instead we followed the more natural approach of 
building $\BA$ \emph{by entry}?  The code is given below.

\begin{lstlisting}[language=Python]
def build_blur_A_by_entry(n=100,z=.025):
    A = np.zeros((n,n));                         # n-by-n matrix of zeros
    s = np.array([(j+.5)/n for j in range(n)])   # s_j points
    t = np.array([(k+.5)/n for k in range(n)])   # t_k points
    for j in range(0,n):
        for k in range(0,n):
            A[j,k] = h_gaussian(s[j],t[k],z)/n   # A(j,k) = h(s_j,t_k)/n
    return A
\end{lstlisting}

This problem explores if there is an efficiency advantage to building $\BA$ by row.
You can time a command in Python by adjusting the following template.

\begin{lstlisting}[language=Python]
import time

t0 = time.time()                                 # start timer
# insert the commands you want to time here
elapsed_time = time.time() - t0                  # show elapsed time 
\end{lstlisting}

\begin{itemize}
\item[(a)] Compute the time it takes to build the $n=1000$ blurring matrix $\BA$ by row 
(use the code \verb|build_blur_A| from Figure~\ref{ME:fig:buildA}).
\item[(b)] Repeat the same timing (be sure to reset \verb|t0 = time.time()|), 
but now building $\BA$ entry-by-entry (use the code \verb|build_blur_A_by_entry|
given above).  
How does the time compare?  
\item[(c)] Repeat~(a) and~(b) with $n=1000$.  How do the timings change when
the dimension $n$ doubles from $500$ to $1000$?  
Can you explain why the timing grows like this?

\item[(d)] Modify \verb|build_blur_A| to build $\BA$ \emph{one column at a time}.
How does the timing of this routine compare to the timing for the the 
row-oriented and entry-by-entry methods? 
\end{itemize}

\vspace*{-1.5em}
\DIVIDER

\setcounter{section}{3} \setcounter{prob}{0}

\item  \label{ME:ex:bw1}
In problem~\ref{ME:ex:fw} you blurred a signal $\Bf$ to get $\Bb=\BA\Bf$,
where $\BA$ uses the Gaussian blurring function~(\ref{ME:gaussker}) with $n=100$
and $z=0.025$.  

Replicate the experiment shown in Figure~\ref{ME:fig:bw} using
this Gaussian kernel.  Following the example in the code below, add random noise
of level $10^{-3} \|\Bb\|$ to $\Bb$ to get $\bnoise$.
Then solve for $\frec = \BAs^{-1}\bnoise$ using the \verb|np.linalg.solve| command
in Python.  
(This approach uses Gaussian elimination to solve $A\frec = \bnoise$ 
for the unknown vector $\frec$, which is preferable to computing 
$\BAs^{-1}$ and multiplying the result against $\bnoise$ to get $\frec$.)

\begin{lstlisting}[language=Python]
n = 100
z = 0.025;
eps = 0.001;

t = np.array([(k+.5)/n for k in range(n)])
A = build_blur_A(n,z) 
f = f_function(t)                   # compute the true signal
b = A@f                             # blur the true signal
s = eps*np.linalg.norm(b)           # standard deviation of noise
e = s*np.random.randn(n)            # generate random noise
bnoise = b+e                        # simulate a noisy observation
frec = np.linalg.solve(A,bnoise)    # solve A*frec=bnoise for frec

plt.plot(t,f,'.',color='blue',markersize=5,label='$f$')
plt.plot(t,frec,'.',markerfacecolor='none',color='blue',\
         markersize=8,label='$f_{rec}$')
plt.legend()
\end{lstlisting}

\item Repeat Exercise~\ref{ME:ex:bw1}, but reduce the noise level down to 
      $10^{-5}\|\Bb\|$, $10^{-7}\|\Bb\|$, and $10^{-9}\|\Bb\|$.  
      Does $\frec$ improve as the level of noise drops?

\item \label{ME:ex:bign}
      Now take $n=500$, still with $z=0.025$.  
      This larger value of $n$ should mean that our 
      blurring operation incurs less discretization error when we 
      convert the calculus problem into a linear algebra problem.
      You might speculate this extra accuracy in the approximation of 
      the blurring integral would lead to a better recovered value
      for $\frec$.
      This question explores whether this conjecture is accurate.
      We want to determine how increasing $n$ affects the accuracy 
      of the recovered signal.

      Repeat Exercise~\ref{ME:ex:bw1} with $n=500$ and noise level 
      $10^{-3}\|\Bb\|$.  Is $\frec$ accurate?

      What if we attempt to compute $\frec$ directly from $\Bb$, with
      no artificial noise added (i.e., \verb|bnoise=b|).  
      Does \verb|frec = np.linalg.solve(A,b)| produce the correct $\Bf$?

\DIVIDER

\setcounter{section}{4} \setcounter{prob}{0}

\item \label{ME:ex:2x2}
Suppose we redefine $\BA$ so that it depends on a parameter $\alpha>0$:
\[ \BA = \left[\begin{array}{cc} 1 & -\alpha \\ 1 & \phantom{-}\alpha \end{array}\right].\]
\begin{enumerate}
\item[(a)] Compute a formula for $\BAs^{-1}$, which will depend upon $\alpha$.  

           How does $\BAs^{-1}$ change as $\alpha \to 0$?
\item[(b)] Using your formula from part~(a), write out $\BAs^{-1}$ for the 
           specific value $\alpha=0.001$.
\item[(c)] Again with $\alpha=0.001$, suppose 
           \[ \Bb_1 = \left[\begin{array}{c}1 \\ 1 \end{array}\right], \qquad
              \Bb_2 = \left[\begin{array}{c}1 \\ 1.01 \end{array}\right]. \]
           Solve for $\Bf_1 = \BAs^{-1}\Bb_1$ and $\Bf_2 = \BAs^{-1}\Bb_2$. 

           How does the difference between the solutions, $\|\Bf_2-\Bf_1\|$, 
           compare to the difference between the data, $\|\Bb_2 - \Bb_1\|$?
\end{enumerate}

\DIVIDER

\setcounter{section}{5} \setcounter{prob}{0}

\item \label{ME:ex:xy}
This exercise investigates the connection between the objective 
function~(\ref{ME:phi}) and the least squares problem~(\ref{ME:LS}).
\begin{enumerate}
\item[(a)] For any vectors $\Bx,\By\in\R^n$, show that 
\[ \left\|\left[\begin{array}{c} \Bx \\ \By \end{array}\right]\right\|^2 
     = \|\Bx\|^2 + \|\By\|^2.\]
\item[(b)] Now show that
\[ \left\| \left[\begin{array}{c} \Bb \\ \Bzero \end{array}\right]  
           - \left[\begin{array}{c} \BA \\ \lambda@\BI \end{array}\right]@\Bf\,\right\|^2  
    = \|\Bb-\BA\Bf\|^2 + \lambda^2@\|\Bf\|^2,\]
    and explain why $\phi(\Bf)$ in~(\ref{ME:phi}) is thus minimized by the same $\Bf$
    that solves the least squares problem~(\ref{ME:LS}).
\end{enumerate}
\item \label{ME:ex:reg2x2}
Suppose that 
\[ \BA = \left[\begin{array}{ll} 1 & 1 \\ 1+10^{-6} & 1-10^{-6} \end{array}\right], \qquad
   \Bb = \left[\begin{array}{l} 1 \\ 1 \end{array}\right].\]
Set $\lambda=10^{-8}$. 
Solve for $\Bf_\lambda$ using (\ref{ME:phi}) using Gaussian elimination
(e.g., via the \verb|np.linalg.solve| command),
and by solving the least squares problem (\ref{ME:LS}) 
(e.g., via the \verb|np.linalg.lstsq| command).
Show your results to at least 8~decimal places.  Do your answers agree?
If not, which one do you think is more accurate?  Why?

\DIVIDER

\setcounter{section}{6} \setcounter{prob}{0}

\item \label{ME:ex:Lcurve}
Using the code above as a starting point, generate an L curve like the
one shown in Figure~\ref{ME:fig:bw_L} by taking the following steps.  
\begin{enumerate}
\item[(a)] Reproduce Figure~\ref{ME:fig:bw_L} using the same parameters:
take $n=500$, generate $\BA$ to use the hat-function blurring kernel~\ref{ME:hatker}
with $z=0.05$,
and set the noise level to be $10^{-5}\|\Bb\|$. 
\item[(b)] Compute the regularized solutions for several special values of $\lambda$;
on top of your L~curve plot from part~(a), 
plot $(\|\bnoise-\BA\flam\|, \|\flam\|)$ as red circles.
\item[(c)] Adapt your code to use the Gaussian kernel~\ref{ME:gaussker}
with $z=0.01$ and noise level $10^{-5}\|\Bb\|$.
Adjust the blurring parameter $z$ and/or the noise level, 
and see how these changes affect the L~curve plot.

Note: as you adjust the parameters in your model, you might need to change
the range of $\lambda$ values you are using (defined in \verb|lamvec|);
it is also a good idea to remove the lines that set the axis limits
(\verb|plt.xlim| and \verb|plt.ylim|), as the location of the L~curve
will change as you adjust the problem.
\end{enumerate}

\DIVIDER

\setcounter{section}{7} \setcounter{prob}{0}

\item \label{ME:ex:geupc}
Using the unpolluted blurred vector $\Bb$,
produce a plot showing the vector $\Bf_{\rm rec}$ one obtains by directly 
solving $\BA\Bf_{\rm rec} = \Bb$
 (e.g., \verb|frec = numpy.linalg.solve(A,b)|) via Gaussian elimination.  
To produce an elongated plot like the ones shown above in Python,  
you can adapt the following code.
\begin{lstlisting}[language=Python]
plt.figure(figsize=(10,1.85))      % create an elongated figure window
plt.plot(t,f,'b.')                 % plot f_true
plt.xlabel('$s_j$',fontsize=14)    % label horizontal axis s_j
plt.ylabel('$f_j$',fontsize=14)    % label vertical axis f_j
plt.tight_layout()
\end{lstlisting}

\medskip
\item Now using the polluted vector $\bnoise$, repeat the last exercise:  
use Gaussian elimination (\verb|numpy.linalg.solve|) to
solve $\BA\Bf_{\rm noise} = \Bb_{\rm noise}$ for $\Bf_{\rm noise}$.
Plot your recovered barcode $\Bf_{\rm noise}$.  
Does your solution at all resemble $\Bf_{\rm true}$?

\medskip
\item  \label{ME:ex:coke_reg}
Explore the solutions obtained from the regularized least squares problem
\[ \min_{\Bf \in \R^n} \|\BA\Bf-\Bb\|^2 + \lambda^2 \|\Bf\|^2\]
for various values of $\lambda>0$.  
Recall from Section~\ref{ME:sec:reg} that you can find $\Bf_\lambda$ 
by solving the least squares problem
\[  \min_{\Bf \in \R^n} \|\wh{\Bb} - \BA_\lambda\Bf\|,\]
where 
\[ \BA_{\lambda} = \left[\begin{array}{c}\BA \\ \lambda\BI\end{array}\right] \in \R^{2n\times n}, \qquad
   \widehat{\Bb} = \left[\begin{array}{c} \Bb \\ \Bzero\end{array}\right]\in\R^{2n}.\]
This least squares problem can be solved in Python via the \verb|np.linalg.lstsq| command.

\begin{enumerate}
\item[(a)] Using the value $\Bb_{\rm noise}$ (\textsl{be sure to use $\Bb_{\rm noise}$}), 
create an L curve plot using parameter values $\lambda = 10^{-6},\ldots, 10^{1}$.
(Use \verb|lam = np.logspace(-6,1,100)| to generate 100 logarithmically-spaced values of $\lambda$
in this range.)  Recall that the L curve is a \verb|plt.loglog| plot with 
$\|\Bb_{\rm noise} - \BA \Bf_\lambda\|$ on the horizontal axis
and $\|\Bf_\lambda\|$ on the vertical axis.
(\textsl{Be sure you are using $\Bb_{\rm noise}$ to compute the norm of the
misfit for the horizontal axis, $\|\Bb_{\rm noise} - \BA \Bf_\lambda\|$.})

\item[(b)] Pick several $\lambda$ values corresponding to the main corner in the L-curve
(or just a bit beyond this point), solve for the regularized solution $\flam$, 
and plot $\flam$.  
Show plots for a few different $\lambda$ values that vary over
several orders of magnitude.  Also show the $\Bf_\lambda$ you get for the
same value of $\lambda$ if you instead use the true data $\Bb$ instead of 
$\Bb_{\rm noise}$.
\end{enumerate}

\smallskip
\item Using your best recovered $\Bf_\lambda$ from Exercise~\ref{ME:ex:coke_reg},
reconstruct the barcode for the can of Coke.  (You should \emph{decode} the specific 
numbers shown at the bottom of the Coke barcode presented earlier; 
it is not enough to just recover the function $f$ itself.)  This is a little tricky
-- do your best, using insight from the structure of the UPC code function to help.

\medskip
\item Edit the \verb|coke_upc| function to create \verb|bnoise| with more noise,
      such as $5\times 10^{-3} \|\Bb\|$ (instead of the $10^{-3}\|\Bb\|$ in the original 
      \verb|coke_upc| code).  Can you still reliably compute the barcode?
      If so, how much more noise can you add, before it becomes very difficult to 
      recover the correct bar code (even with the use of regularization)?

\medskip
\item \label{ME:ex:mystery1}
The file \verb|mystery_bnoise1.csv| specifies the 
blurred, noisy $\Bb_{\rm noise}$ barcode for a mystery product 
(again with {\tt n=570} and the same $\BA$ as for the Coke barcode).

You can load the file with the simple command
\begin{lstlisting}[language=Python]
bnoise = np.loadtxt('mystery_bnoise1.csv')
\end{lstlisting}

\begin{enumerate}
\item[(a)] Use the techniques from Exercise~\ref{ME:ex:coke_reg} to attempt 
           to recover the function $f$ associated with the mystery barcode.
           Show samples of the results from your various attempts.
\item[(b)] Decode the UPC from your best recovered solution $\flam$ from part~(a).
           (Be sure to explain your process.)
           What is the product described by the mystery UPC?  
           (Once you have found the correct UPC, a search on Amazon 
            or \url{https://www.barcodelookup.com} should turn up the product.)
\end{enumerate}

\smallskip
\item  \label{ME:ex:mystery2}
Repeat the deblurring process for the blurred barcode 
stored in \verb|mystery_bnoise2.csv|, which is a bit noisier that the previous mystery 
UPC (but again uses {\tt n=570} and the same $\BA$ as for the two previous barcodes).

To what product does this blurry barcode correspond?  Include your evidence.

\DIVIDER

\setcounter{section}{8} \setcounter{prob}{0}

\item \label{ME:ex:comp_svs}
      Compute the singular values of the blurring matrix~(\ref{ME:blurmat})
      for the Gaussian kernel~(\ref{ME:gaussker}) with $n=570$ and $z=0.01$, 
      as used in the barcode decoding experiments in the exercises at the 
      end of section~\ref{ME:sec:UPC}. 
      How does the decay of the singular values for this kernel compare with what
      we observed for the hat function kernel in Figure~\ref{ME:fig:blur2_svd}?
      What insight does your plot give you about the sensitivity of the UPC
      deblurring operation?

\medskip
\item Repeat the experiments in this section for the Coke UPC decoding 
example described earlier, blurring with the Gaussian kernel as in the
last problem.  Compare the results you obtain using the unpolluted 
vector $\Bb$, as well as $\bnoise$ with different levels of noise.

\medskip
\item \label{ME:ex:trunc_svd}
The discussion in this section might cause you to ask:  
if small singular values cause problems in the formula~(\ref{ME:frec_sv}) 
for $\frec$, why not just chop those troublesome terms off of the sum, 
stopping after $k<n$ terms:
\begin{equation} \label{ME:frec_tsvd}
   \Bf_k = \sum_{j=1}^k {\Bu_j^T\bnoise \over \sigma_j}@\Bv_j.
\end{equation}
This \emph{truncated SVD solution} can be a very effective alternative
to the regularization method presented in this manuscript.

\medskip
Repeat the Coke UPC decoding operation using the truncated SVD solution $\Bf_k$.
Can you correctly decode the UPC from $\Bf_k$?
Try a range of $k$ values.  Which $k$ values give smoother or rougher solutions?
\end{exercises}


\section*{Modeling Projects}

\medskip
\begin{enumerate}
\item[M1.] Find a UPC barcode for some product of your choice.  
Make your own version of Exercises~\ref{ME:ex:mystery1} and~\ref{ME:ex:mystery2}.  
Turn your barcode in an $\Bf$ vector.  
Apply the blur, $\Bb = \BA\Bf$, add some noise to get $\bnoise$, and
then give the result to a friend to decode.

\medskip
\item[M2.] All our barcode deblurring work has been based on taking $n=570$,
which means that each of the 95~units that make up a barcode is discretized
with 6~entries in the vector $\Bf\in\R^{570}$ (notice that $95\times 6 = 570$).

What if you rework everything using only 2~entries per unit 
(giving $\Bb, \Bf\in\R^{190}$ and $\BA\in\R^{190\times 190}$).
Reducing the size of the matrices and vectors will make the 
linear algebra \emph{faster}.  
Does it make the barcode deblurring easier or harder?  Show evidence.

What if you rework everything using 10~entries per unit 
(giving $\Bb, \Bf\in\R^{950}$ and $\BA\in\R^{950\times 950}$).
Enlarging the size of the matrices and vectors will make the 
linear algebra \emph{slower}.  
Does it make the barcode deblurring easier or harder?  Show evidence.

\medskip
\item[M3.] The previous modeling project asks you to experiment with 
blurring matrices of different dimension.  Increasing $n$ leads to a
more high-fidelity model, with a more accurate approximation to the
blurring integral~(\ref{ME:midpt}).  

Construct the blurring matrix $\BA$ in~(\ref{ME:blurmat}),
using the hat function kernel~(\ref{ME:hatker}) with $z=0.05$,
but now with $n=1000$.  Compute the singular values, and compare
them to those shown for $n=500$ in Figure~\ref{ME:fig:blur2_svd}.
What happens to the smallest singular values?  

Repeat the experiments of Section~\ref{ME:sec:SVD} with $n=1000$.
Is the recovered solution $\frec = \BA^{-1}\bnoise$ more or less
susceptible to noise in $\bnoise$, compared to the $n=500$ case?
Use regularization to compute $\flam$ for several $\lambda$ values.
Does regularization work?

Now repeat these experiments with $n=2000$ (or larger).  
Can you deduce a trend in performance as $n$ increases?

\medskip
\item[M4.] The ``images'' we have deblurred in this manuscript have all been
one-dimensional.  Such primitive ``images'' make the modeling a 
bit clearer to explain, but are perhaps less satisfying to study
than true two-dimensional images.  Two dimensions requires two independent
variables, and so the true image becomes $f(t^{(1)},t^{(2)})$, 
with its blurry counterpart $b(s^{(1)},s^{(2)})$.  
The analog of the blurring operation~(\ref{ME:conv}) takes the form
\begin{equation} \label{ME:eq:2dblur}
 b(s^{(1)},s^{(2)}) 
   = \int_0^1 \int_0^1 
    h\left(\left[\!\begin{array}{c} s^{(1)} \\ s^{(2)}\end{array}\!\right],
          \left[\!\begin{array}{c} t^{(1)} \\ t^{(2)}\end{array}\!\right]\right) 
     f(t^{(1)},t^{(2)})\,\dop t^{(1)} \dop t^{(2)}.
\end{equation}

\begin{enumerate}
\item Generalize one of the blurring functions~(\ref{ME:avgker}), 
      (\ref{ME:hatker}), and (\ref{ME:gaussker}) to obtain a
      suitable $h$ for this two dimensional setting.

\item Discretize $s^{(1)}, s^{(2)} \in[0,1]$ 
      into $n$ points $s_\ell = (\ell-1/2)/n$ for $\ell=1,\ldots, n$, 
      and approximate the double integral~(\ref{ME:eq:2dblur})
      with a double sum:
\begin{equation} \label{ME:eq:2dsum}
      b(s_\ell,s_m) 
      \approx  
      b_{\ell,m} = \frac{1}{n^2} \sum_{j=1}^n \sum_{k=1}^n 
       h\left(\left[\!\!\begin{array}{c} s_\ell \\ s_m\end{array}\!\!\right],
          \left[\!\begin{array}{c} t_j \\ t_k \end{array}\!\right]\right),
\end{equation}
      where $t_j = (j-1/2)/n$.
 
      Explain how to arrange the equation~(\ref{ME:eq:2dsum}) for 
      $\ell,m=1,\ldots, n$ into the form of a linear system 
      $\Bb = \BA\Bf$, where $\BA\in\R^{n^2\times n^2}$ and 
\begin{eqnarray*}
     \Bb &=& [\begin{array}{ccccccc}
      b_{1,1} & \cdots & b_{n,1} & \cdots & b_{1,n} & \cdots & b_{n,n} 
       \end{array}]^T \in \R^{n^2\times 1}  \\[3pt]
     \Bf &=& [\begin{array}{ccccccc}
      f(t_1,t_1) & \cdots & f(t_n,t_1) & \cdots & f(t_1,t_n) & \cdots & f(t_n,t_n) 
       \end{array}]^T \in \R^{n^2\times 1}.
\end{eqnarray*}
      To make things concrete, write everything out for the case $n=3$.

\item Find a grayscale image, and use an image editing program to adjust
      its size to $n\times n$ pixels, for some modest value of $n$.
      Take $f(t_j,t_k)$ to be the value of the $(j,k)$ pixel, so that, 
      for a grayscale image, $f(t_j,t_k)$ will be an integer between
      0 (black) and 255 (white).  Plot the image as a heat map.

\item Code up the matrix $\BA$ in part~(b) for a general value of $n$,
      Take your $f$ from part~(c), and reshape it into the vector 
      $\Bf\in\R^{n^2\times 1}$.
      Blur the image:  $\Bb = \BA\Bf$.  Reshape $\Bb$ into an $n\times n$ matrix,
      and visualize the result as a heat map.

\item Attempt to deblur the image by solving $\frec = \BA^{-1}\Bb$.
      How does the process perform?  What if noise is added to $\Bb$?
      Explore the use of regularization to stabilize the deblurring process.
\end{enumerate}
Note:  When working with two dimensional images, the linear algebra 
      becomes quite a bit more expensive.
      (For example, if $n=1000$, then
      $n^2 = 10^6$, and your $\BA$ would be a million-by-a million matrix.)
      For these experiments, try using $n=40$ or $n=60$.

For more details about the two-dimensional blurring process in the
spirit of this manuscript (and hints for some of the above problems), 
see~\cite[Section~8.3]{Emb3606}.
\end{enumerate}

\section*{Acknowledgements}
This manuscript was written as a submission for the volume
\emph{Build Your Course in Mathematical Modeling: A Modular Textbook},
developed by Dr.~Shelley Rohde Poole in coordination with the 
SIAM Activity Group on Applied Mathematics Education.
I am grateful for many helpful suggestions from Shelley Rohde Poole,
and from Mahesh Banavar, Seo-Eun Choi, and Gabriel Soto, which have
improved this manuscript.


\end{document}